# An improved parametric model for hysteresis loop approximation

Rostislav V. Lapshin[1, 2]

[1]*Solid Nanotechnology Laboratory, Institute of Physical Problems, Zelenograd, Moscow, 124460, Russian Federation*

[2]*Department of Integral Electronics and Microsystems, Moscow Institute of Electronic Technology, Zelenograd, Moscow, 124498, Russian Federation*

E-mail: rlapshin@gmail.com

A number of improvements have been added to the existing analytical model of hysteresis loop defined in parametric form. In particular, three phase shifts are included in the model, which permits to tilt the hysteresis loop smoothly by the required angle at the split point as well as to smoothly change the curvature of the loop. As a result, the error of approximation of a hysteresis loop by the improved model does not exceed 1%, which is several times less than the error of the existing model. The improved model is capable of approximating most of the known types of rate-independent symmetrical hysteresis loops encountered in the practice of physical measurements. The model allows building smooth, piecewise-linear, hybrid, minor, mirror-reflected, inverse, reverse, double and triple loops. One of the possible applications of the model developed is linearization of a probe microscope piezoscanner. The improved model can be found useful for the tasks of simulation of scientific instruments that contain hysteresis elements.

## I. INTRODUCTION

The phenomenon of hysteresis is widespread in nature, it is often met in many fields of science and engineering including instruments used in scientific research.[1, 2, 3] There are a number of quite complicated analytical models describing this phenomenon.[4] One of the simple ones is the analytical model suggested in Ref. 1. With that model, a family of hysteresis loops is described by the following parametric equations

$$\begin{aligned} x(\alpha) &= a\cos^m\alpha + b_x\sin^n\alpha, \\ y(\alpha) &= b_y\sin\alpha, \end{aligned} \qquad (1)$$

where $\alpha$ is a real parameter ($\alpha$=0…2π); $a$ is $x$ coordinate of the split point (see Fig. 1); $b_x$, $b_y$ are the saturation point coordinates; $m$ is a positive odd integer ($m$=1, 3, 5, …) defining the curvature of the hysteresis loop; $n$ is a positive integer defining the type of the hysteresis loop and its curvature. With $n$=1, the Leaf loop type is formed; with $n$=2, 4, 6, … – the Crescent (Boomerang), and with $n$=3, 5, 7, … – the Classical. With increasing the parameter $\alpha$, the movement along the loop occurs in the counterclockwise direction, with decreasing – clockwise. The start point ($\alpha$=0) and the end point ($\alpha$=2π) of a loop are both at the split point $a$.

The key distinctive feature of the model (1) is its simplicity. The model is intuitive, it allows quickly creating hysteresis loops of a required type and easily determining their key parameters $a$, $b_x$, $b_y$, $m$, $n$. The 1.5-6%[1] approximation accuracy of the model (1) is quite enough for most practical tasks. However, there are cases when a higher accuracy is required. The improved model[5] suggested in the article approximates hysteresis loops with an error 1% or less.

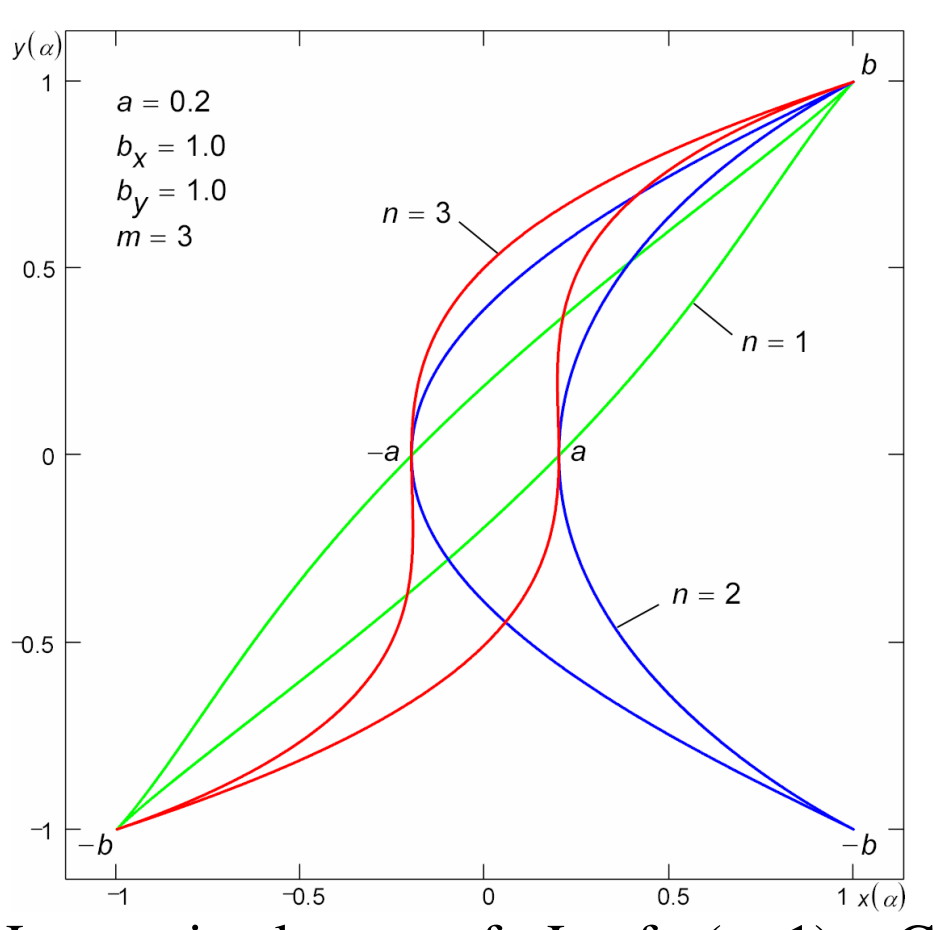

Fig. 1. Hysteresis loops of Leaf ($n$=1), Crescent (Boomerang, $n$=2), and Classical ($n$=3) types. The area of all the three loops is the same.

The model (1) covers most of the known types of rate-independent symmetrical smooth hysteresis loops. The improved model allows controlling the tilt and curvature of smooth loops more accurately (see Sec. II.A.2). Besides the smooth loops, the improved model allows building various piecewise-linear loops (see Sec. II.B) as well as hybrid loops (see Sec. II.B.1.b) in which rectilinear sections are combined with curvilinear ones. Moreover, the use of the improved model can help create continuously drawn double (see Sec. II.C) and triple (see Sec. II.D) loops (both self-crossing and non-self-crossing) out of smooth, piecewise-linear, and hybrid loops as well as of their combinations. Section II.E provides formulas for area calculation of the hysteresis loops. Overall, the article presents a general approach applicable for approximation of a large number of varieties of hysteresis loops.

To simplify the construction, analysis, and identification of hysteresis loops under consideration, supplementary material is provided in the form of Mathcad® worksheets (MathSoft, USA).[6] In order not to overcomplicate the article, the derivation of some formulas is omitted. The detailed derivation can be found in the supplementary material.

## II. DESCRIPTION OF THE IMPROVED MODEL

### A. Smooth hysteresis loops

#### *1. Additional representations of hysteresis loop*

#### *a. Representation in the form of the sum of an unsplit loop and a splitting curve*

The hysteresis loop (1) can always be represented as a sum of two parametric curves

$$\begin{aligned} x(\alpha) &= x_1(\alpha) + x_2(\alpha), \\ y(\alpha) &= y_1(\alpha) + y_2(\alpha), \end{aligned} \tag{2}$$

where $x_1(\alpha)=b_x \sin^n\alpha$, $y_1(\alpha)=b_y \sin\alpha$ is the unsplit loop; $x_2(\alpha)=a\cos^m\alpha$, $y_2(\alpha)=0$ is the splitting curve. This representation is useful while considering transformations that change tilt and/or curvature of the hysteresis loop (1): first, the unsplit loop is subjected to tilting/curving, after that the obtained result is split by simply adding the splitting curve.[6] Moreover, by altering/adding a direction of action of the splitting curve, it is possible to create double (see Sec. II.C.3) and triple (see Sec. II.D.3) hysteresis loops.

#### *b. Representation in the form of a frequency spectrum*

Using the de Moivre's formula, the generating function $x(\alpha)$ in model (1) can also be represented as a sum of cosines and sines having multiple frequencies (frequency spectrum)

$$\begin{aligned} x(\alpha) &= \frac{a}{2^{m-1}} \sum_{k=0}^{\frac{m-1}{2}} C_m^k \cos((m-2k)\alpha) + \frac{b_x}{2^{n-1}} \sum_{k=0}^{\frac{n-1}{2}} (-1)^{\frac{n-1}{2}+k} C_n^k \sin((n-2k)\alpha), \\ y(\alpha) &= b_y \sin\alpha, \end{aligned} \tag{3}$$

where $C_l^k$ is a binomial coefficient ($k$, $l$ are positive integers); $C_l^k = l!/[k!(l-k)!]$ if $0 \le k \le l$, otherwise $C_l^k = 0$ (below in (6), (25), $C_l^k = 0$ if $k$ is a real number). The equations (3) are valid for odd $n$; equations for even $n$ are given in the supplementary material. For example, equations describing the hysteresis loop Classical (see Fig. 1) in accordance with (3) with $m$=$n$=3 are as follows

$$
\begin{aligned}
x(\alpha) &= \frac{a}{4}\left[3\cos\alpha + \cos(3\alpha)\right] + \frac{b_x}{4}\left[3\sin\alpha - \sin(3\alpha)\right], \\
y(\alpha) &= b_y \sin\alpha.
\end{aligned} \tag{4}
$$

It is easy to see that notation of the generating function $x(\alpha)$ in the form (3) is, in fact, a decomposition of $x(\alpha)$ into odd harmonics of the Fourier series

$$
\begin{aligned}
x(\alpha) &= \frac{1}{2}A_0 + \sum_{k=1}^{l}\left(A_k \cos(k\alpha) + B_k \sin(k\alpha)\right), \\
y(\alpha) &= b_y \sin\alpha,
\end{aligned} \tag{5}
$$

where the Fourier coefficients $A_k$, $B_k$ are determined by the algebraic formulas

$$
\begin{aligned}
A_k &= \frac{a}{2^{m-1}} C_m^{\frac{m-k}{2}}, \\
B_k &= (-1)^{\left\lfloor \frac{k-1}{2} \right\rfloor} \frac{b_x}{2^{n-1}} C_n^{\frac{n-k}{2}}.
\end{aligned} \tag{6}
$$

The amplitude $A_0$ of the zero-frequency component ($k$=0) and the amplitudes $A_k$ and $B_k$ of all even harmonics ($k$=2, 4, 6, …) in (5) are equal to zero. The value of $l$ is set to the largest of the $m$ and $n$ power. The floor function in the expression for $B_k$ in (6) is optional and used only to avoid the generation of complex numbers when $k$ is even. For example, the Fourier coefficients of the generating function $x(\alpha)$ of the Classical loop ($m$=$n$=3) shown in Fig. 1 are as follows: $A_k$=(0, 3$a$/4, 0, $a$/4), $B_k$=(0, 3$b_x$/4, 0, -$b_x$/4).

Having the Fourier coefficients $A_k$, $B_k$ (6), the generating function $x(\alpha)$ can also be represented as

$$
x(\alpha) = \sum_{k=1}^{l} Am_k \cos(k\alpha - \varphi_k), \tag{7}
$$

where amplitudes $Am_k$ and phases $\varphi_k$ of the harmonics are determined by the formulas

$$
\begin{aligned}
Am_k &= \sqrt{A_k^2 + B_k^2} = \sqrt{\left(\frac{a}{2^{m-1}} C_m^{\frac{m-k}{2}}\right)^2 + \left(\frac{b_x}{2^{n-1}} C_n^{\frac{n-k}{2}}\right)^2}, \\
\tan\varphi_k &= \frac{B_k}{A_k} = (-1)^{\left\lfloor \frac{k-1}{2} \right\rfloor} 2^{m-n} \frac{C_n^{\frac{n-k}{2}} b_x}{C_m^{\frac{m-k}{2}} a}.
\end{aligned} \tag{8}
$$

For example, the hysteresis loop Classical ($m$=$n$=3) shown in Fig. 1 has the amplitudes $Am_k$=(0, $3\sqrt{a^2+b_x^2}/4$, 0, $\sqrt{a^2+b_x^2}/4$) and the phases $\varphi_k$=(0, $b_x/a$, 0, -$b_x/a$).

The generating function $x(\alpha)$ can also be represented in the exponential form

$$
x(\alpha) = \sum_{k=-l}^{l} C_k e^{ik\alpha}, \tag{9}
$$

where $i$ is the imaginary unit, and the complex Fourier coefficient $C_k$ is determined as follows

$$C_k = \begin{cases} \frac{1}{2}(A_k - iB_k),\ k = 1,3,5,..., \\ \frac{1}{2}(A_k + iB_k),\ k = -1,-3,-5,..., \\ 0,\ k = 0,\pm 2,\pm 4,\pm 6,.... \end{cases} \tag{10}$$

Having the complex Fourier coefficient $C_k$, the generating function $x(\alpha)$ can also be represented as

$$x(\alpha) = \sum_{k=-l}^{l} Am_k e^{i(k\alpha - \varphi_k)} = \sum_{k=-l}^{l} Am_k \cos(k\alpha - \varphi_k), \tag{11}$$

where amplitudes $Am_k$ and phases $\varphi_k$ of the harmonics are determined by the formulas

$$\begin{aligned} Am_k &= \sqrt{\mathrm{Re}^2(C_k) + \mathrm{Im}^2(C_k)}, \\ \varphi_k &= -\mathrm{Arg}(C_k) = -\arctan\frac{\mathrm{Im}(C_k)}{\mathrm{Re}(C_k)}, \end{aligned} \tag{12}$$

where Re, Im, and Arg are the real part, the imaginary part, and the argument of a complex number, respectively.

Representation of the generating function $x(\alpha)$ in the form of a frequency spectrum (3), (5), (7) and (11) allows to synthesize hysteresis loops of a custom shape, tilt and curvature by changing the amplitude and phase of harmonic components and by adding/excluding harmonic components with certain values of amplitude and phase.[6] This approach permits to build the smooth hysteresis loops with almost any shape.

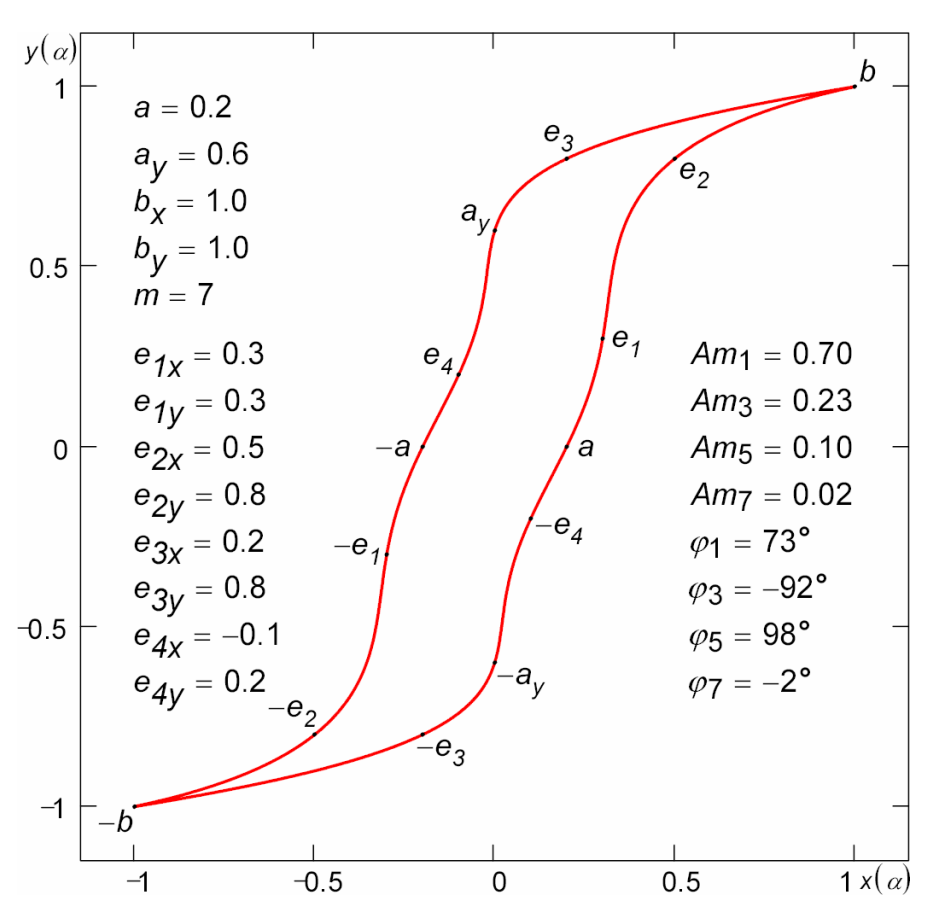

Fig. 2. Classical hysteresis loop of a complex shape, whose generating function $x(\alpha)$ is formed by summing the first four odd harmonics. The loop is drawn exactly through 14 predefined points. The loop area depends on the amplitude and phase of the 1st harmonic only.

Fig. 2 shows an example of a Classical hysteresis loop synthesis with a complex shape[7] using harmonic components. The generating function $x(\alpha)$ of the loop is formed by summing the first four odd harmonics ($m$=7). The hysteresis loop is drawn exactly through 14 predefined points, six of them are the key points of the loop $\pm a$, $\pm a_y$, $\pm b$ (where $a_y$ is a vertical splitting, see Secs. II.C.3, II.D.3). The amplitudes and phases of the harmonic components are determined numerically by solving a system of nonlinear equations.[6] The loop area depends on the amplitude and phase of the 1st harmonic only; the rest of harmonics do not affect the loop area (see Sec. II.E.1).

### *2. Using phase shifts*

Among the main modifications to the previously proposed model of the hysteresis loop is the introduction of phase shifts $\Delta\alpha_1$, $\Delta\alpha_2$, $\Delta\alpha_3$

$$\begin{aligned} x(\alpha) &= \hat{a}\cos^m(\alpha + \Delta\alpha_1) + \hat{b}_x \sin^n(\alpha + \Delta\alpha_2), \\ y(\alpha) &= b_y \sin(\alpha + \Delta\alpha_3), \end{aligned} \tag{13}$$

where $\hat{a}$, $\hat{b}_x$ are corrected parameters of $a$, $b_x$, respectively. First, let us consider the influence of each of the three phase shifts $\Delta\alpha_1$, $\Delta\alpha_2$, $\Delta\alpha_3$ separately.

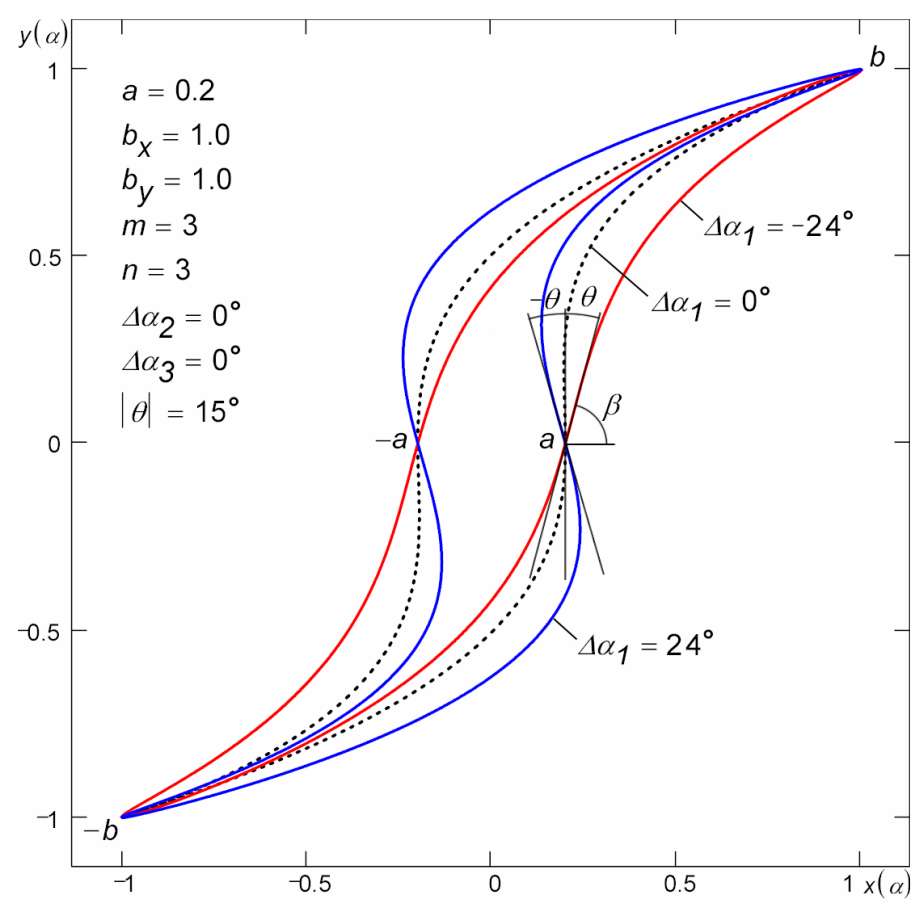

Fig. 3. Tilting hysteresis loop of the Classical type with the phase shift $\Delta\alpha_1$. Tilting by angle $\pm\theta$ at the split point $a$ is provided with a shift by $\mp\Delta\alpha_1$. Loop tilting results in an increase in the loop area. The area of the loop with a positive slope is equal to the area of the loop with a negative slope.

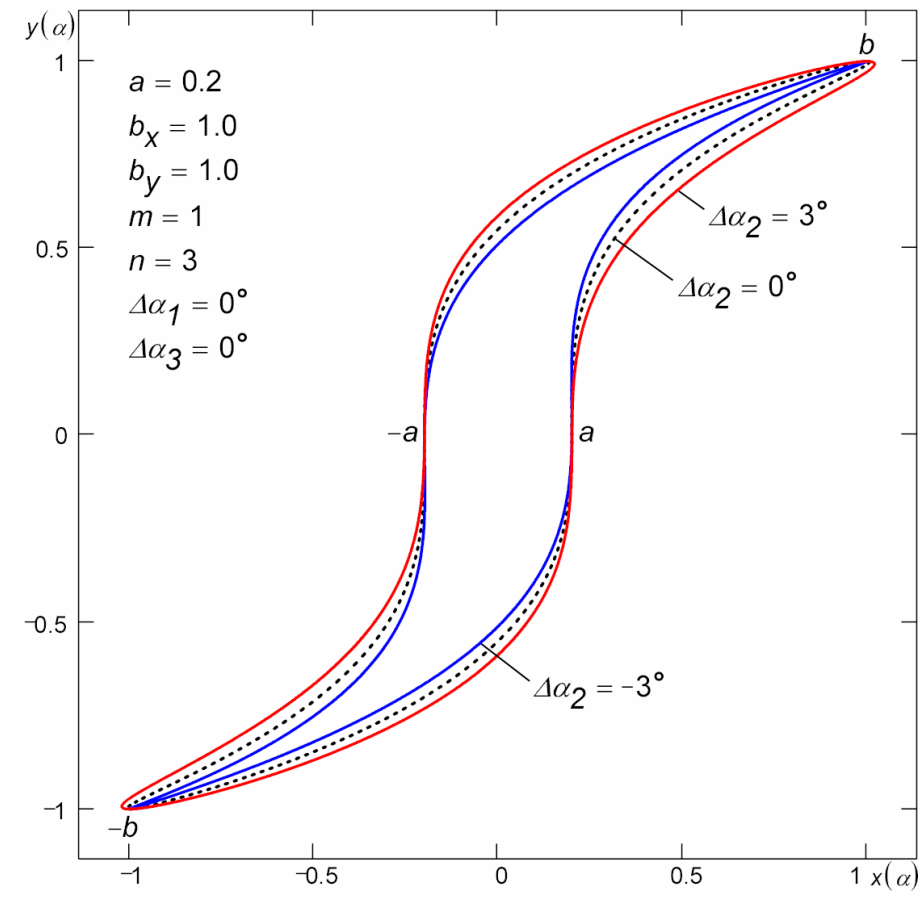

Fig. 4. Continuous change in the curvature of the hysteresis loop by the phase shift $\Delta\alpha_2$. The loop area increases with an increase in the phase shift $\Delta\alpha_2$.

***a. Loop tilting by the phase shift $\Delta\alpha_1$***

The phase shift $\Delta\alpha_1$ allows to tilt the hysteresis loop gradually by changing slope angle $\beta=\pi/2-\theta$ of a tangent to the loop at the split point $a$ (see Fig. 3). In the model suggested earlier, loop tilting was available by rotation of coordinate system by angle $\theta$ and predistortion of the loop parameters $a$, $b_x$, $b_y$ by their rotation in the opposite direction.[1] Below are more correct formulas for loop tilting by rotation, which exclude a slight displacement of the split point from the given position that existed previously ($\Delta\alpha_1=\Delta\alpha_2=\Delta\alpha_3=0$):[6]

$$\begin{aligned}\bar{x}(\alpha) &= x(\alpha)+\sin\theta\left(b_x\sin\theta+b_y\cos\theta\right)\left(\sin\alpha-\sin^n\alpha\right),\\ \bar{y}(\alpha) &= y(\alpha)+\sin\theta\left(b_x\cos\theta-b_y\sin\theta\right)\left(\sin\alpha-\sin^n\alpha\right).\end{aligned} \tag{14}$$

It can be seen from transformations (14) that the tilted loop represents by itself a result of addition of the initial loop with some curve, which provides the tilt of the initial loop (see Sec. II.A.4).

Since the introduction of the phase shift $\Delta\alpha_1$ leads to a change of the coordinates of the points of splitting $a$ and saturation $b_x$, a correction of coordinates of these points is required. In the improved model (13), the corrected parameters $\hat{a}$, $\hat{b}_x$ are found from the following simple system of equations composed for the split point $\alpha=0$ and the saturation point $\alpha=\pi/2$ ($\Delta\alpha_2=\Delta\alpha_3=0$)

$$\begin{aligned}&\hat{a}\cos^m(0+\Delta\alpha_1)+\hat{b}_x\sin^n 0=a,\\ &\hat{a}\cos^m\left(\frac{\pi}{2}+\Delta\alpha_1\right)+\hat{b}_x\sin^n\frac{\pi}{2}=b_x,\end{aligned} \tag{15}$$

whence

$$\begin{aligned}&\hat{a}=\frac{a}{\cos^m\Delta\alpha_1},\\ &\hat{b}_x=b_x+a\tan^m\Delta\alpha_1\end{aligned} \tag{16}$$

can be easily determined. Leaf type loops ($n$=1) with $m$=1 do not depend on the phase shift $\Delta\alpha_1$. The phase shift $\Delta\alpha_1$ required for loop tilting by the preset angle $\theta$ at the split point $a$ ($\alpha$=0) is calculated by the formula[6]

$$\Delta\alpha_1=-\arctan\frac{b_y\tan\theta}{ma}. \tag{17}$$

Unsplit loops ($a$=0) cannot be tilted with the phase shift $\Delta\alpha_1$.

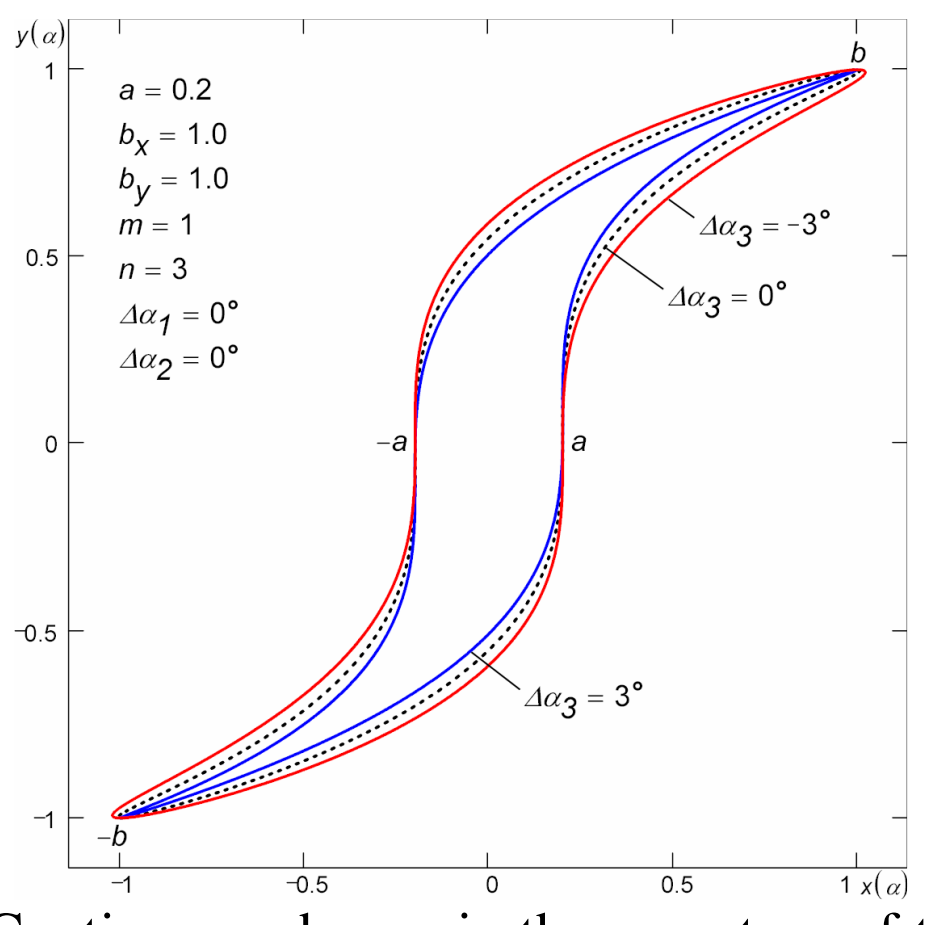

Fig. 5. Continuous change in the curvature of the hysteresis loop by the phase shift $\Delta\alpha_3$. The loop area decreases with an increase in the phase shift $\Delta\alpha_3$.

***b. Changing loop curvature by the phase shift*** $\Delta\alpha_2$

The phase shift $\Delta\alpha_2$ allows for changing the curvature of a hysteresis loop (see Fig. 4). Unlike the parameter $m$, the phase shift $\Delta\alpha_2$ provides a continuous change of the curvature of the loop. Like the phase shift $\Delta\alpha_1$, in case of the phase shift $\Delta\alpha_2$, corrected coordinates $\hat{a}$, $\hat{b}_x$ of the split and the saturation points, respectively, should be determined. To do that, a system of equations composed for the split point $\alpha$=0 and the saturation point $\alpha$=π/2 should be solved ($\Delta\alpha_1$=$\Delta\alpha_3$=0)

$$\begin{aligned} &\hat{a}\cos^m 0+\hat{b}_x\sin^n\left(0+\Delta\alpha_2\right)=a,\\ &\hat{a}\cos^m\frac{\pi}{2}+\hat{b}_x\sin^n\left(\frac{\pi}{2}+\Delta\alpha_2\right)=b_x, \end{aligned} \tag{18}$$

whence

$$\begin{aligned} &\hat{a}=a-b_x\tan^n\Delta\alpha_2,\\ &\hat{b}_x=\frac{b_x}{\cos^n\Delta\alpha_2}. \end{aligned} \tag{19}$$

Leaf type loops ($n$=1) with $m$=1 do not depend on the phase shift $\Delta\alpha_2$.

***c. Changing loop curvature by the phase shift*** $\Delta\alpha_3$

Like the phase shift $\Delta\alpha_2$, the phase shift $\Delta\alpha_3$ allows for continuous change of the curvature of a hysteresis loop (see Fig. 5). Composing equations for the split point $\alpha$=0 and the saturation point $\alpha$=π/2, like it has been done above, one can obtain ($\Delta\alpha_1$=$\Delta\alpha_2$=0)

$$\begin{aligned} &\hat{a}\cos^m\left(0-\Delta\alpha_3\right)+\hat{b}_x\sin^n\left(0-\Delta\alpha_3\right)=a,\\ &\hat{a}\cos^m\left(\frac{\pi}{2}-\Delta\alpha_3\right)+\hat{b}_x\sin^n\left(\frac{\pi}{2}-\Delta\alpha_3\right)=b_x, \end{aligned} \tag{20}$$

whence the corrected parameters can be found as follows

$$\begin{aligned} &\hat{a}=\frac{a\cos^n\Delta\alpha_3+b_x\sin^n\Delta\alpha_3}{\sin^{m+n}\Delta\alpha_3+\cos^{m+n}\Delta\alpha_3},\\ &\hat{b}_x=\frac{b_x\cos^m\Delta\alpha_3-a\sin^m\Delta\alpha_3}{\sin^{m+n}\Delta\alpha_3+\cos^{m+n}\Delta\alpha_3}. \end{aligned} \tag{21}$$

The shape of Leaf type loops ($n$=1) with $m$=1 does not depend on the phase shift $\Delta\alpha_3$.

The effect of the phase shift $\Delta\alpha_3=\Delta\alpha$ is opposite to the one of $\Delta\alpha_2=\Delta\alpha$ (compare Fig. 5 with Fig. 4). Moreover, the loop built with $\Delta\alpha_2=\Delta\alpha$ is somewhat different from the loop built with $\Delta\alpha_3=-\Delta\alpha$. The difference grows as $|\Delta\alpha|$ increases. The phase shift $\Delta\alpha_3$ causes an offset of parameter $\alpha$ by value $\Delta\alpha_3$. As a result, parameter $\alpha$ becomes equal to $-\Delta\alpha_3$ at the split point $a$ and to $\pi/2-\Delta\alpha_3$ at the saturation point $b$.

***d. Using several phase shifts simultaneously***

In the general case when all the three phase shifts $\Delta\alpha_1$, $\Delta\alpha_2$, and $\Delta\alpha_3$ are used, the system of equations composed for the split point $\alpha$=0 and for the saturation point $\alpha$=π/2 will look as follows

$$
\begin{aligned}
&\hat{a}\cos^{m}\left(0+\Delta\alpha_1-\Delta\alpha_3\right)+\hat{b}_x\sin^{n}\left(0+\Delta\alpha_2-\Delta\alpha_3\right)=a,\\
&\hat{a}\cos^{m}\left(\frac{\pi}{2}+\Delta\alpha_1-\Delta\alpha_3\right)+\hat{b}_x\sin^{n}\left(\frac{\pi}{2}+\Delta\alpha_2-\Delta\alpha_3\right)=b_x.
\end{aligned}
\tag{22}
$$

Solving the system (22), the searched for corrected parameters are found

$$
\begin{aligned}
\hat{a}&=\frac{a\cos^{n}\left(\Delta\alpha_2-\Delta\alpha_3\right)-b_x\sin^{n}\left(\Delta\alpha_2-\Delta\alpha_3\right)}{\sin^{m}\left(\Delta\alpha_1-\Delta\alpha_3\right)\sin^{n}\left(\Delta\alpha_2-\Delta\alpha_3\right)+\cos^{m}\left(\Delta\alpha_1-\Delta\alpha_3\right)\cos^{n}\left(\Delta\alpha_2-\Delta\alpha_3\right)},\\
\hat{b}_x&=\frac{a\sin^{m}\left(\Delta\alpha_1-\Delta\alpha_3\right)+b_x\cos^{m}\left(\Delta\alpha_1-\Delta\alpha_3\right)}{\sin^{m}\left(\Delta\alpha_1-\Delta\alpha_3\right)\sin^{n}\left(\Delta\alpha_2-\Delta\alpha_3\right)+\cos^{m}\left(\Delta\alpha_1-\Delta\alpha_3\right)\cos^{n}\left(\Delta\alpha_2-\Delta\alpha_3\right)}.
\end{aligned}
\tag{23}
$$

By substituting the obtained parameters $\hat{a}$, $\hat{b}_x$ into (13), it is easy to show that with $m$=1, loops of the Leaf type ($n$=1) do not depend on the phase shifts $\Delta\alpha_1$, $\Delta\alpha_2$ and the shape of these loops does not depend on the phase shift $\Delta\alpha_3$. Formulas (16), (19), (21) are particular cases of formulas (23).

Although the effect of the phase shift $\Delta\alpha_3$ is opposite to the effect of the phase shift $\Delta\alpha_2$, they do not completely neutralize each other when used jointly $\Delta\alpha_2=\Delta\alpha_3=\Delta\alpha$ ($\Delta\alpha_1=0$). When the two phase shifts are set equal to each other $\Delta\alpha_2=\Delta\alpha_3=\Delta\alpha$, their joint action is opposite to the action of the phase shift $\Delta\alpha_1=\Delta\alpha$. Therefore, in order to tilt a hysteresis loop in the split point, the phase shifts $\Delta\alpha_2=\Delta\alpha_3=-\Delta\alpha$ can be used instead of the phase shift $\Delta\alpha_1=\Delta\alpha$. Since the phase shifts $\Delta\alpha_2$ and $\Delta\alpha_3=-\Delta\alpha_2$ produce similar effect, the practical application of the model (13) may involve only two of them: $\Delta\alpha_1$ and $\Delta\alpha_2$ (see Sec. III) or $\Delta\alpha_1$ and $\Delta\alpha_3$. When using the phase shifts $\Delta\alpha_1$ and $\Delta\alpha_2$, formulas for the corrected parameters $\hat{a}$ and $\hat{b}_x$ are obtained from formulas (23) by substituting $\Delta\alpha_3=0$; when using the phase shifts $\Delta\alpha_1$ and $\Delta\alpha_3$ – by substituting $\Delta\alpha_2=0$.

It is noteworthy that with $\Delta\alpha_1=\Delta\alpha_2=\Delta\alpha_3=\Delta\alpha\neq0$, where $\Delta\alpha$ is an arbitrary real number, the corrected parameters $\hat{a}$, $\hat{b}_x$ of the improved model (13) degenerate into the parameters $a$, $b_x$ of the original model (1), respectively. Thus, in the case under consideration, the hysteresis loops built by the model (13) coincide in shape with the loops built by the model (1), but the parameter $\alpha$ of the former is shifted by $\Delta\alpha$ as compared to the latter. The conclusion about the degeneration of the model (13) into the model (1) also comes out of the following consideration. Since the effect of the joint usage of the equal phase shifts $\Delta\alpha_2=\Delta\alpha_3=\Delta\alpha$ is opposite to the effect of the phase shift $\Delta\alpha_1=\Delta\alpha$, the joint usage of three identical phase shifts $\Delta\alpha_1=\Delta\alpha_2=\Delta\alpha_3=\Delta\alpha$ should lead to the formation of loop (1).

Like the original model (1), the improved model (13) can also be represented as a sum of harmonic oscillations

$$
\begin{aligned}
&x(\alpha)=\frac{\hat{a}}{2^{m-1}}\sum_{k=0}^{\frac{m-1}{2}}C_m^k\cos\left((m-2k)\left(\alpha+\Delta\alpha_1\right)\right)+\frac{\hat{b}_x}{2^{n-1}}\sum_{k=0}^{\frac{n-1}{2}}(-1)^{\frac{n-1}{2}+k}C_n^k\sin\left((n-2k)\left(\alpha+\Delta\alpha_2\right)\right),\\
&y(\alpha)=b_y\sin\left(\alpha+\Delta\alpha_3\right).
\end{aligned}
\tag{24}
$$

The Fourier coefficients $A_k$, $B_k$ can also be directly extracted from the representation (24):[6]

$$
\begin{aligned}
A_k&=\frac{\hat{a}}{2^{m-1}}C_m^{\frac{m-k}{2}}\cos\left(k\Delta\alpha_1\right)+(-1)^{\left\lfloor\frac{k-1}{2}\right\rfloor}\frac{\hat{b}_x}{2^{n-1}}C_n^{\frac{n-k}{2}}\sin\left(k\Delta\alpha_2\right),\\
B_k&=\frac{-\hat{a}}{2^{m-1}}C_m^{\frac{m-k}{2}}\sin\left(k\Delta\alpha_1\right)+(-1)^{\left\lfloor\frac{k-1}{2}\right\rfloor}\frac{\hat{b}_x}{2^{n-1}}C_n^{\frac{n-k}{2}}\cos\left(k\Delta\alpha_2\right).
\end{aligned}
\tag{25}
$$

Having these coefficients, it is easy to determine the amplitudes $Am_k$ and the phases $\varphi_k$ of harmonics of the generating function $x(\alpha)$ (see Eqs. (8), (12)).

According to definition (13), in order to produce a harmonic signal at the output $y(\alpha)$ of a hysteresis element, a signal of the form $x(\alpha)=\hat{a}\cos^m(\alpha+\Delta\alpha_1)+\hat{b}_x\sin^n(\alpha+\Delta\alpha_2)$ (that is a sum of harmonics (7) in accordance with spectral representation (24)) should be applied to its input. The proper initial phase shift $\Delta\alpha_3$ of the harmonic signal is obtained by setting the corresponding corrected parameters $\hat{a}$, $\hat{b}_x$ calculated by formulas (23). That being said, it follows that the hysteresis loop (13) has a filtering capability.

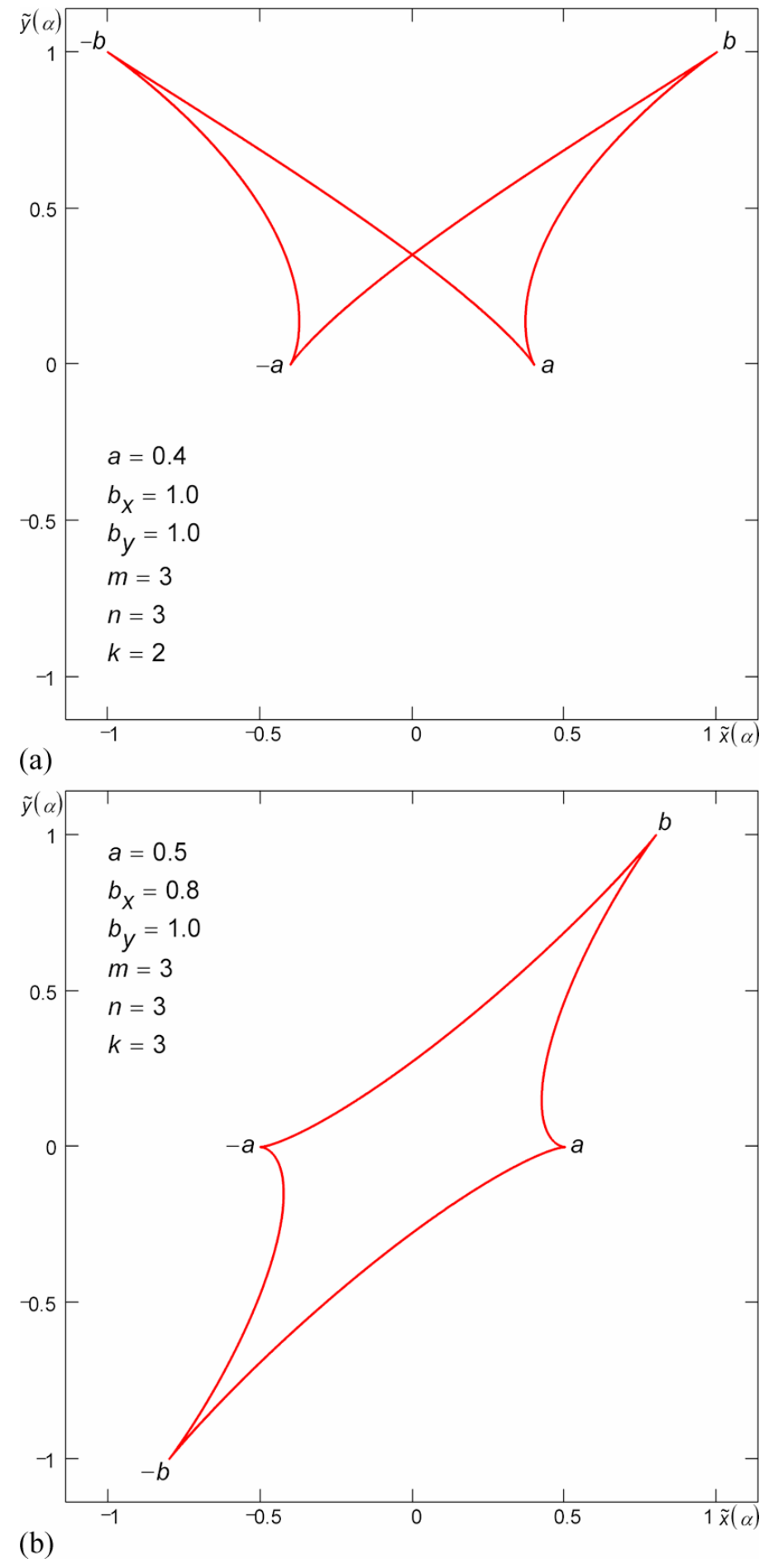

Fig. 6 Hysteresis loop (a) Bat (Butterfly), (b) Astro.

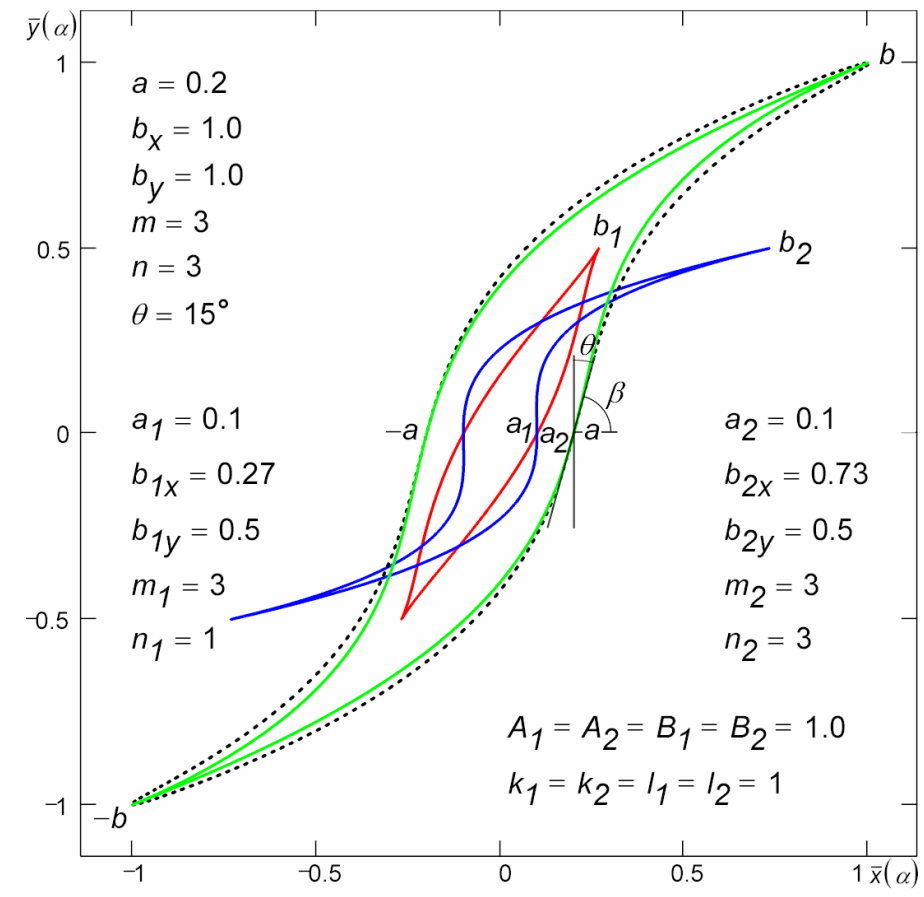

Fig. 7. Tilted Classical loop (green) as a result of the addition of two hysteresis loops – the loop 1 Leaf (red) and the loop 2 Classical (blue). The Classical loop tilted by the same angle $\theta$=15° at the split point $a$ by means of the phase shift $\Delta\alpha_1$ is shown for comparison in the dotted line.

***3. Hysteresis loops of Bat and Astro types***

In Ref. 1, the hysteresis loop of Bat (Butterfly) type was obtained by taking the absolute value of generating function $y(\alpha)$, i. e., $\tilde{y}(\alpha)=\left|b_y\sin\alpha\right|$. A similar type of loop can also be obtained by raising $\sin\alpha$ in $y(\alpha)$ to an even power $k$=2, 4, …, i. e.,

$$\begin{aligned}\tilde{x}(\alpha)&=x(\alpha),\\ \tilde{y}(\alpha)&=b_y\sin^k\alpha.\end{aligned}\qquad(26)$$

The Bat type hysteresis loop with $k$=2 is shown in Fig. 6(a). Odd powers $k$=3, 5, … will yield loops of Astro[8, 9] type. The hysteresis loop of Astro type with $k$=3 is shown in Fig. 6(b).

***4. Arithmetic operations on hysteresis loops***

The hysteresis loops (1)-(5), (7), (9), (11), (13) and their components can be added, subtracted, multiplied by a number, raised to a power, and used as arguments in some functions

$$\begin{aligned}\bar{x}(\alpha)&=\sum_i A_i x_i^{k_i}(\alpha),\\ \bar{y}(\alpha)&=\sum_i B_i y_i^{l_i}(\alpha),\end{aligned}\qquad(27)$$

where $A_i$, $B_i$ are real factors; $k_i$, $l_i$ are positive integer and usually odd powers. The arithmetic operations (27) defined on the loops allow changing the shape, tilt and curvature of hysteresis loops as well as making double loops (see Sec. II.C.4). Addition of loops of different types is admissible.

An example of the addition of two hysteresis loops of different types – the Leaf loop and the Classical loop

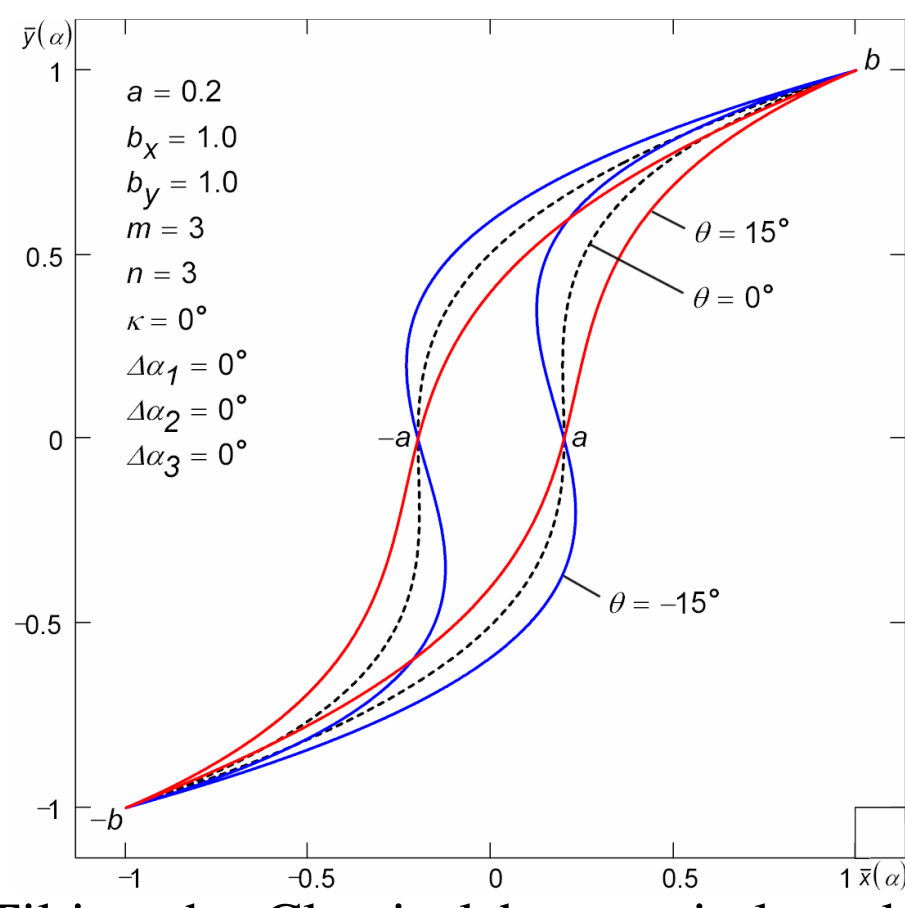

Fig. 8. Tilting the Classical hysteresis loop by skewing the coordinate system by angle $\theta$ along the $x$ axis. The area of all the loops is the same for any oblique angle $\theta$.

($a_2$=$a$-$a_1$, $b_{2x}$=$b_x$-$b_{1x}$, $b_{2y}$=$b_y$-$b_{1y}$) is shown in Fig. 7. The parameters $a$, $b_x$, $b_y$, $m$, $\theta$ of the resulting loop are given in the top left corner of the figure. It should be noted that by adding with a Leaf type loop, Classical loops can be tilted at the split point $a$ by angle $\theta$ ($\beta$=π/2-$\theta$ is a slope angle of tangent at point $a$). By taking into account the preset angle $\theta$, the loop parameter $b_{1x}$ (or $b_{1y}$, $b_{2y}$, $b_y$) can be found according to the formula[6]

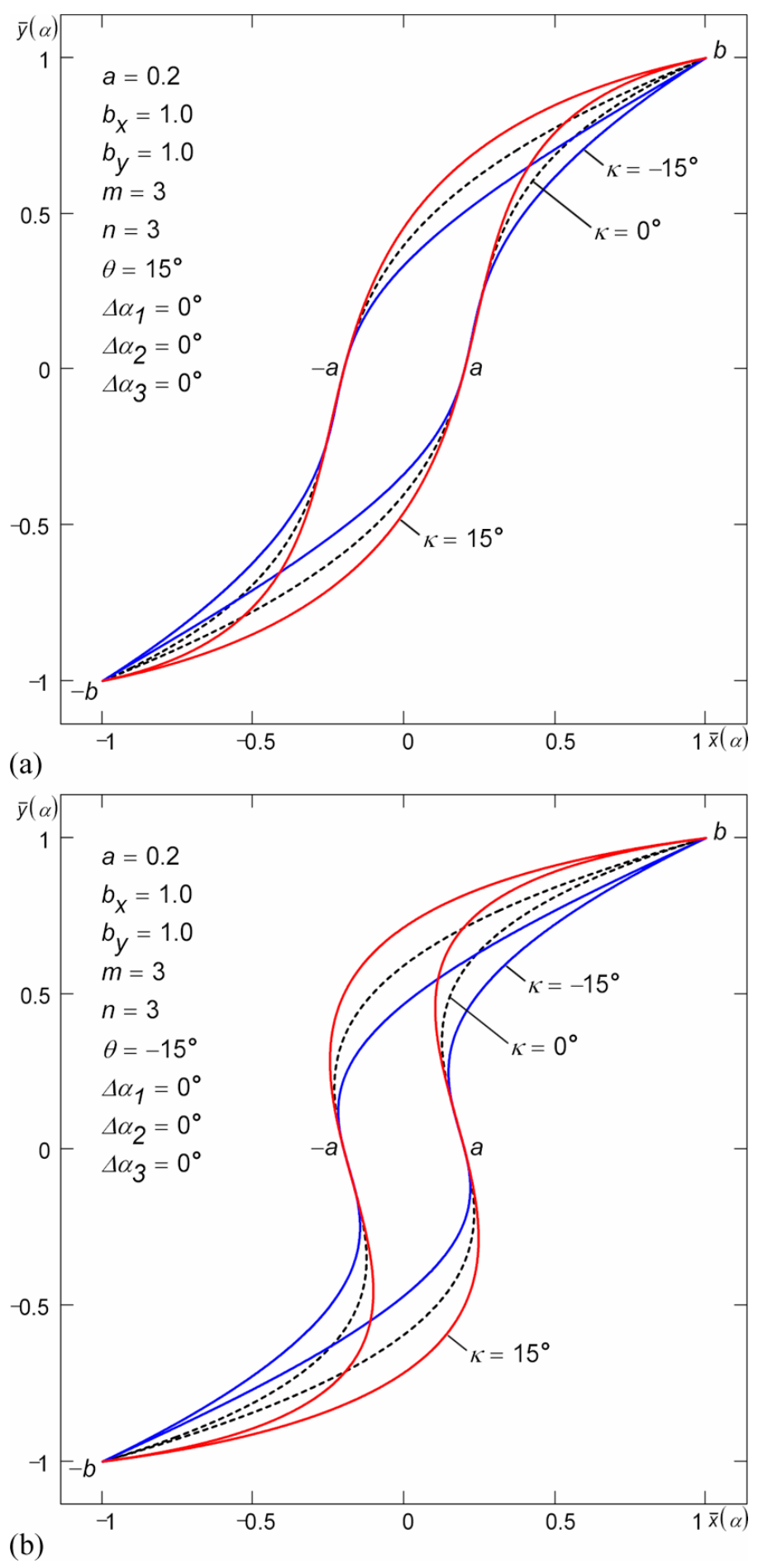

Fig. 9. Changing the curvature of the Classical hysteresis loop by skewing the coordinate system by angle $\kappa$ along the $y$ axis. Tilting the loop at the split point by angle (a) $\theta$=15°, (b) $\theta$=-15°. Loops tilted by any angle $\theta$ have the same areas provided that other parameters of the loop are the same.

$$b_{1x} = \left(b_{1y} + b_{2y}\right)\tan\theta = b_y \tan\theta. \tag{28}$$

In the general case, any hysteresis loop can be decomposed into an infinite number of loops, since according to the above method, each loop from the pair of added loops can always be represented as a pair of some other loops.

It is shown in the supplementary material that the result of the following transformation ($\Delta\alpha_1$=$\Delta\alpha_2$=$\Delta\alpha_3$=0)

$$\begin{aligned} \bar{x}(\alpha) &= x(\alpha) + b_y \tan\theta\left(\sin\alpha - \sin^n\alpha\right), \\ \bar{y}(\alpha) &= y(\alpha), \end{aligned} \tag{29}$$

built on skewing of the coordinate system by angle $\theta$ along the $x$ axis is equivalent to tilting the loop by addition with a Leaf loop (see Fig. 7). Fig. 8 demonstrates the Classical hysteresis loops tilted at the split point according to transformation (29).

An additional skewing of the coordinate system by angle $\kappa$ along the $y$ axis[6]

$$\begin{aligned} \bar{x}(\alpha) &= x(\alpha) + \tan\theta\left(b_x \tan\kappa + b_y\right)\left(\sin\alpha - \sin^n\alpha\right), \\ \bar{y}(\alpha) &= y(\alpha) + b_x \tan\kappa\left(\sin\alpha - \sin^n\alpha\right) \end{aligned} \tag{30}$$

allows continuously changing the loop curvature. Fig. 9 shows tilted Classical loops having various curvatures that are built according to transformations (30). Transformations (29) are particular case of transformations (30).

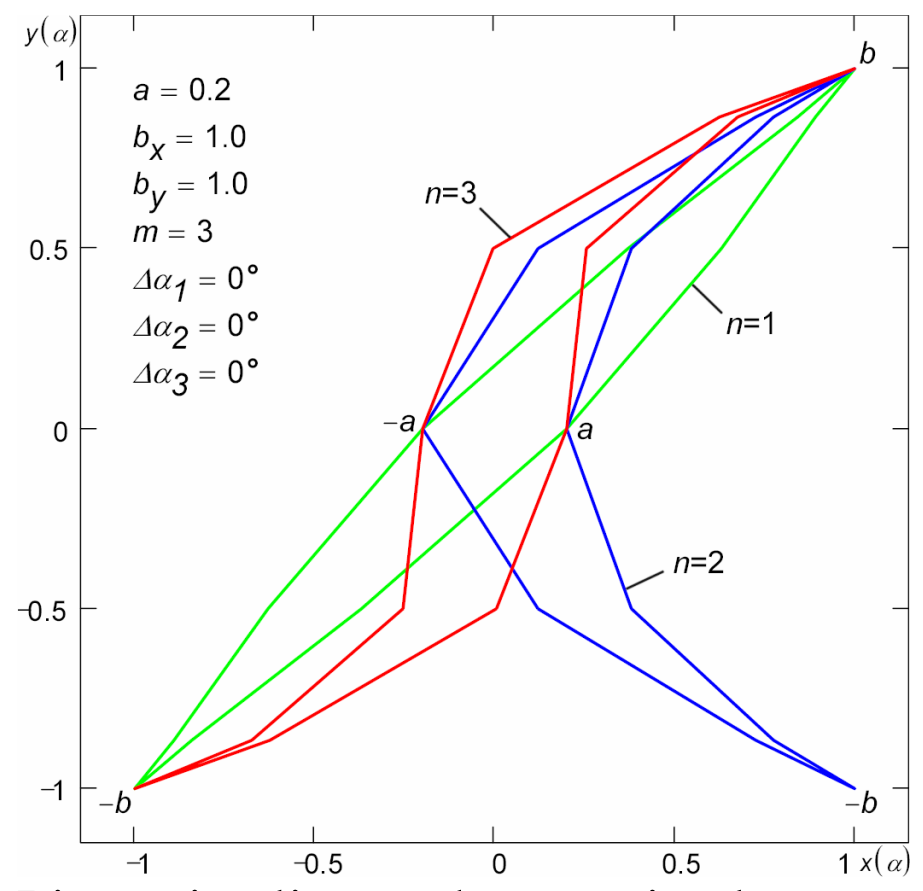

Fig. 10. Piecewise-linear hysteresis loops of Leaf ($n$=1), Crescent (Boomerang, $n$=2), and Classical ($n$=3) types. The number of division intervals $k$=12 (six-linear loop). The areas of the loops are approximately the same.

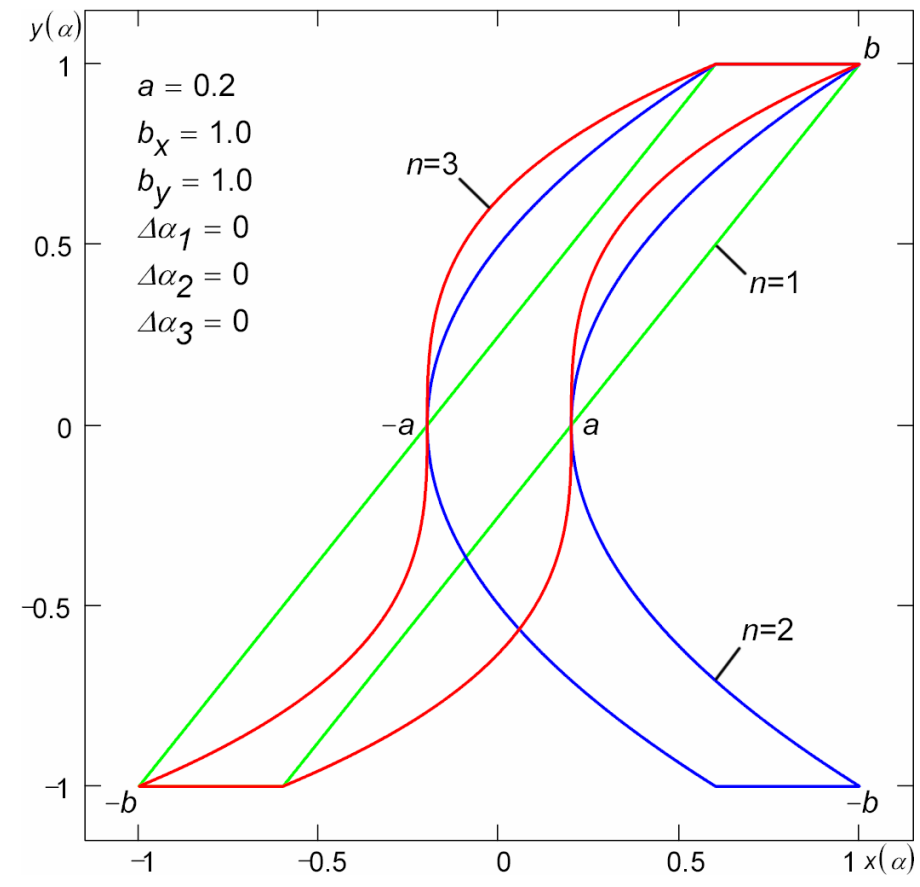

Fig. 11. Piecewise-linear hysteresis loop Leaf (Play without Whiskers, $n$=1), hybrid Crescent (hybrid Boomerang, $n$=2), and hybrid Classical ($n$=3) built on trapezoidal pulses. The area of all the three loops is the same.

## B. Piecewise-linear and hybrid hysteresis loops

Within the model under consideration, the simplest way of building piecewise-linear loops is as follows. The period $T$=2π of change of parameter $\alpha$ is divided into the required number $k$ of intervals, where $k$ is an even integer ($k$≥4). Continuous values $\alpha$ in formula (13) are then replaced with values changing with step $T/k$ so obtaining the $x$, $y$ coordinates of the points of a piecewise-linear hysteresis loop. Those coordinates are then connected with line segments. An example of piecewise-linear loops of three types Leaf, Crescent (Boomerang), and Classical built for $k$=12 (six-linear loop, $k$/2=6) is given in Fig. 10.

Other methods of obtaining piecewise-linear hysteresis loops imply replacing the sine and the cosine in the generating functions $x(\alpha)$ and $y(\alpha)$ of the model (13) by unit amplitude trapezoidal, triangular, and rectangular pulses and their combinations.

### *1. Loops built on trapezoidal pulses*

Replacing the sine and the cosine in the generating functions $x(\alpha)$ and $y(\alpha)$ of the model (1) with unit-amplitude trapezoidal pulses $\mathrm{trp_s}$ and $\mathrm{trp_c}$, respectively, we can produce piecewise-linear hysteresis loops built on trapezoidal pulses

$$\begin{aligned} x(\alpha) &= a\,\mathrm{trp_c^{\mathit{m}}}\,\alpha + (b_x - a)\,\mathrm{trp_s^{\mathit{n}}}\,\alpha, \\ y(\alpha) &= b_y\,\mathrm{trp_s}\,\alpha, \end{aligned} \tag{31}$$

where $\mathrm{trp_c}(\alpha)=\mathrm{trp_s}(\alpha+T/4)$; $T$ is the pulse period. The expressions defining the trapezoidal pulses trp are given in Ref. 1 and in the supplementary material. The simplest loops are obtained when the upper $d$ and the lower $D$ bases of the trapezoidal pulses relate as $D$=3$d$ ($T$=$d$+$D$=4$d$). In this case, the shape of loops is independent of $m$. Subtracting the splitting $a$ from $b_x$ in (31) allows us to move the saturation point $b$ from the middle of the horizontal section of the loop, where $\alpha$=$T$/4, to its canonical position (see Fig. 11); here parameter $\alpha$ takes the value of $T$/8.

Fig. 11 shows loops of three types: Leaf (Play without Whiskers), Crescent (Boomerang), and Classical, which are built according to formulas (31) ($D$=3$d$). The loop Leaf is piecewise-linear. The loops Crescent and Classical are hybrid loops since they combine rectilinear sections and curvilinear ones. If $b_x$=$a$, the loops (31) degenerate into rectangular loop (Non-ideal Relay without Whiskers).

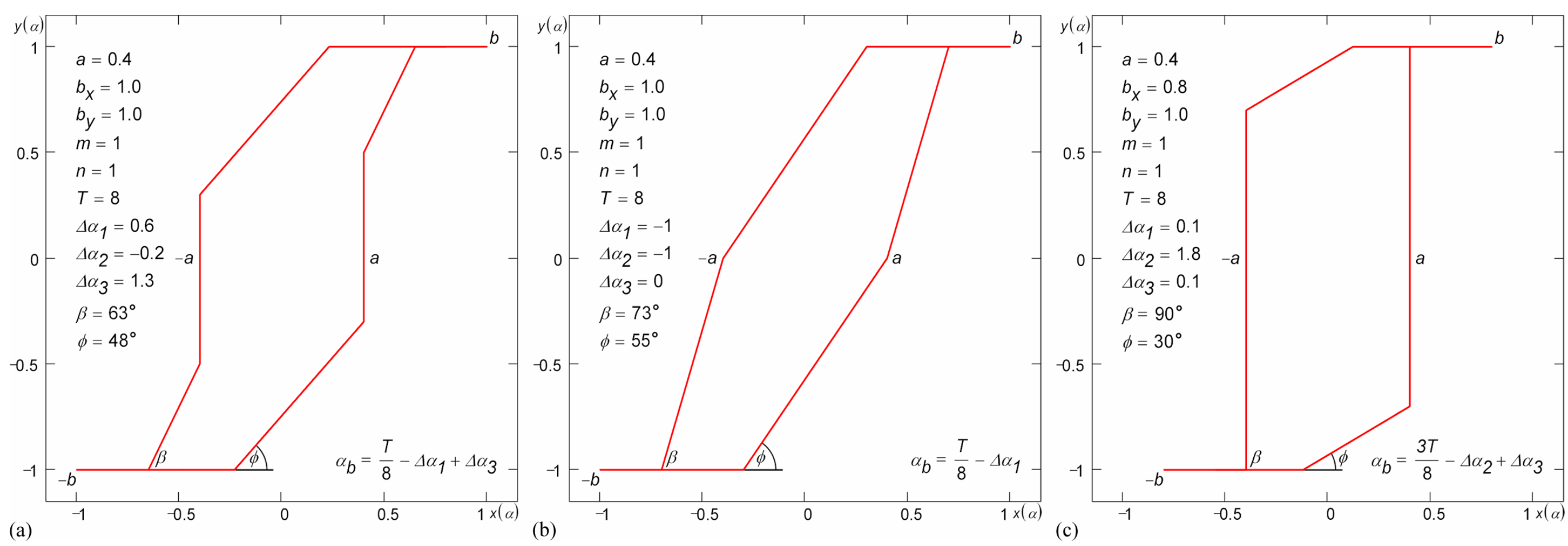

Fig. 12. Piecewise-linear hysteresis loops: (a) Play-Relay-Play (tetra-linear loop), (b) Play-Play (tri-linear loop), (c) Play-Relay (tri-linear loop) built on trapezoidal pulses using the phase shifts.

***a. Piecewise-linear loops***

Taking into account the phase shifts $\Delta\alpha_1$, $\Delta\alpha_2$, $\Delta\alpha_3$, equations (31) are written as

$$\begin{aligned} x(\alpha) &= \hat{a}\,\mathrm{trp}_\mathrm{c}^m(\alpha+\Delta\alpha_1)+\hat{b}_x\,\mathrm{trp}_\mathrm{s}^n(\alpha+\Delta\alpha_2), \\ y(\alpha) &= b_y\,\mathrm{trp}_\mathrm{s}(\alpha+\Delta\alpha_3). \end{aligned} \tag{32}$$

The corrected parameters $\hat{a}$ and $\hat{b}_x$ are determined by the following formulas[6]

$$\begin{aligned} \hat{a} &= \frac{a\,\mathrm{trp}_\mathrm{s}^n(\alpha_b+\Delta\alpha_2-\Delta\alpha_3)-b_x\mathrm{trp}_\mathrm{s}^n(\Delta\alpha_2-\Delta\alpha_3)}{\mathrm{trp}_\mathrm{c}^m(\Delta\alpha_1-\Delta\alpha_3)\mathrm{trp}_\mathrm{s}^n(\alpha_b+\Delta\alpha_2-\Delta\alpha_3)-\mathrm{trp}_\mathrm{c}^m(\alpha_b+\Delta\alpha_1-\Delta\alpha_3)\mathrm{trp}_\mathrm{s}^n(\Delta\alpha_2-\Delta\alpha_3)}, \\ \hat{b}_x &= \frac{b_x\mathrm{trp}_\mathrm{c}^m(\Delta\alpha_1-\Delta\alpha_3)-a\,\mathrm{trp}_\mathrm{c}^m(\alpha_b+\Delta\alpha_1-\Delta\alpha_3)}{\mathrm{trp}_\mathrm{c}^m(\Delta\alpha_1-\Delta\alpha_3)\mathrm{trp}_\mathrm{s}^n(\alpha_b+\Delta\alpha_2-\Delta\alpha_3)-\mathrm{trp}_\mathrm{c}^m(\alpha_b+\Delta\alpha_1-\Delta\alpha_3)\mathrm{trp}_\mathrm{s}^n(\Delta\alpha_2-\Delta\alpha_3)}, \end{aligned} \tag{33}$$

where $\alpha_b$ is equal to $T/8-\Delta\alpha_1+\Delta\alpha_3$, $T/8-\Delta\alpha_1$, $3T/8-\Delta\alpha_2+\Delta\alpha_3$ and other values[6] depending on the applied ranges of phase shifts $\Delta\alpha_1$, $\Delta\alpha_2$, $\Delta\alpha_3$. The formulas (33) for the corrected parameters $\hat{a}$ and $\hat{b}_x$ of the piecewise-linear loops (32) are derived using the same logic as in the case of smooth loops. With $\alpha_b=T/4$, the working formulas for the corrected parameters $\hat{a}$ and $\hat{b}_x$ can be obtained from formulas (16), (19), (21), (23) valid for smooth loops (13) by simply replacing the sine and the cosine functions with trapezoidal pulses $\mathrm{trp}_\mathrm{s}$ and $\mathrm{trp}_\mathrm{c}$, respectively.

Fig. 12 shows the piecewise-linear hysteresis loops built by equations (32) using the corrected parameters (33) ($m$=$n$=1). As an example, the loop Play-Play with different values of phase shift $\Delta\alpha_3$ is demonstrated in Fig. 13. According to the definition (32), when a signal of the form $x(\alpha)=\hat{a}\,\mathrm{trp}_\mathrm{c}^m(\alpha+\Delta\alpha_1)+\hat{b}_x\,\mathrm{trp}_\mathrm{s}^n(\alpha+\Delta\alpha_2)$ is applied to the input of a piecewise-linear hysteresis element, trapezoidal pulses are produced at the output $y(\alpha)$ with the predefined phase shift $\Delta\alpha_3$.

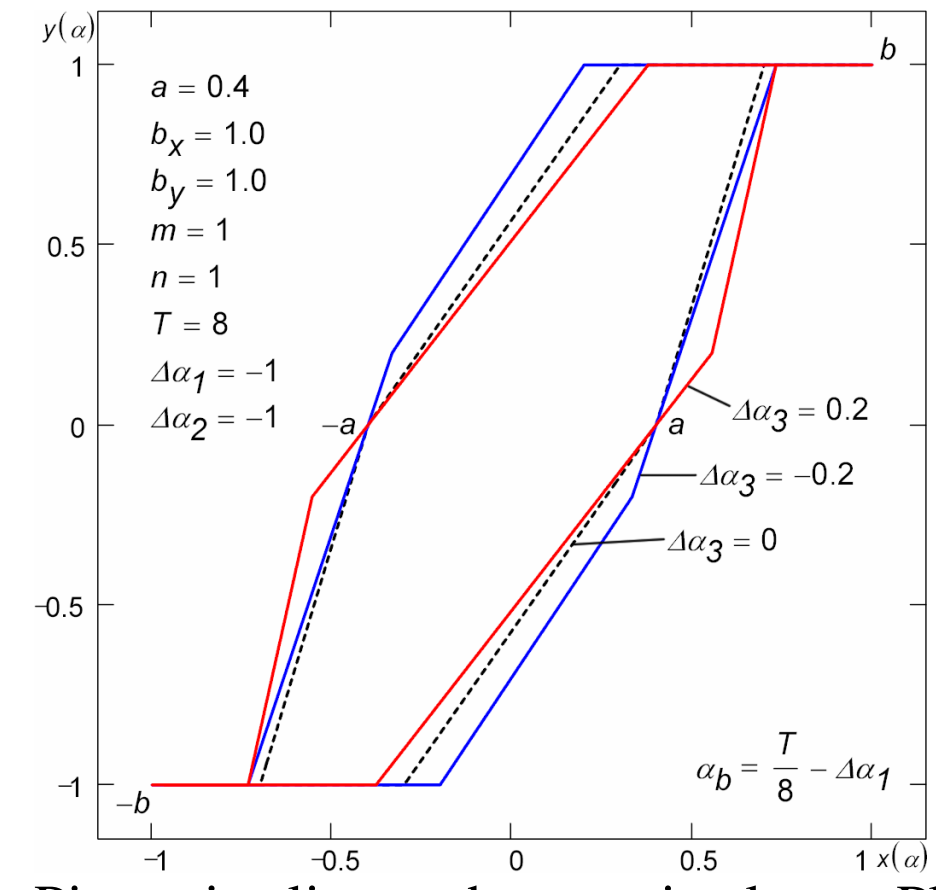

Fig. 13. Piecewise-linear hysteresis loop Play-Play with different values of phase shift $\Delta\alpha_3$.

To obtain piecewise linear loops with gain/attenuation $\gamma$, one should add to (32) an extra term (curve) responsible for gain/attenuation ($m$=$n$=1)[6]

$$\begin{aligned} \bar{x}(\alpha) &= x(\alpha), \\ \bar{y}(\alpha) &= y(\alpha)+\tan\gamma\left[x(\alpha)-b_x\,\mathrm{trp}_\mathrm{s}(\alpha+\Delta\alpha_3)\right]. \end{aligned} \tag{34}$$

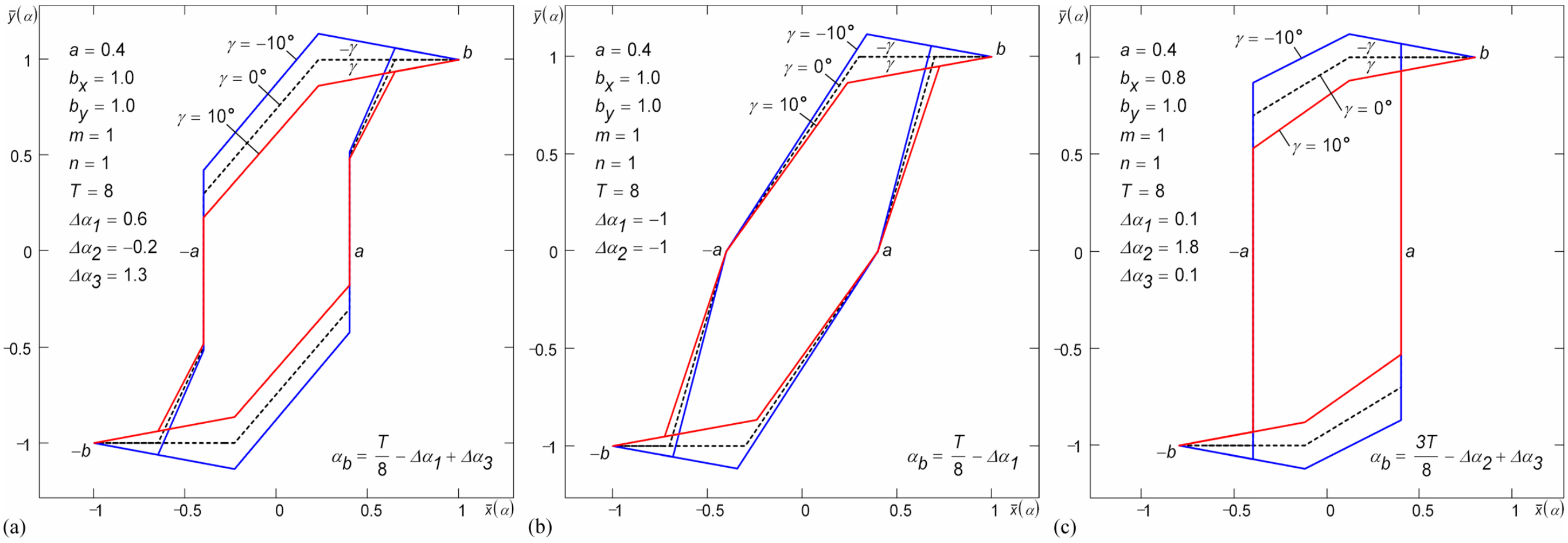

Fig. 14. Piecewise-linear hysteresis loops with gain/attenuation $\gamma$. (a) Play-Relay-Play, (b) Play-Play, (c) Play-Relay built on trapezoidal pulses using the phase shifts.

Fig. 14 shows loops with gain ($\gamma$>0) and with attenuation (negative gain $\gamma$<0) built by Eqs. (34). The loops in Fig. 14 with no gain ($\gamma$=0) correspond to the loops shown in Fig. 12.

Loops of type Play without Whiskers (Leaf) (see Fig. 11) can not be tilted by rotation of coordinate system at the split point. Nevertheless, a gain/attenuation $\gamma$ can be applied to these loops by using the rotation ($m$=$n$=1, $\Delta\alpha_1=\Delta\alpha_2=\Delta\alpha_3=0$)[6]

$$\bar{x}(\alpha)=\frac{a\tan\beta\,\mathrm{trp_c}\,\alpha+\left(b_y-b_x\tan\gamma\right)\mathrm{trp_s}\,\alpha}{\tan\beta-\tan\gamma},$$
$$\bar{y}(\alpha)=\tan\beta\left(\bar{x}(\alpha)-a\,\mathrm{trp_c}\,\alpha\right), \tag{35}$$

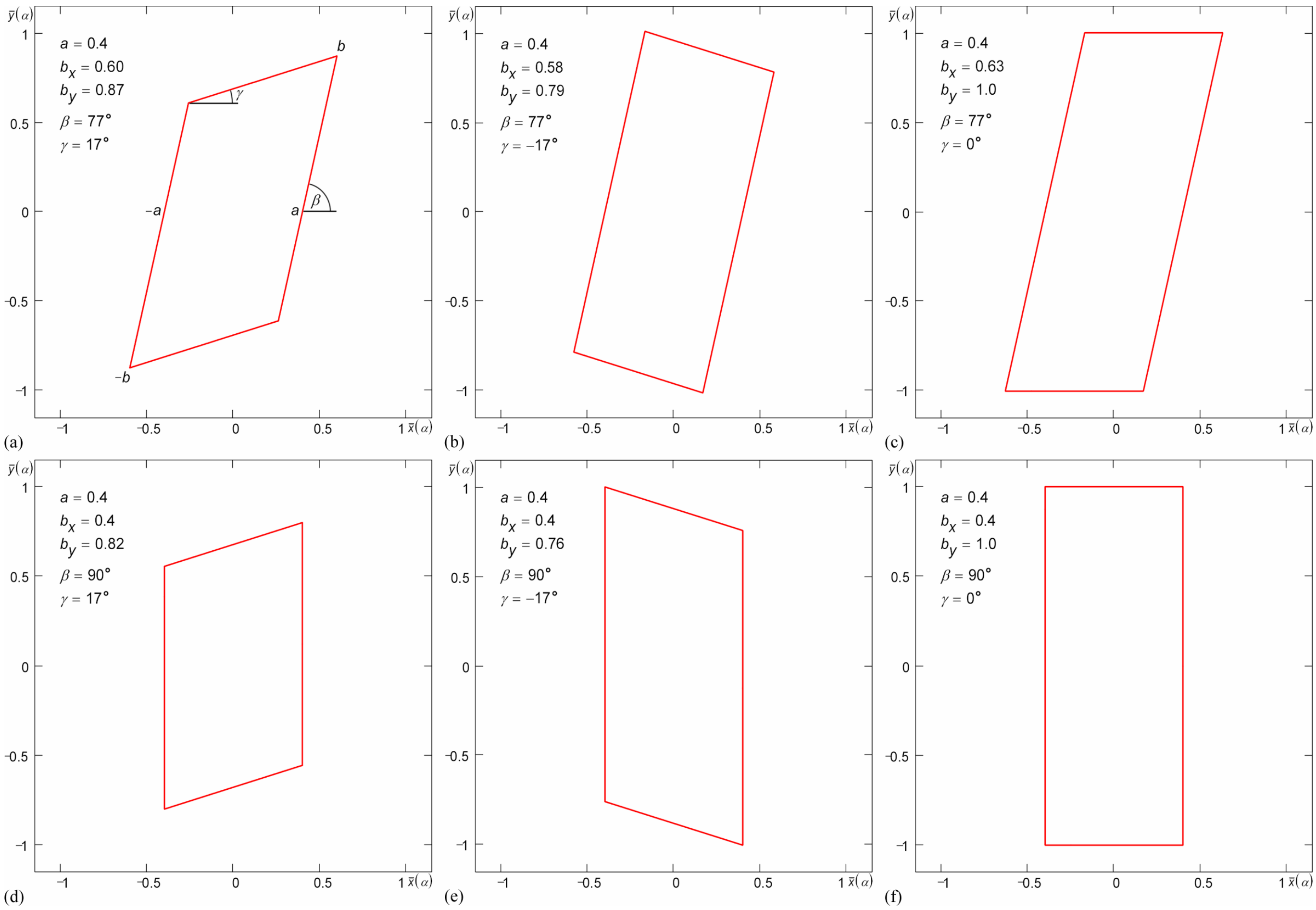

Fig. 15. Piecewise-linear hysteresis loops: (a) Play with Gain without Whiskers (parallelogram, bi-linear, elastic-plastic loop), (b) Play with Attenuation without Whiskers, (c) Play without Whiskers, (d) Non-ideal Relay with Gain without Whiskers, (e) Non-ideal Relay with Attenuation without Whiskers, (f) Non-ideal Relay without Whiskers (rectangular loop) built on trapezoidal pulses.

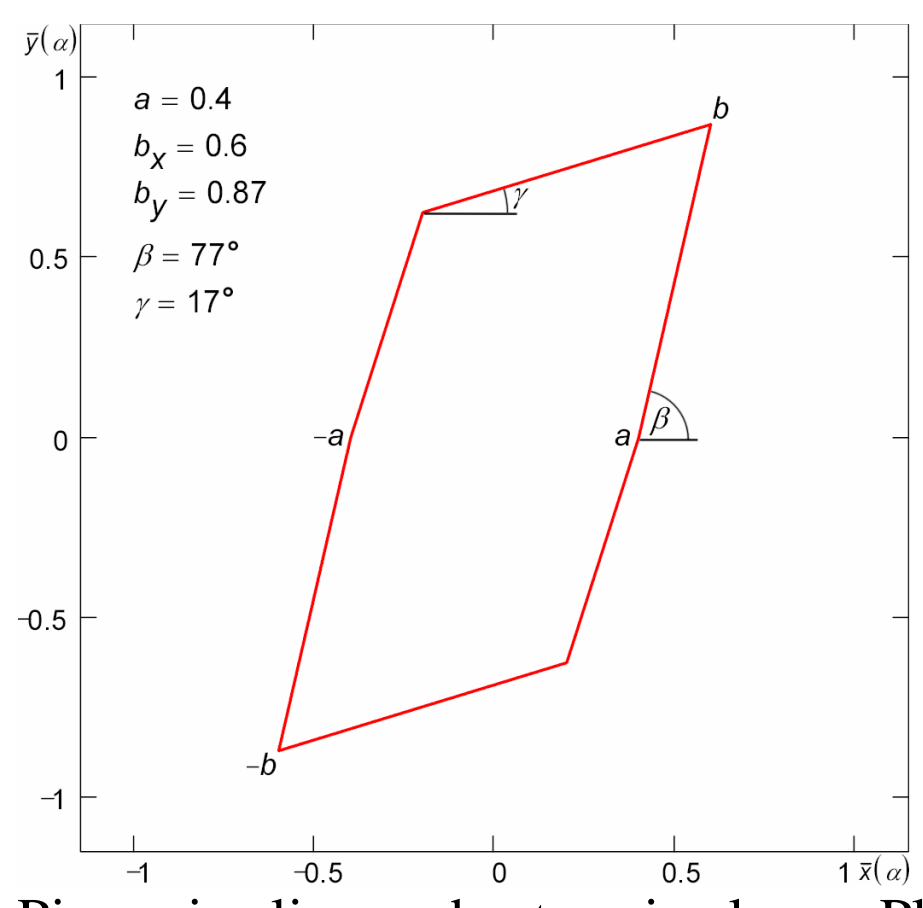

Fig. 16. Piecewise-linear hysteresis loop Play-Play with Gain without Whiskers built on trapezoidal pulses.

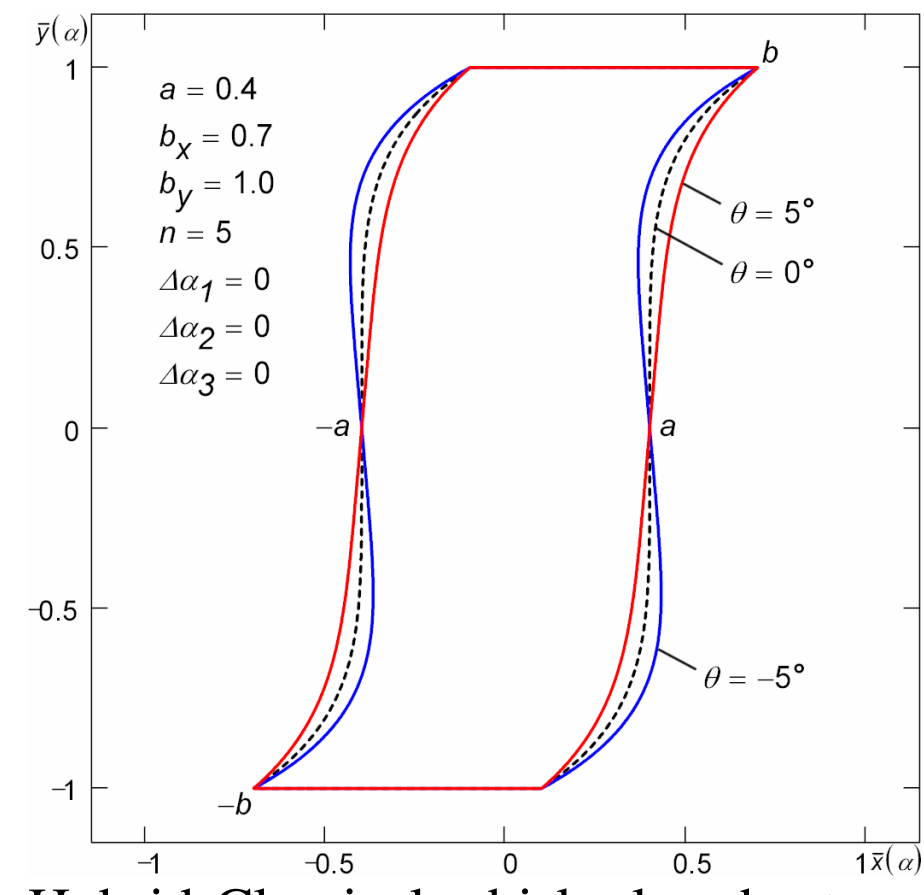

Fig. 17. Hybrid Classical whiskerless hysteresis loops with the specified slope $\beta=\pi/2-\theta$ at the split point. The loops are built on trapezoidal pulses. The area of all the loops is the same regardless of $\theta$.

where the loop parameters $a$, $b_x$, $b_y$, $\beta$ are interrelated as $\tan\beta=b_y/(b_x-a)$. The hysteresis loops built on trapezoidal pulses according to equations (35) are presented in Fig. 15.

Applying a variation of rotation of coordinate system, the following transformations are found ($m=n=1$, $\Delta\alpha_1=\Delta\alpha_2=\Delta\alpha_3=0$)[6]

$$\begin{aligned}\bar{x}(\alpha)&=x(\alpha),\\ \bar{y}(\alpha)&=y(\alpha)+a\tan\gamma\left(\mathrm{trp_c}\,\alpha\left|\mathrm{trp_s}\,\alpha\right|-\mathrm{trp_s}\,\alpha\right),\end{aligned}\qquad(36)$$

according to which loops Play-Play with Gain without Whiskers[10] can be built on trapezoidal pulses (see Fig. 16).

***b. Hybrid loops***

Tilting in the split point by the phase shift $\Delta\alpha_1$ is not appropriate for the hybrid Classical loop ($n=3, 5, \ldots$) shown in Fig. 11 because of corners (derivative jumps) that appear at the curvilinear sections of the loop and because the generating functions $x(\alpha)$, $y(\alpha)$ acquire horizontal sections with the same values of the parameter $\alpha$.[6] The following transformation makes it possible to obtain a loop with the desired tilt $\beta=\pi/2-\theta$ in the split point while keeping the curvilinear sections of the loop smooth and the rectilinear sections of the loop horizontal ($m=1$, $\Delta\alpha_1=\Delta\alpha_2=\Delta\alpha_3=0$):[6]

$$\begin{aligned}\bar{x}(\alpha)&=x(\alpha)+b_y\sin\theta\left(\mathrm{trp_s}\,\alpha-\mathrm{trp_s^{\mathit{n}}}\,\alpha\right),\\ \bar{y}(\alpha)&=y(\alpha).\end{aligned}\qquad(37)$$

Transformation (37) for the hybrid loop is similar to the transformation (29) for the smooth loop. Fig. 17 shows hybrid classical whiskerless hysteresis loops tilted according to transformation (37).[11, 12, 13]

By applying a series of skewing linear transformations,[6] a universal formula is obtained ($m=1$, $\Delta\alpha_1=\Delta\alpha_2=\Delta\alpha_3=0$)

$$\begin{aligned}\bar{x}(\alpha)&=x(\alpha)+\tan\theta\left[(b_x-a)\tan\kappa-b_x\tan\gamma+b_y\right]\left(\mathrm{trp_s}\,\alpha-\mathrm{trp_s^{\mathit{n}}}\,\alpha\right),\\ \bar{y}(\alpha)&=y(\alpha)+\left[(b_x-a)\tan\kappa-b_x\tan\gamma\right]\left(\mathrm{trp_s}\,\alpha-\mathrm{trp_s^{\mathit{n}}}\,\alpha\right)+a\tan\gamma\left(\mathrm{trp_c}\,\alpha\,\mathrm{trp_s^{\mathit{k}}}\,\alpha-\mathrm{trp_s^{\mathit{n}}}\,\alpha\right),\end{aligned}\qquad(38)$$

which allows building hybrid classical whiskerless loops[12, 13] (see Fig. 18) on the base of trapezoidal pulses. The resulting loops have the required slope $\beta$ at the split point ($\beta=\pi/2-\theta$), the required inclination (gain/attenuation) $\gamma$ of the rectilinear section and the required curvature $\kappa$ of the curvilinear section. Introducing additional parameters $k=2, 4, \ldots$ and $\kappa$ besides $n=3, 5, \ldots$ to control the curvature is dictated by the fact that the use of the phase shifts $\Delta\alpha_2$, $\Delta\alpha_3$ would lead to the appearance of undesirable corners in the curvilinear sections of the loop. In formula (38), the function $\mathrm{trp}_\mathrm{s}^k\,\alpha$ could be replaced with $\left|\mathrm{trp}_\mathrm{s}^k\,\alpha\right|$, where $k$ is any positive number. In case $\gamma=\kappa=0$, formulas (38) degenerate into formulas (37).

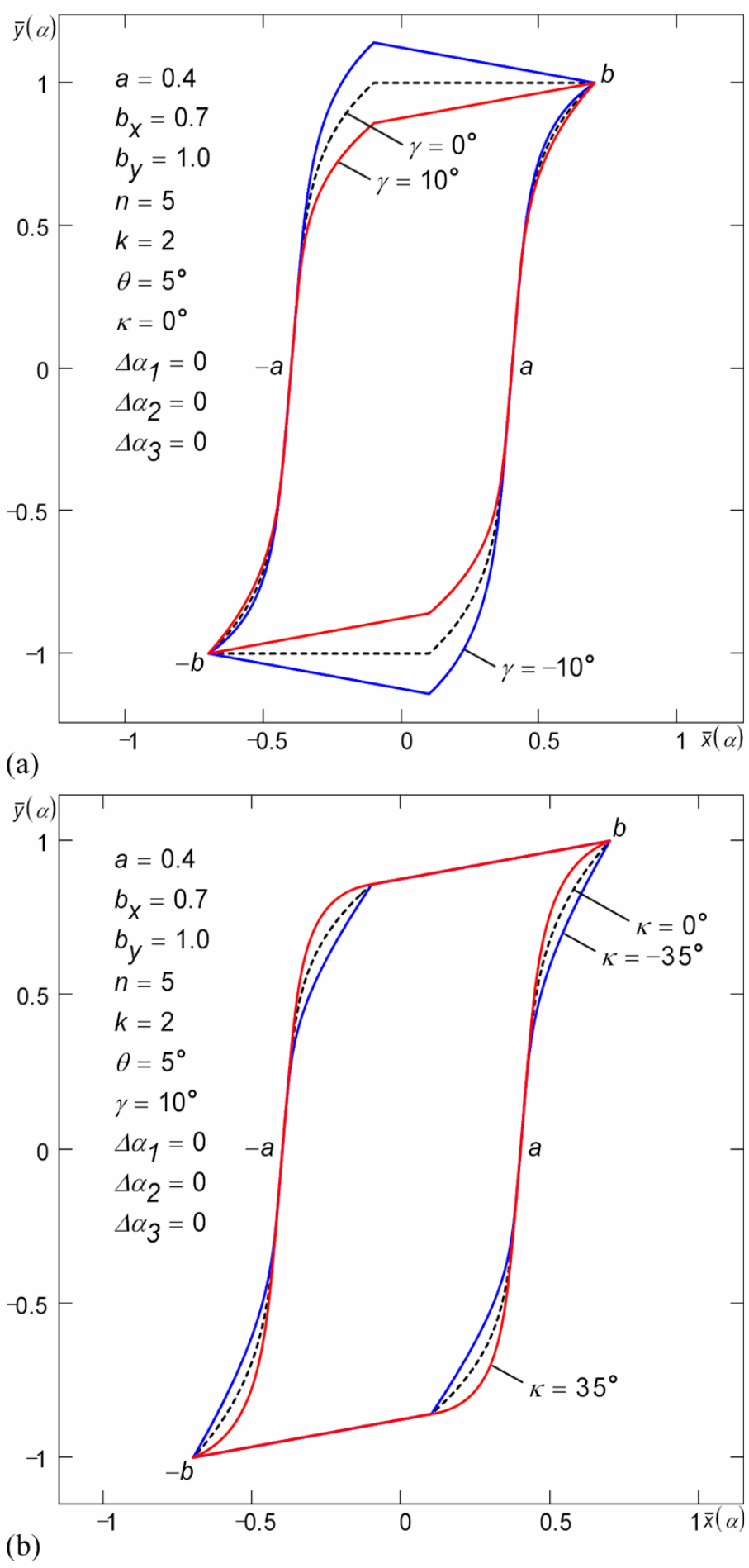

Fig. 18. Hybrid Classical whiskerless hysteresis loops with specified slope $\beta=\pi/2-\theta$, gain/attenuation $\gamma$, and curvature $\kappa$. (a) Various gains $\gamma$ for fixed $\beta$ and $\kappa$; (b) various curvatures $\kappa$ for fixed $\beta$ and $\gamma$. The loops are built on trapezoidal pulses.

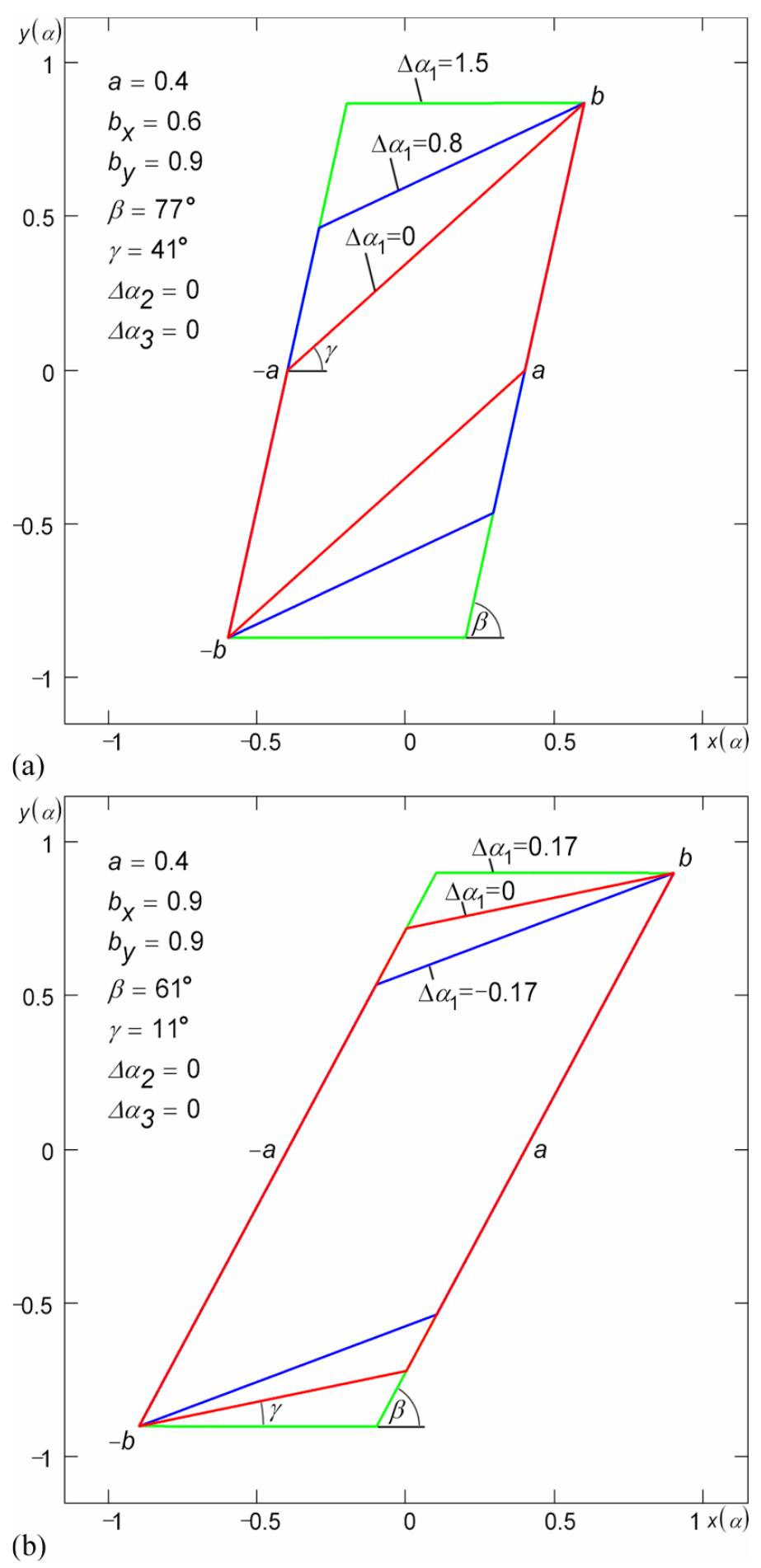

Fig. 19. Piecewise-linear loop Leaf (bi-linear loop, Play without Whiskers) at different values of phase shift $\Delta\alpha_1$ built on triangular pulses by using (a) base, (b) inverse to base equations.

If necessary, the formulas (31)-(38) can be expressed through triangular pulses by using the following representation of trapezoidal pulses as a sum of two triangular pulses ($D=3d$, $T=d+D$):

$$\begin{aligned}\mathrm{trp}_\mathrm{s}(\alpha)&=\mathrm{tri}_\mathrm{s}(\alpha+\frac{T}{8})+\mathrm{tri}_\mathrm{s}(\alpha-\frac{T}{8}),\\ \mathrm{trp}_\mathrm{c}(\alpha)&=\mathrm{tri}_\mathrm{c}(\alpha+\frac{T}{8})+\mathrm{tri}_\mathrm{c}(\alpha-\frac{T}{8}).\end{aligned}\tag{39}$$

### *2. Loops built on triangular pulses*

Besides trapezoidal pulses trp, formulas (13) can operate with triangular pulses tri,[1, 6] which are particular cases of trapezoidal pulses ($d$=0, $T=D$):

$$\begin{aligned}x(\alpha)&=\hat{a}\,\mathrm{tri}_\mathrm{c}^m(\alpha+\Delta\alpha_1)+\hat{b}_x\,\mathrm{tri}_\mathrm{s}^n(\alpha+\Delta\alpha_2),\\ y(\alpha)&=b_y\,\mathrm{tri}_\mathrm{s}(\alpha+\Delta\alpha_3).\end{aligned}\tag{40}$$

The working formulas for the corrected parameters $\hat{a}$ and $\hat{b}_x$ are obtained from formulas (16), (19), (21), (23) valid for smooth loops (13) by simply replacing the sine and the cosine functions with triangular pulses $\mathrm{tri}_\mathrm{s}$ and $\mathrm{tri}_\mathrm{c}$, respectively.

Fig. 19(a) presents a piecewise-linear loop Leaf ($m=n=1$) built according to formulas (40) with various positive values of phase shift $\Delta\alpha_1$ ($\Delta\alpha_2=\Delta\alpha_3=0$). With $\Delta\alpha_1=0.8$, the loop is a variety of the loop Play with Gain without Whiskers (see Fig. 20(b)); with $\Delta\alpha_1$=1.5, the loop is a variety of the loop Play without Whiskers (see Fig. 20(f)). As can be seen from Fig. 19(a), it is possible to adjust the gain $\gamma$ using the phase shift $\Delta\alpha_1$.

The angles $\beta$ and $\gamma$ shown in the figure are determined by the formulas

$$\begin{aligned} \tan\beta &= \frac{b_y}{\hat{b}_x-\hat{a}} = \frac{b_y\tan\gamma}{b_y-2\hat{a}\tan\gamma} = \frac{b_y\,\mathrm{tri_c}\,\Delta\alpha_1}{a\left(\mathrm{tri_s}\,\Delta\alpha_1-1\right)+b_x\,\mathrm{tri_c}\,\Delta\alpha_1}, \\ \tan\gamma &= \frac{b_y}{\hat{b}_x+\hat{a}} = \frac{b_y\tan\beta}{b_y+2\hat{a}\tan\beta} = \frac{b_y\,\mathrm{tri_c}\,\Delta\alpha_1}{a\left(\mathrm{tri_s}\,\Delta\alpha_1+1\right)+b_x\,\mathrm{tri_c}\,\Delta\alpha_1}. \end{aligned} \tag{41}$$

Using formulas (41), equations (40) can be expressed, if necessary, through the angles $\beta$ and/or $\gamma$. With $b_x$=$a$, the loop (40) degenerates into the loop variety Non-ideal Relay with Gain without Whiskers (see Fig. 20(h)).

Just as it was done above for loops (35) built on trapezoidal pulses, the gain/attenuation $\gamma$ of loops (40) can be changed by rotating the coordinate system ($m$=$n$=1, $\Delta\alpha_1$=$\Delta\alpha_2$=$\Delta\alpha_3$=0)[6]

$$\begin{aligned} \bar{x}(\alpha) &= b_x\,\mathrm{tri_s}\,\alpha + \left(\frac{2a\tan\beta}{\tan\beta-\tan\gamma}-b_x\right)\mathrm{tri_c}\,\alpha, \\ \bar{y}(\alpha) &= b_y\,\mathrm{tri_s}\,\alpha + \frac{\left(a+b_x\right)\tan\gamma-b_y}{\tan\beta-\tan\gamma}\tan\beta\,\mathrm{tri_c}\,\alpha, \end{aligned} \tag{42}$$

where the loop parameters $a$, $b_x$, $b_y$, $\beta$ are interrelated as $\tan\beta$=$b_y$/($b_x$-$a$). With the help of equations (42), it is possible to build piecewise-linear loops shown in Figs. 20(b), (d), (f), (h), (j), (l).

Another method to describe the Leaf type loop is based on function inverse to (40) ($m$=$n$=1)[6]

$$\begin{aligned} x(\alpha) &= b_x\,\mathrm{tri_s}\left(\alpha+\Delta\alpha_3\right), \\ y(\alpha) &= \hat{b}_y\left[\frac{\hat{a}}{\hat{a}-b_x}\,\mathrm{tri_c}\left(\alpha+\Delta\alpha_1\right)+\mathrm{tri_s}\left(\alpha+\Delta\alpha_2\right)\right]. \end{aligned} \tag{43}$$

The corrected parameters $\hat{a}$ and $\hat{b}_y$ in (43) are determined by the following formulas[6]

$$\begin{aligned} \hat{a} &= \frac{b_x\,\mathrm{tri_s}\left(\alpha_a+\Delta\alpha_2-\Delta\alpha_3\right)}{\mathrm{tri_c}\left(\alpha_a+\Delta\alpha_1-\Delta\alpha_3\right)+\mathrm{tri_s}\left(\alpha_a+\Delta\alpha_2-\Delta\alpha_3\right)}, \\ \hat{b}_y &= \frac{b_y\,\mathrm{tri_c}\left(\alpha_a+\Delta\alpha_1-\Delta\alpha_3\right)}{\mathrm{tri_c}\left(\alpha_a+\Delta\alpha_1-\Delta\alpha_3\right)\mathrm{tri_c}\left(\Delta\alpha_2-\Delta\alpha_3\right)+\mathrm{tri_s}\left(\alpha_a+\Delta\alpha_2-\Delta\alpha_3\right)\mathrm{tri_s}\left(\Delta\alpha_1-\Delta\alpha_3\right)}, \end{aligned} \tag{44}$$

where $\alpha_a$=$aT$/(4$b_x$) is the value of parameter $\alpha$ at the split point $a$. With $b_x$>2$a$ ($\Delta\alpha_1$=$\Delta\alpha_2$=$\Delta\alpha_3$=0), the loop (43) is a variety of the loop Play with Gain without Whiskers (see Fig. 20(b)). With $b_x$<2$a$ ($\Delta\alpha_1$=$\Delta\alpha_2$=$\Delta\alpha_3$=0), the loop (43) is a variety of the loop Play with Attenuation without Whiskers (see Fig. 20(d)). With $b_x$=2$a$ ($\Delta\alpha_1$=$\Delta\alpha_2$=$\Delta\alpha_3$=0), the loop (43) degenerates into the loop variety Play without Whiskers (see Fig. 20(f)).

Fig. 19(b) shows a Leaf piecewise-linear loop built by formulas (43) at different values of the phase shift $\Delta\alpha_1$ ($b_x$>2$a$, $\Delta\alpha_2$=$\Delta\alpha_3$=0). As can be seen from the figure, it is possible to adjust the gain $\gamma$ using the phase shift $\Delta\alpha_1$. The angles $\beta$ and $\gamma$ shown in the figure are determined by the formulas[6]

$$\begin{aligned} \tan\beta &= \frac{\hat{b}_y}{b_x-\hat{a}} = \frac{b_x\tan\gamma}{b_x-2\hat{a}} = \frac{b_y\left[\mathrm{tri_c}\left(\alpha_a+\Delta\alpha_1\right)+\mathrm{tri_s}\,\alpha_a\right]}{b_x\left[\mathrm{tri_c}\left(\alpha_a+\Delta\alpha_1\right)+\mathrm{tri_s}\,\alpha_a\,\mathrm{tri_s}\,\Delta\alpha_1\right]}, \\ \tan\gamma &= \frac{\hat{b}_y\left(b_x-2\hat{a}\right)}{b_x\left(b_x-\hat{a}\right)} = \frac{\left(b_x-2\hat{a}\right)\tan\beta}{b_x} = \frac{b_y\left[\mathrm{tri_c}\left(\alpha_a+\Delta\alpha_1\right)-\mathrm{tri_s}\,\alpha_a\right]}{b_x\left[\mathrm{tri_c}\left(\alpha_a+\Delta\alpha_1\right)+\mathrm{tri_s}\,\alpha_a\,\mathrm{tri_s}\,\Delta\alpha_1\right]}. \end{aligned} \tag{45}$$

Using formulas (45), equations (43) can be expressed, if necessary, through the angles $\beta$ and/or $\gamma$. If we set the phase shift $\Delta\alpha_1$ equal to -1.0 ($\gamma$=41°), -0.73 ($\gamma$=25°), and -0.5 ($\gamma$=0°) then we can built by formulas (43) the same family of loops having the same parameters $a$, $b_x$, $b_y$ ($b_x$<2$a$, $\Delta\alpha_2$=$\Delta\alpha_3$=0) as in Fig. 19(a), where the loops have been built by formulas (40). Setting nonzero phase shifts $\Delta\alpha_1$, $\Delta\alpha_2$, $\Delta\alpha_3$, we can build the piecewise-linear loops shown in Figs. 20(b), (d)-(f), (k), (l) by equations (43).

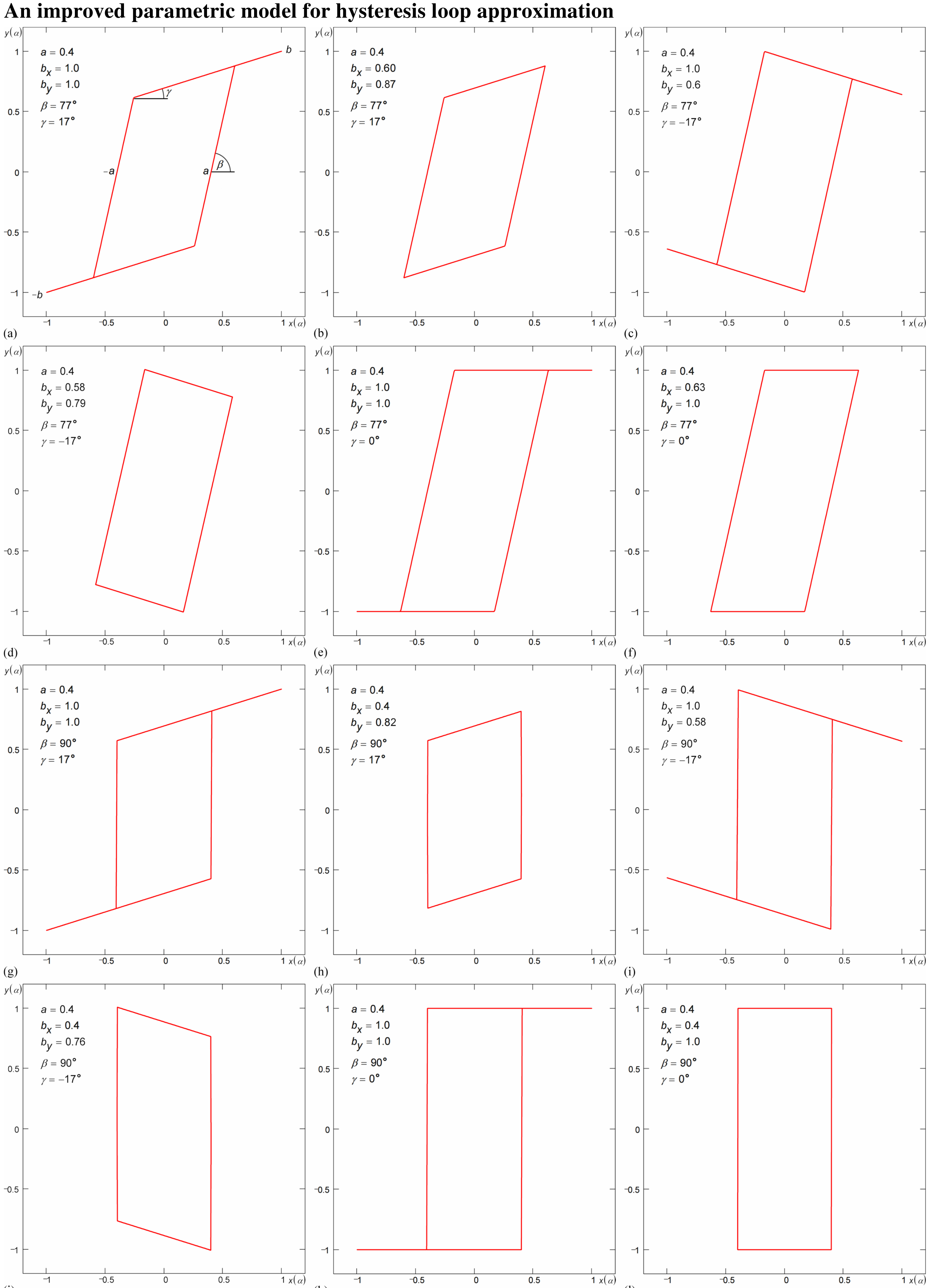

Fig. 20. Piecewise-linear hysteresis loops: (a) Play with Gain, (b) Play with Gain without Whiskers (parallelogram loop), (c) Play with Attenuation, (d) Play with Attenuation without Whiskers, (e) Play (backlash), (f) Play without Whiskers, (g) Non-ideal Relay with Gain, (h) Non-ideal Relay with Gain without Whiskers, (i) Non-ideal Relay with Attenuation, (j) Non-ideal Relay with Attenuation without Whiskers, (k) Non-ideal Relay (Schmitt trigger), (l) Non-ideal Relay without Whiskers (rectangular loop) built on triangular pulses.

Moving the split parameter $a$ in formulas (43) into the argument of the triangular pulse function, a system of equations can be obtained that describe the piecewise-linear hysteresis loop Play with Gain ($\Delta\alpha_1$=$\Delta\alpha_2$=$\Delta\alpha_3$=0)[6]

$$x(\alpha)=b_x \operatorname{tri_s} \alpha,$$
$$y(\alpha)=\left(b_y-b_x\tan\gamma\right)\left[\operatorname{tri_s}\left(\alpha-\frac{\alpha_a\tan\beta}{\tan\beta-\tan\gamma}+\frac{T}{8}\right)-\operatorname{tri_c}\left(\alpha-\frac{\alpha_a\tan\beta}{\tan\beta-\tan\gamma}+\frac{T}{8}\right)\right]+b_x\tan\gamma \operatorname{tri_s} \alpha. \qquad (46)$$

The loop parameters $b_x$, $b_y$, $\beta$, $\gamma$ relate as follows: $2b_y/b_x=\tan\beta+\tan\gamma$. This relationship makes equations (46) much less universal than equations (47) or (54) presented below. Nevertheless, equations (46) can be used for building loops similar to the ones shown in Figs. 20(a)-(f).

***a. The universal formula of a piecewise-linear loop***

In the present study, the following general expression is obtained, which describe a piecewise-linear hysteresis loop Play with Gain[14] (see Fig. 20(a))[6]

$$x(\alpha)=b_x \operatorname{tri_s} \alpha,$$
$$y(\alpha)=\left(b_y-b_x\tan\gamma\right)\operatorname{trp_s}\left(\alpha-\frac{\alpha_a\tan\beta}{\tan\beta-\tan\gamma}\right)+b_x\tan\gamma \operatorname{tri_s} \alpha, \qquad (47)$$

in which the upper base $d$ of the trapezoidal pulses $\operatorname{trp_s}$ is defined according to the formula

$$d=\frac{T\left(b_x\tan\beta-b_y\right)}{2b_x\left(\tan\beta-\tan\gamma\right)}, \qquad (48)$$

and the lower base $D$ according to the formula $D$=$T$–$d$. Step-by-step derivation of formulas (47) is given in the supplementary material.

The piecewise-linear loop (47) is built using a combination of triangular and trapezoidal pulses. If necessary, the trapezoidal pulses $\operatorname{trp_s}$ in (47) can be represented as a sum of two triangular pulses[6]

$$\operatorname{trp_s}(\alpha)=\frac{\operatorname{tri_s}\left(\alpha+\frac{d}{2}\right)+\operatorname{tri_s}\left(\alpha-\frac{d}{2}\right)}{\operatorname{tri_c}\left(\frac{\alpha_a\tan\beta}{\tan\beta-\tan\gamma}+\frac{d}{2}\right)+\operatorname{tri_c}\left(\frac{\alpha_a\tan\beta}{\tan\beta-\tan\gamma}-\frac{d}{2}\right)}. \qquad (49)$$

Equations (47) are universal since they allow to obtain the whole set of piecewise-linear loops of Play and Relay shown in Fig. 20. For example, letting $\gamma$=0° in (47), (48), we obtain a system of equations for building a hysteresis loop of Play (backlash) type[14, 15] (see Fig. 20(e))

$$x(\alpha)=b_x \operatorname{tri_s} \alpha,$$
$$y(\alpha)=b_y \operatorname{trp_s}(\alpha-\alpha_a), \qquad (50)$$

in which the upper base $d$ of the trapezoidal pulses $\operatorname{trp_s}$ is defined by the formula

$$d=\frac{T\left(b_x\tan\beta-b_y\right)}{2b_x\tan\beta}. \qquad (51)$$

Assuming $\beta$=90° in (47), (48), one can obtain the system of equations for building a Non-ideal Relay with Gain[16] hysteresis loop (see Fig. 20(g))

$$x(\alpha)=b_x \operatorname{tri}_s \alpha,$$
$$y(\alpha)=\left(b_y - b_x \tan\gamma\right)\operatorname{rect}_s\left(\alpha-\alpha_a\right)+b_x \tan\gamma \operatorname{tri}_s \alpha, \tag{52}$$

in which the upper base $d$ and the lower base $D$ of the trapezoidal pulses are the same and equal to $T/2$. The latter means degeneration of the trapezoidal pulses in the generating function $y(\alpha)$ into the rectangular pulses $\operatorname{rect}_s$ with a 50% duty cycle.

Assuming $\gamma=0°$ in (52), one can obtain the system of equations for building a Non-ideal Relay hysteresis loop (see Fig. 20(k))

$$x(\alpha)=b_x \operatorname{tri}_s \alpha,$$
$$y(\alpha)=b_y \operatorname{rect}_s\left(\alpha-\alpha_a\right). \tag{53}$$

To obtain loops without whiskers (see Figs. 20(b), (d), (f), (h), (j), (l)), one of the parameters $a$, $b_x$, $b_y$, $\beta$ in formulas (47), (50), (52), (53) should be expressed through the remaining parameters according to the relation $\tan\beta=b_y/(b_x-a)$. To obtain loops with attenuation[17, 18] (see Figs. 20(c), (d), (i), (j)), the negative angle $\gamma$ should be set in formulas (47), (52).

Fig. 21 shows examples of decomposition of piecewise-linear loops Play with Gain ($b_{1x}=a_1b_x/a$, $b_{1y}=a_1b_y/a$, $\beta_1=\beta$, $\gamma_1=\gamma$, $a_2=a-a_1$, $b_{2x}=b_x-b_{1x}$, $b_{2y}=b_y-b_{1y}$, $\beta_2=\beta$, $\gamma_2=\gamma$), Play ($b_{1x}=a_1b_x/a$, $\tan\beta_1=Tb_{1y}/[(T-2d)b_{1x}]$, $\gamma_1=0$, $a_2=a-a_1$, $b_{2x}=b_x-b_{1x}$, $b_{2y}=b_y-b_{1y}$, $\tan\beta_2=Tb_{2y}/[(T-2d)b_{2x}]$, $\gamma_2=0$), Non-ideal Relay with Gain ($b_{1x}=a_1b_x/a$, $\beta_1=90°$, $\gamma_1=\gamma$, $a_2=a-a_1$, $b_{2x}=b_x-b_{1x}$, $b_{2y}=b_y-b_{1y}$, $\beta_2=90°$, $\gamma_2=\gamma$) and Non-ideal Relay ($b_{1x}=a_1b_x/a$, $\beta_1=90°$, $\gamma_1=0$, $a_2=a-a_1$, $b_{2x}=b_x-b_{1x}$, $b_{2y}=b_y-b_{1y}$, $\beta_2=90°$, $\gamma_2=0$) into two loops of the same type as the loop under decomposition.[6] Since each of the loops of the decomposition, in turn, can be decomposed in a similar way then, like the smooth loops, piecewise linear loops can be represented in general case as an infinite sum of loops.

Fig. 22 shows examples of building more complex piecewise-linear loops: Play-Play-Play with Gain[14] ($p=2.4$, $b_{1x}=b_x/p$, $b_{1y}=b_y/p$, $\beta_1=\beta$, $\gamma_1=\gamma$, $a_2=a-a_1$, $b_{2x}=b_x-b_{1x}$, $b_{2y}=b_y-b_{1y}$, $\beta_2=\beta$, $\gamma_2=\gamma$, where $p$ is an arbitrary number in general case), Play-Play-Play ($p=2.4$, $b_{1x}=b_x/p$, $b_{1y}=b_y/p$, $\beta_1=\beta$, $\gamma_1=0$, $a_2=a-a_1$, $b_{2x}=b_x-b_{1x}$, $b_{2y}=b_y-b_{1y}$, $\beta_2=\beta$, $\gamma_2=0$), Play-Relay-Play with Gain[19] ($\beta_1=\beta$, $\gamma_1=0$, $a_2=a-a_1$, $b_{2x}=b_x-b_{1x}$, $b_{2y}=b_y-b_{1y}$; $a_1$, $\beta_2$ and $\gamma_2$ are determined by numerical solution of a system of

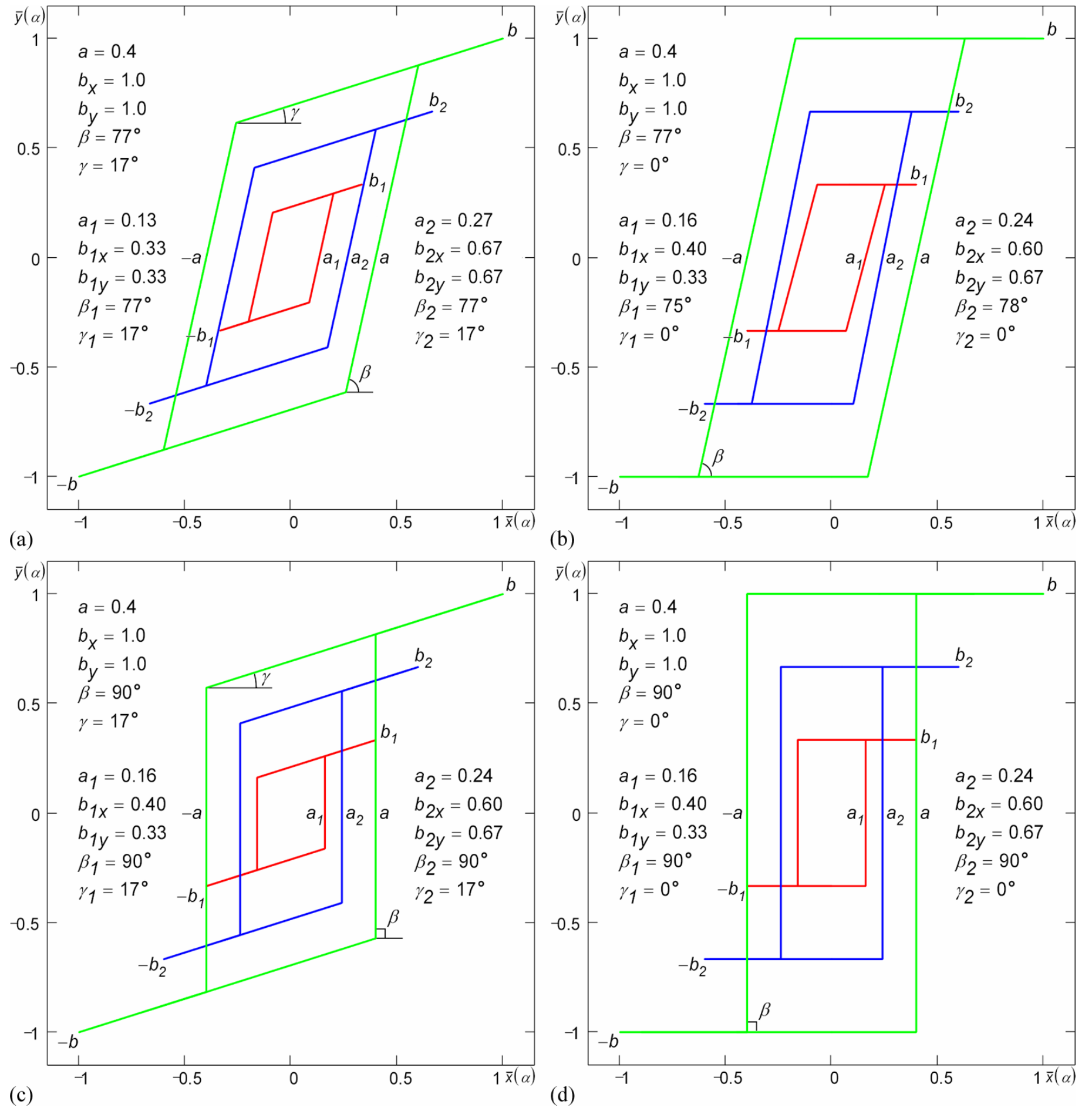

Fig. 21. Decomposition of the piecewise-linear hysteresis loop (green) (a) Play with Gain, (b) Play, (c) Non-ideal Relay with Gain, (d) Non-ideal Relay into two loops of the same type as the loop under decomposition.

equations[6]), Play-Relay-Play ($\beta_1$=$\beta$, $\gamma_1$=0, $a_2$=$a$-$a_1$, $b_{2x}$=$b_x$-$b_{1x}$, $b_{2y}$=$b_y$-$b_{1y}$, $\gamma_2$=0; $a_1$ and $\beta_2$ are determined by numerical solution of a system of equations[6]). The loops of the Play-Play-Play type are obtained by combining a pair of loops of Play type; the loops of the Play-Relay-Play type are obtained by combining loops of Non-ideal Relay type and Play type.[6] The details of composition/decomposition of the piecewise-linear hysteresis loops presented in Figs. 21 and 22 are given in the supplementary material, variations of these loops with attenuation and/or without whiskers can also be found there.

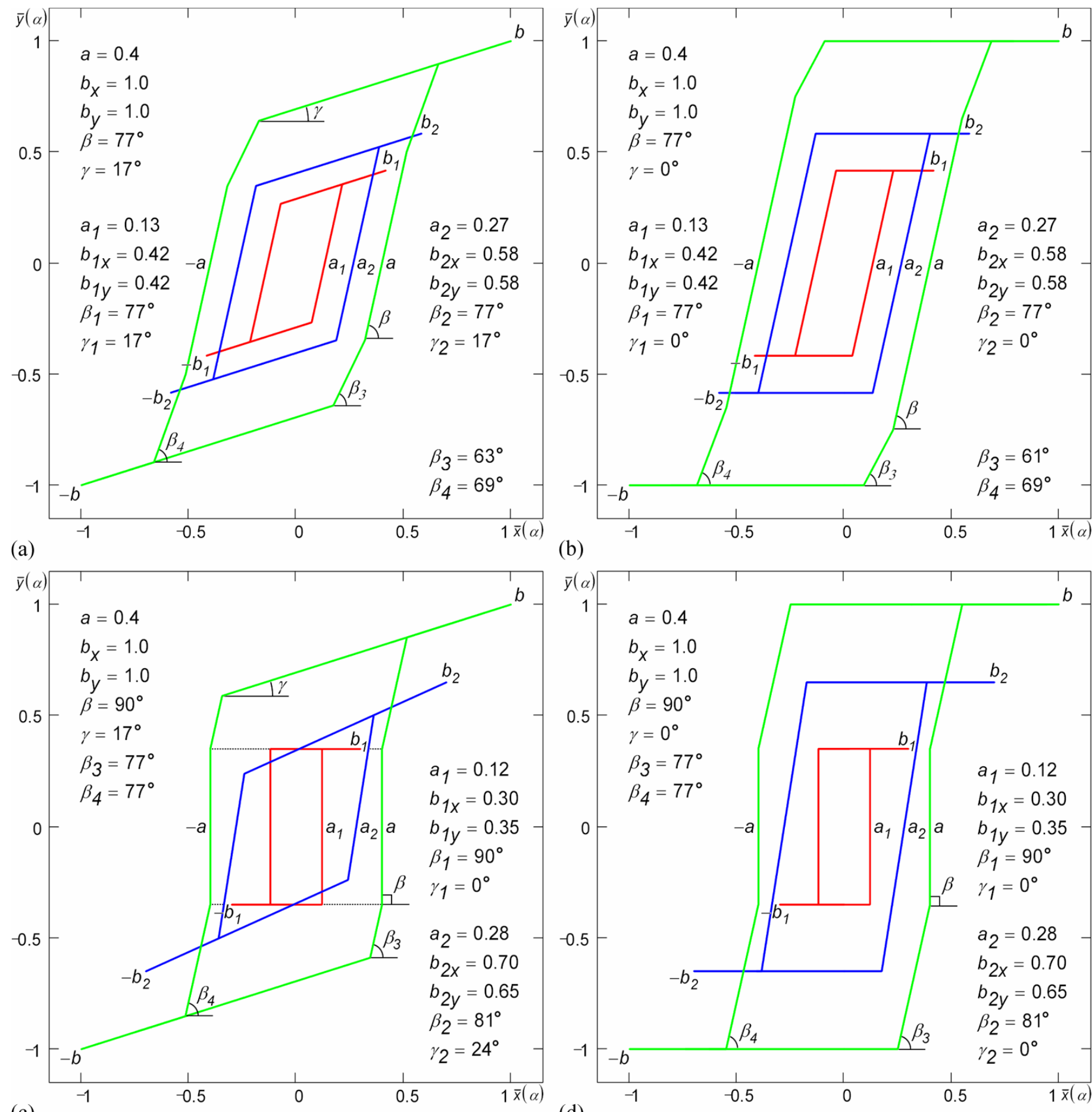

Fig. 22. Piecewise-linear hysteresis loops (a) Play-Play-Play with Gain, (b) Play-Play-Play, (c) Play-Relay-Play with Gain, (d) Play-Relay-Play obtained by addition of two loops of type (a), (b) Play, (c), (d) Relay and Play.

Detection of the loops shown in Fig. 22 during physical measurements may indicate that in the system under consideration there are two different hysteresis processes superimposed on each other.

***b. Generating function with a threshold element***

In addition to formulas (47), to describe the piecewise-linear hysteresis loop Play with Gain (see Fig. 20(a)), one can apply equations, in which the combination of triangular pulses in the generating function $y(\alpha)$ is passed through a threshold element $\mathrm{H_r}$[6]

$$x(\alpha) = b_x \operatorname{tri_s} \alpha,$$

$$y(\alpha) = 2(b_y - b_x \tan\gamma)\left[\mathrm{H_r}\left(\frac{a\tan\beta \operatorname{tri_c}\alpha + b_x\tan\gamma - b_y}{(a - b_x)\tan\beta + b_x\tan\gamma} + \operatorname{tri_s}\alpha\right) - \frac{1}{2}\right] + b_x\tan\gamma \operatorname{tri_s}\alpha. \qquad (54)$$

The threshold element $\mathrm{H_r}$ is a real (non-ideal) unit step function defined as follows

$$\mathrm{H_r}(t) = \begin{cases} 0, & t < 0, \\ \dfrac{t}{t_f}, & 0 \le t \le t_f, \\ 1, & t > t_f, \end{cases} \qquad (55)$$

where $t_f$ is a "front duration" of the real step function. The parameter $t_f$ is determined by the loop parameters $a$, $b_x$, $b_y$, $\beta$, $\gamma$ according to the following formula[6]

$$t_f = \frac{2(b_x \tan\gamma - b_y)}{(a - b_x)\tan\beta + b_x\tan\gamma}, \qquad (56)$$

and, vice versa: any of the mentioned parameters can be determined by this formula by fixing the front duration $t_f$

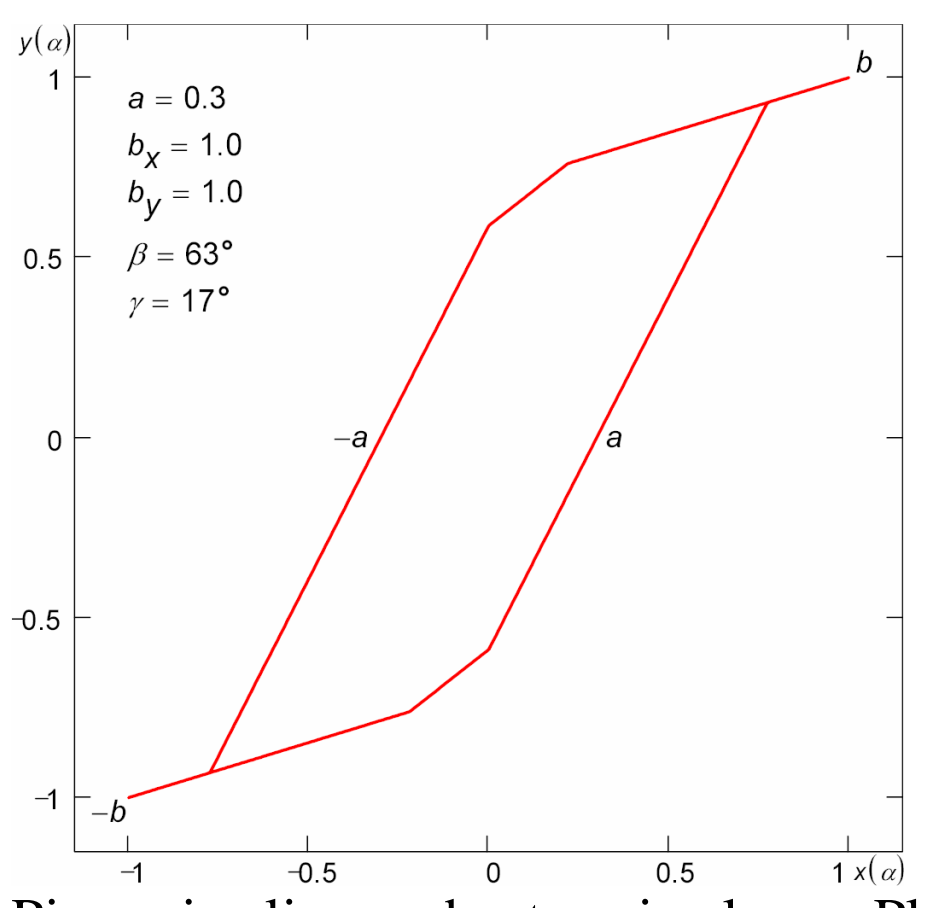

Fig. 23. Piecewise-linear hysteresis loop Play-Play with Gain. The additional corners occur when the splitting $a$ and the tilting $\beta$ become less of a certain value.

of the step function. With $t_f \to 0$ ($\beta \to 90°$), the function $\mathrm{H_r}$ degenerates into the ideal unit step function H (Heaviside function). In that case, Non-ideal Relay loops are formed (see Figs. 20(g)-(l)).

As a rule, equations (54) allow to obtain the whole set of types of the piecewise-linear loops of Play and Relay shown in Fig. 20. However, for some values of parameters $a$, $b_y$, $\beta$, $\gamma$, additional corners (see Fig. 23, cf. Fig. 20(a)) may appear turning the loop Play into the loop Play-Play[10].

If necessary, the formulas (40)-(47), (50), (52)-(54) can be expressed through trapezoidal pulses by using the following representation of triangular pulses as a sum of two trapezoidal pulses ($D$=3$d$, $T$=$d$+$D$)

$$\begin{aligned}
\mathrm{tri_s}(\alpha) &= \frac{1}{2}\mathrm{trp_s}\Big(\alpha + \frac{T}{8}\Big) + \frac{1}{2}\mathrm{trp_s}\Big(\alpha - \frac{T}{8}\Big),\\
\mathrm{tri_c}(\alpha) &= \frac{1}{2}\mathrm{trp_c}\Big(\alpha + \frac{T}{8}\Big) + \frac{1}{2}\mathrm{trp_c}\Big(\alpha - \frac{T}{8}\Big).
\end{aligned} \tag{57}$$

***3. Shifted loops***

In general, the shifted piecewise-linear loops (see Fig. 24) are obtained by adding two piecewise-linear loops (47). In such loops, the split parameters $a_1$, $a_2$ set the horizontal size $w$ of the shifted section; the saturation parameter $b_{1y}$ defines vertical position $h$ of this section; the gain $\gamma$ sets the slope of the section; and the angle $\beta_1$ ($\beta_2$=$\beta_1$) defines such a slope of the loops being added that the total slope $\beta$ of the resultant loop would be equal to the preset slope.

For the loop Shifted Play with Gain[19, 20] (see Fig. 24(a)), the splitting $a$ can be found by the formula

$$a = a_h - \frac{h}{\tan\beta} + \frac{w}{2}\left(\frac{\tan\gamma}{\tan\beta} - 1\right) \tag{58}$$

and the splitting $a_1$ – by the formula

$$a_1 = \frac{\left(a_h - \dfrac{w}{2}\right)(\tan\beta_1 - \tan\gamma) + b_x \tan\gamma - 2b_{1y}}{2\tan\beta_1}. \tag{59}$$

The other parameters $a_2$, $b_{1y}$, $\beta_1$ are determined numerically by solving the following system of equations ($b_{1x}$=$b_x$/2, $\gamma_1$=$\gamma$, $b_{2x}$=$b_x$-$b_{1x}$, $b_{2y}$=$b_y$-$b_{1y}$, $\beta_2$=$\beta_1$, $\gamma_2$=$\gamma$)

$$\begin{aligned}
&\bar{x}\left(\alpha_{e_2}(a_2) - \frac{T}{2}\right) - \bar{x}\left(\alpha_{d_1}\right) = w,\\
&\bar{y}\left(\alpha_{a_h}, b_{1y}\right) = h,\\
&\frac{\bar{y}\left(\alpha_{d_2}(\beta_1)\right) - \bar{y}\left(\alpha_{e_2}(\beta_1) - \dfrac{T}{2}\right)}{\bar{x}\left(\alpha_{d_2}(\beta_1)\right) - \bar{x}\left(\alpha_{e_2}(\beta_1) - \dfrac{T}{2}\right)} = \tan\beta,
\end{aligned} \tag{60}$$

where $\bar{x}(\alpha)$, $\bar{y}(\alpha)$ are generating functions of the shifted loop; $\alpha_{d_1}$, $\alpha_{a_h}$ are the values of parameter $\alpha$ in points $d_1$, $a_h$ of the shifted loop, respectively; $\alpha_{e_2}(a_2)$, $\alpha_{e_2}(\beta_1)$, $\alpha_{d_2}(\beta_1)$ are the dependences of parameter $\alpha$ in points $e_2$, $d_2$ of the shifted loop from variables $a_2$ and $\beta_1$, respectively; $\bar{y}(\alpha_{a_h}, b_{1y})$ is the dependence of the generating function on $b_{1y}$ in point $a_h$ ($\alpha=\alpha_{a_h}$).[6]

Figs. 24(b)-(d) shows examples of the shifted loops: Play,[14] Non-ideal Relay with Gain,[21] and Non-ideal Relay[22] which are special cases of the loop Shifted Play with Gain. The other variations of the shifted loops with attenuation and/or without whiskers are given in the supplementary material. The detection of the shifted loop during the physical measurements may indicate the simultaneous existence of two separate hysteresis processes superimposed on each other in the system under consideration (for example, there is a film of two layers which magnetic properties/states are different[20, 22]).

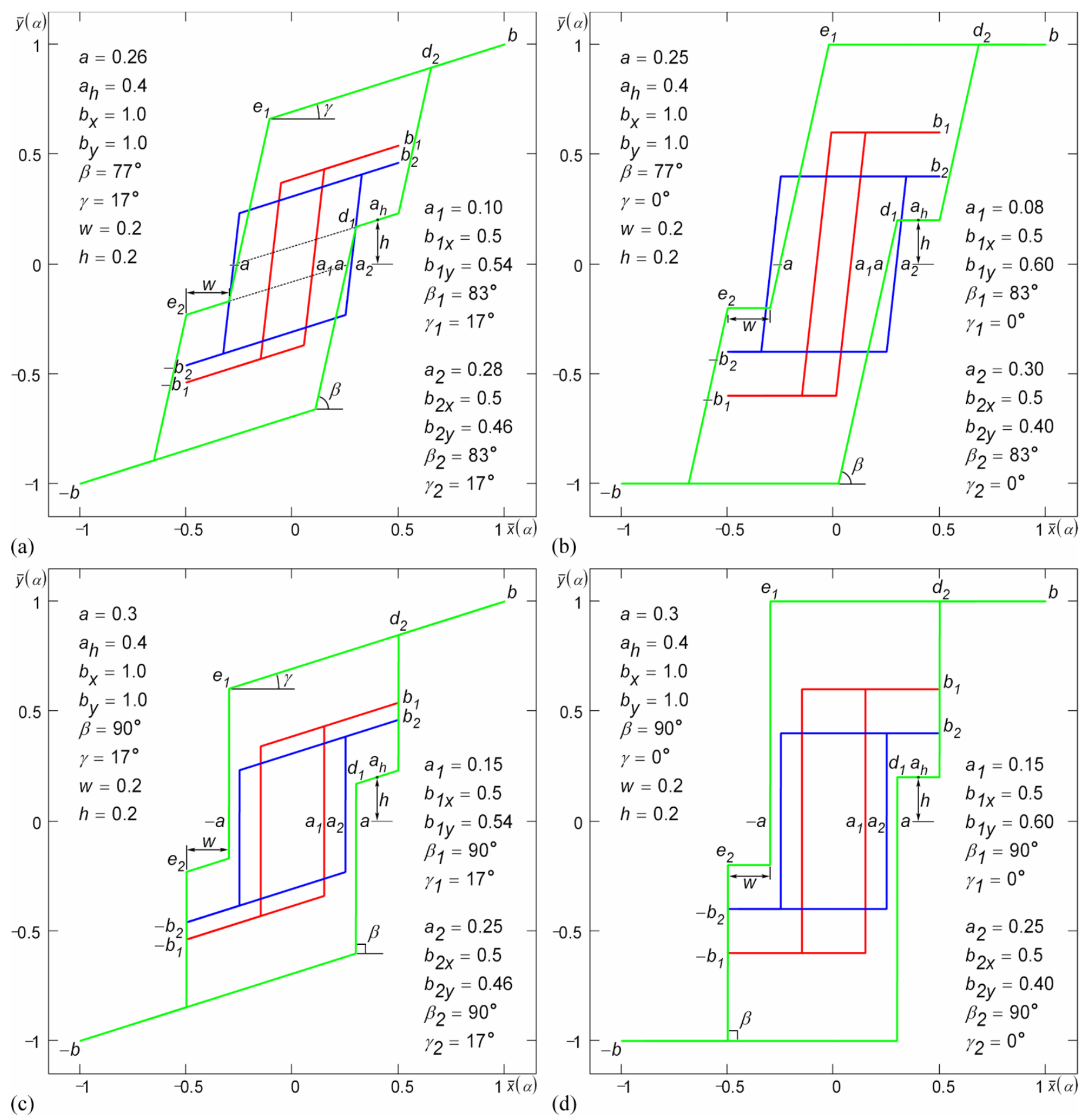

Fig. 24. Shifted piecewise-linear hysteresis loops (a) Play with Gain, (b) Play, (c) Non-ideal Relay with Gain, (d) Non-ideal Relay obtained by addition of two loops of type (a), (b) Play, (c), (d) Relay.

### *a. Three-level loops*

The three-level loops are the special case of the shifted loops (see Fig. 25). Three-level loops appear when negative splitting $a$ is set and the next equation is satisfied[6]

$$w = \frac{2a\tan\beta}{\tan\gamma - \tan\beta}. \tag{61}$$

The description of three-level loops with attenuation and/or with no whiskers is given in the supplementary material. The detection of the three-level loop during the physical measurements may serve as a sign of simultaneous action of two separate hysteresis processes superimposed on each other in the system under consideration.

## C. Double hysteresis loops

### *1. Linking loops at the saturation point*

To make the description of double hysteresis loops more accurate, the method suggested in Ref. 1 has been improved. In particular, according to the improved method, the equations for a double smooth loop non-self-crossing in the origin of coordinates (0-shaped loop) are as follows ($\alpha$=0…2π):

$$\bar{\bar{x}}(\alpha)=x\left(2\alpha-\frac{\pi}{2}\operatorname{sgn}(\pi-\alpha)-\Delta\alpha_3\right)+b_x\operatorname{sgn}(\pi-\alpha)=(2\operatorname{rect}\alpha-1)\left(x\left(2\alpha-\Delta\alpha_3-\frac{\pi}{2}\right)+b_x\right),$$
$$\bar{\bar{y}}(\alpha)=y\left(2\alpha-\frac{\pi}{2}\operatorname{sgn}(\pi-\alpha)-\Delta\alpha_3\right)+b_y\operatorname{sgn}(\pi-\alpha)=(2\operatorname{rect}\alpha-1)\left(y\left(2\alpha-\Delta\alpha_3-\frac{\pi}{2}\right)+b_y\right),\qquad(62)$$

where sgn$\alpha$=$\alpha$/|$\alpha$| is the signum function; rect$\alpha$=H($\alpha$)-H($\alpha$-π) is a π-wide rectangular pulse. Double loop (62) is formed by linking two loops in the saturation point $b$ where the generating function $y(\alpha)$ reaches its maximum value $b_y$. Fig. 26(a) shows an example of double smooth non-self-crossing loop[23] built by formulas (62).

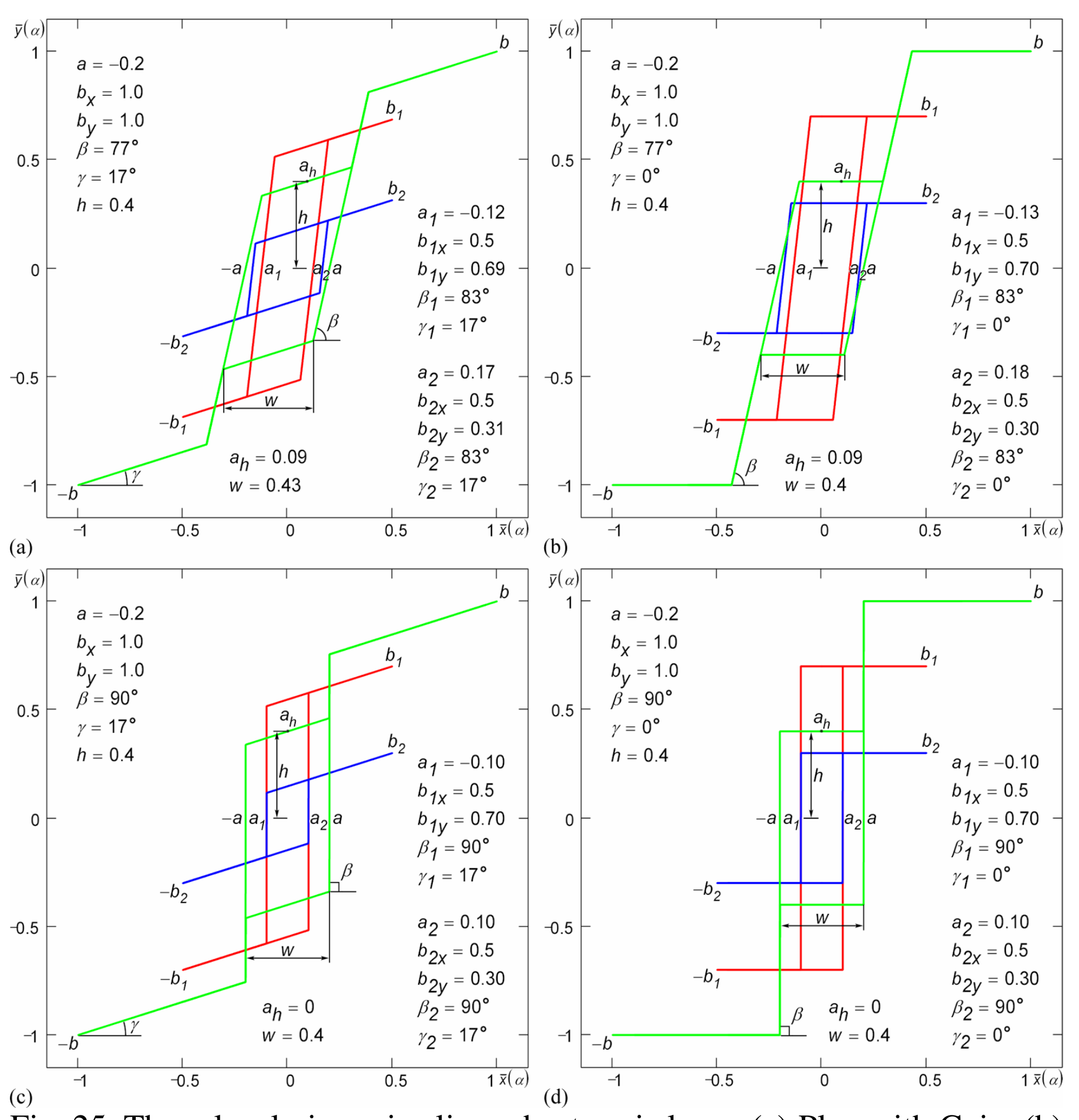

Fig. 25. Three-level piecewise-linear hysteresis loops (a) Play with Gain, (b) Play, (c) Non-ideal Relay with Gain, (d) Non-ideal Relay obtained by addition of two loops of type (a), (b) Play, (c), (d) Relay.

According to (62), movement along the double loop starts at the point (0, 0), goes on counterclockwise first by the top loop then by the bottom one and finishes at the point (0, 0). The double loop (62) differs from the one used before in that the starting point of the second loop coincides with the end point of the first one, and the end point of the second loop with the starting point of the first one. Thus, the double loop (62) can be drawn continuously, "without lifting the pencil from the paper", as the parameter $\alpha$ changes from 0 to 2π.

The loop with a self-crossing in the origin of coordinates (8-shaped loop) can be built according to the formulas ($\alpha$=0…2π):

$$\bar{\bar{x}}(\alpha)=x\left(\left(2\alpha-\frac{\pi}{2}\right)\operatorname{sgn}(\pi-\alpha)-\Delta\alpha_3\right)+b_x\operatorname{sgn}(\pi-\alpha)$$
$$=\operatorname{rect}\alpha\left(x\left(2\alpha-\Delta\alpha_3-\frac{\pi}{2}\right)+b_x\right)+(1-\operatorname{rect}\alpha)\left(x\left(\frac{\pi}{2}-\Delta\alpha_3-2\alpha\right)-b_x\right),$$
$$\bar{\bar{y}}(\alpha)=y\left(\left(2\alpha-\frac{\pi}{2}\right)\operatorname{sgn}(\pi-\alpha)-\Delta\alpha_3\right)+b_y\operatorname{sgn}(\pi-\alpha)$$
$$=\operatorname{rect}\alpha\left(y\left(2\alpha-\Delta\alpha_3-\frac{\pi}{2}\right)+b_y\right)+(1-\operatorname{rect}\alpha)\left(y\left(\frac{\pi}{2}-\Delta\alpha_3-2\alpha\right)-b_y\right).\qquad(63)$$

The loops (62) and (63) have no appearance differences.

Double non-self-crossing and self-crossing piecewise-linear loops are built by formulas (62) and (63), respec-tively, replacing $\pi$ with $T/2$ ($\alpha$=0…$T$). An example of double piecewise-linear non-self-crossing loop[16, 19, 24] built by formulas (62) with $\Delta\alpha_3$=0 is shown in Fig. 26(b). Self-crossing loop has the same appearance.

***2. Linking loops at the point $x_{max}$***

Besides linking two loops in the saturation point $b$, the loops can be linked at the point, where generating function $x(\alpha)$ reaches its maximum $x_{max}$. In order to determine the value of parameter $\alpha_{max}$ for which $x(\alpha_{max})=x_{max}$, the equation $dx(\alpha)/d\alpha=0\,|_{\alpha=\alpha_{max}}$ should be numerically solved. In the expanded form it can be written as follows

$$m\hat{a}\sin(\alpha_{max}+\Delta\alpha_1)\cos^{m-1}(\alpha_{max}+\Delta\alpha_1)-n\hat{b}_x\cos(\alpha_{max}+\Delta\alpha_2)\sin^{n-1}(\alpha_{max}+\Delta\alpha_2)=0. \qquad (64)$$

Double non-self-crossing loops linking in the $x_{max}$ point are built according to the formulas

$$\begin{aligned}\bar{\bar{x}}(\alpha)&=x\left(2\alpha-\frac{\pi}{2}(\operatorname{sgn}(\pi-\alpha)+1)+\alpha_{max}\right)+x(\alpha_{max})\operatorname{sgn}(\pi-\alpha)=(1-2\operatorname{rect}\alpha)\big(x(2\alpha+\alpha_{max})-x(\alpha_{max})\big),\\ \bar{\bar{y}}(\alpha)&=y\left(2\alpha-\frac{\pi}{2}(\operatorname{sgn}(\pi-\alpha)+1)+\alpha_{max}\right)+y(\alpha_{max})\operatorname{sgn}(\pi-\alpha)=(1-2\operatorname{rect}\alpha)\big(y(2\alpha+\alpha_{max})-y(\alpha_{max})\big).\end{aligned} \qquad (65)$$

Double self-crossing loops linking in the $x_{max}$ point are built according to the formulas

$$\begin{aligned}\bar{\bar{x}}(\alpha)&=x\left(\left(2\alpha-\frac{\pi}{2}\right)\operatorname{sgn}(\pi-\alpha)+\alpha_{max}-\frac{\pi}{2}\right)+x(\alpha_{max})\operatorname{sgn}(\pi-\alpha)\\ &=(1-\operatorname{rect}\alpha)\big(x(\alpha_{max}-2\alpha)-x(\alpha_{max})\big)-\operatorname{rect}\alpha\big(x(2\alpha+\alpha_{max})-x(\alpha_{max})\big),\\ \bar{\bar{y}}(\alpha)&=y\left(\left(2\alpha-\frac{\pi}{2}\right)\operatorname{sgn}(\pi-\alpha)+\alpha_{max}-\frac{\pi}{2}\right)+y(\alpha_{max})\operatorname{sgn}(\pi-\alpha)\\ &=(1-\operatorname{rect}\alpha)\big(y(\alpha_{max}-2\alpha)-y(\alpha_{max})\big)-\operatorname{rect}\alpha\big(y(2\alpha+\alpha_{max})-y(\alpha_{max})\big).\end{aligned} \qquad (66)$$

A double non-self-crossing loop built by formulas (65) is shown in Fig. 27. A double self-crossing loop built by formulas (66) has the same appearance.

***3. Replacement of the horizontal splitting with the vertical one***

The simplest way to get double non-self-crossing loops of the Propeller type[25] is based on representation (2), where instead of the splitting $a$ ($a\approx0$) along $x$, the splitting $a_y$ along $y$ is performed. The splitting curve equations in this case are as follows: $x_2(\alpha)=a\cos^m\alpha$, $y_2(\alpha)=a_y\cos^m\alpha$. From equation $y(\alpha_{d_{1,2}})=b_y/2$ written for a half-height of the saturation point $b_y$ (any desired value can be used instead of half-height), one can find the following dependencies ($m$=1)

$$\alpha_{d_{1,2}}(a_y)=\arccos\frac{b_y\left(a_y\pm\sqrt{4a_y^2+3b_y^2}\right)}{2\left(a_y^2+b_y^2\right)}. \qquad (67)$$

After that, by solving the equation numerically

$$x\big(\alpha_{d_1}(a_y)\big)-x\big(\alpha_{d_2}(a_y)\big)=2a_d, \qquad (68)$$

we determine the vertical split value $a_y$ such that the horizontal split value $a_d$ of the double loop at the half-height of the saturation point $b_y$ is equal to the specified value. Fig. 28 shows a double loop Propeller obtained according to the described method. The working expressions for the case $m$=3, 5, … are given in the supplementary material.

***4. Pinching a loop in the origin of coordinates***

Double non-self-crossing loops of Propeller type can be formed by "pinching" a loop in the origin of coordinates (see Fig. 29(a)). The pinching is achieved by raising the generating functions $x(\alpha)$, $y(\alpha)$ to odd powers $k$ and $l$ according to expressions (27). The parameter $\alpha_d$ can be found from equation $\bar{y}(\alpha_d)=b_y/2$ written for a half-height of the saturation point $b_y$ (any desired value can be used instead of half-height)

$$\alpha_d = \arcsin\frac{1}{\sqrt[l]{2Bb_y^{l-1}}}. \tag{69}$$

Solving the equation

$$\bar{x}(\alpha_d) - \bar{x}(\pi - \alpha_d) = 2a_d, \tag{70}$$

a split value $a$ can be determined such that the split value $a_d$ of the double loop at the half-height of the saturation point $b_y$ is equal to the specified value. The analytical solution of equation (70) for $k$=3 is as follows

$$a = \frac{\sqrt[3]{2A\left(a_d + \sqrt{4A^2 b_x^6 \sin^{6n}\alpha_d + a_d^2}\right)^2} - 2Ab_x^2 \sin^{2n}\alpha_d}{\cos^m \alpha_d \sqrt[3]{4A^2\left(a_d + \sqrt{4A^2 b_x^6 \sin^{6n}\alpha_d + a_d^2}\right)}}. \tag{71}$$

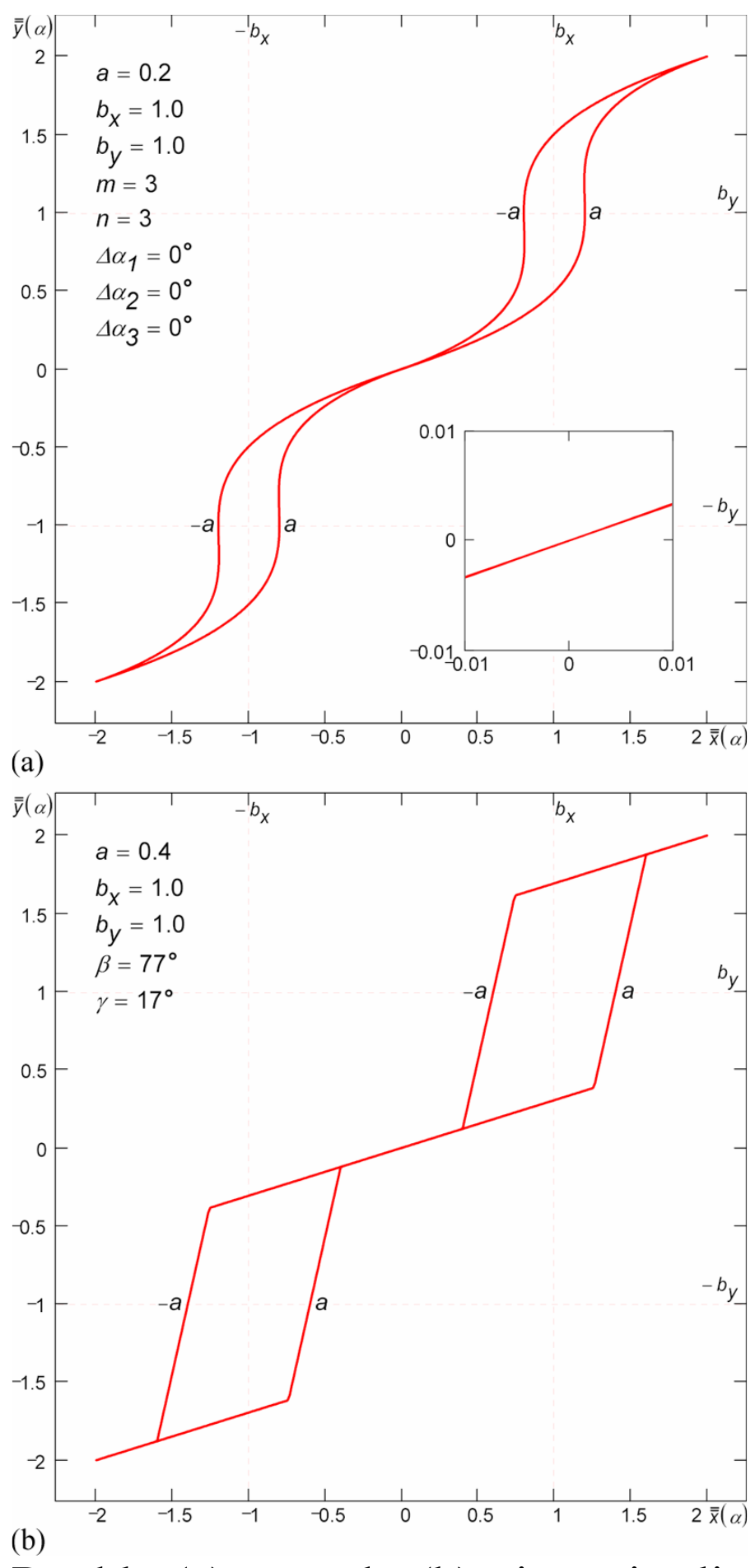

Fig. 26. Double (a) smooth, (b) piecewise-linear hysteresis loop formed by linking two loops (a) Classical, (b) Play with Gain in the saturation point $b$. In the linking point, the loop can be made both non-self-crossing and self-crossing.

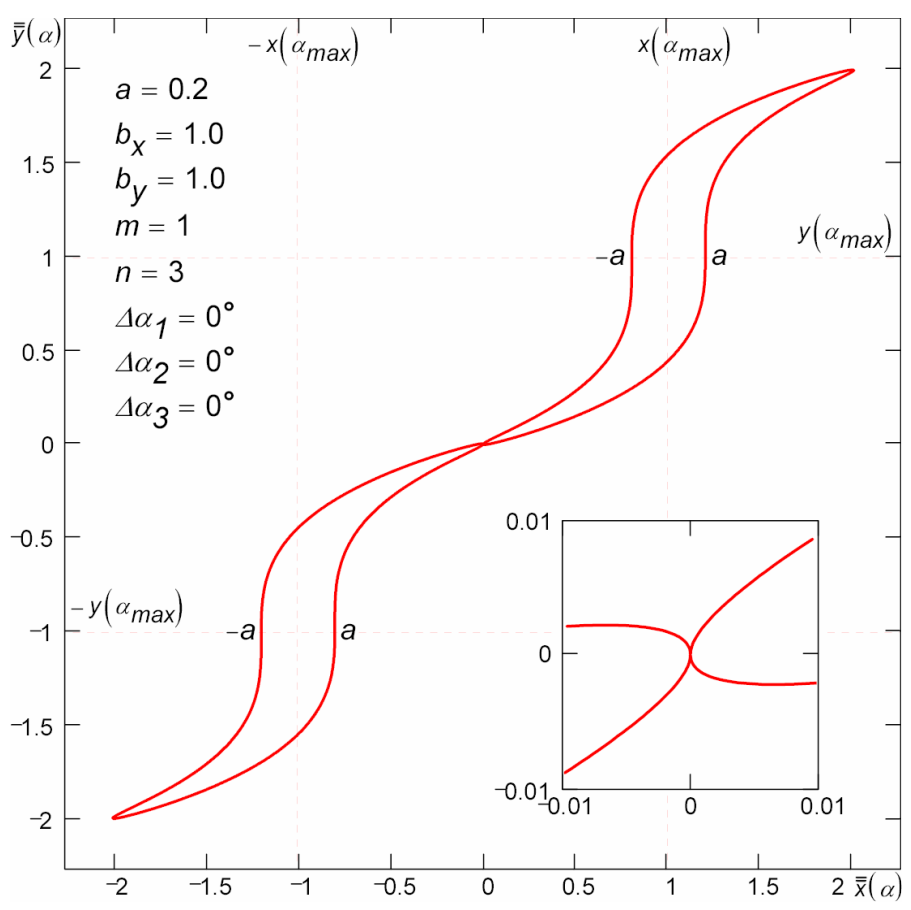

Fig. 27. Double smooth hysteresis loop formed by linking two loops Classical in the point $x_{max}$. In the linking point, the loop can be made both non-self-crossing and self-crossing.

There is another method of building double non-self-crossing loops by means of loop pinching. This method implies setting the splitting $a$ to zero in formula (13) and the phase shift $\Delta\alpha_2$ (or $\Delta\alpha_3$) to a non-zero value. With such loops, the shift $\Delta\alpha_2$ (or $\Delta\alpha_3$) serves as a splitting $a$. Fig. 29(b) shows an example of a double non-self-crossing loop of Propeller type[25, 26] obtained by pinching. The parameter $\alpha_d$=π/6 can be found from equation $y(\alpha_d)$=$b_y$/2 written for a half-height of the saturation point $b_y$ (any desired value can be used instead of half-height). Solving numerically the equation

$$x(\alpha_d) - x(\pi - \alpha_d) = 2a_d, \tag{72}$$

it is possible to determine such a value of the phase shift $\Delta\alpha_2$ (or $\Delta\alpha_3$) that the split value $a_d$ of the double loop at the half-height of the saturation point $b_y$ is equal to the specified value.

***5. Double piecewise-linear loop as a special case of the piecewise-linear shifted loop***

By setting splitting $a$ to 0, the shifted loop (see Sec. II.B.3) degenerates into the double non-self-crossing loop

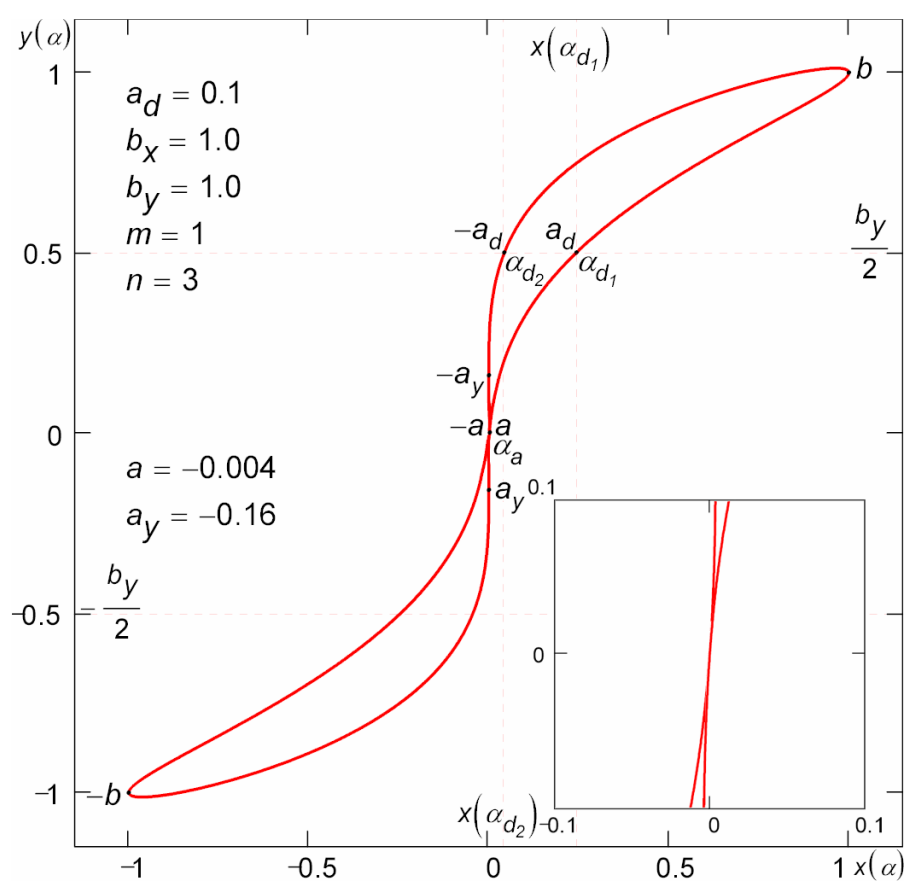

Fig. 28. Double non-self-crossing hysteresis loop of Propeller type formed by replacement of the horizontal splitting with the vertical one.

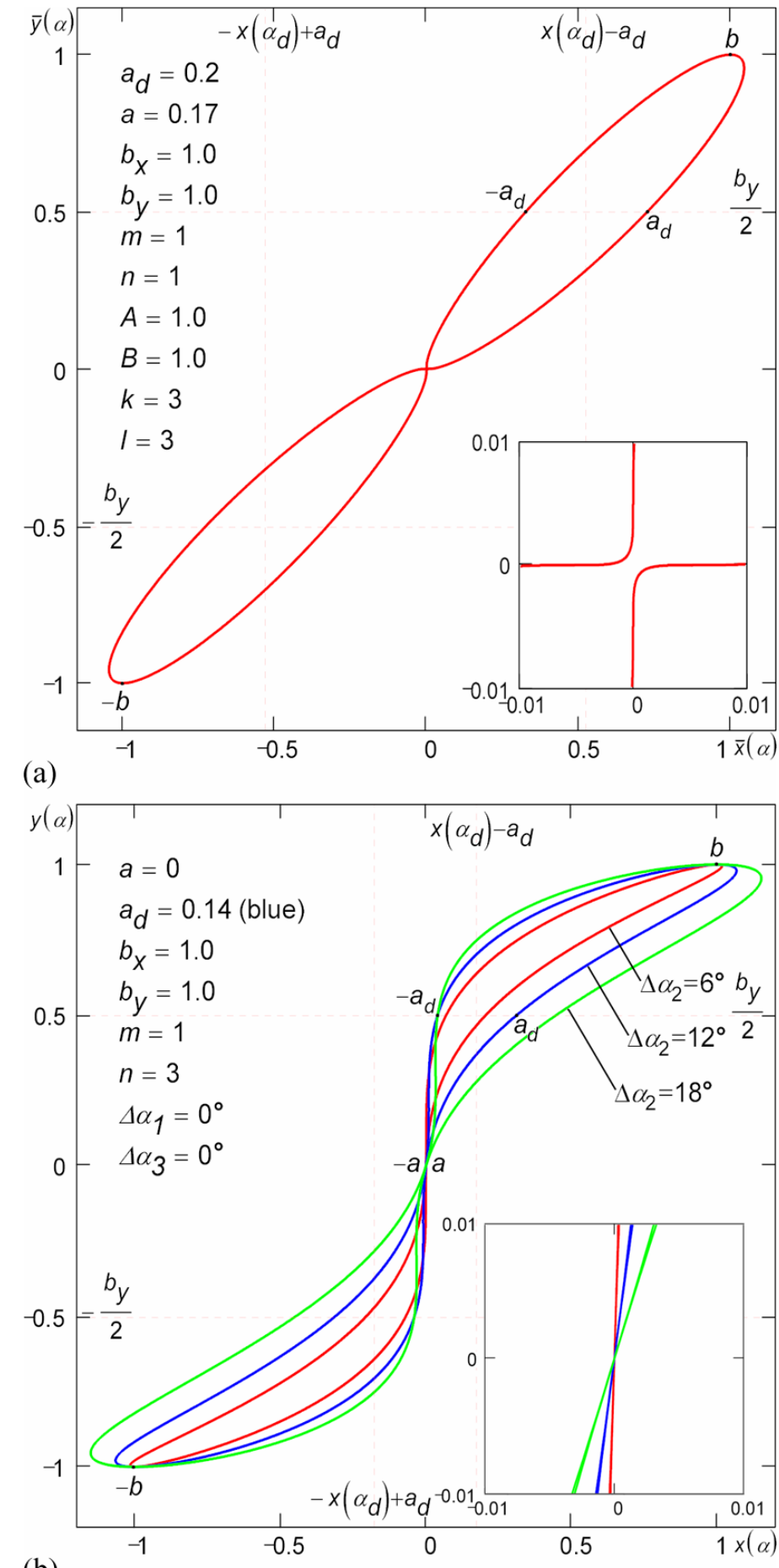

Fig. 29. Double non-self-crossing hysteresis loops of Propeller type formed by pinching a loop at the point of coordinate origin by means of (a) exponentiation, (b) zero splitting $a$ and a phase shift $\Delta\alpha_2$.

Play with Gain[27] (see Fig. 30(a)). Beside zero splitting, setting the height of the shifted section according to the formula $h=w\tan\gamma/2$ allows to obtain the double non-self-crossing loop Play with Gain with no bridge connection[19, 24] (see Fig. 30(b)). Double non-self-crossing loops Play, Relay with Gain and Relay having bridge connection and having no bridge connection as well as variations of these loops with attenuation and/or without whiskers are given in the supplementary material. Due to the way of formation of this double loop, it follows that detection of such loop during physical measurements may serve as a sign of simultaneous action of two separate hysteresis processes superimposed on each other in the system under consideration.

## D. Triple hysteresis loops

### *1. Linking loops at the saturation points*

A triple loop can be assembled of three loops – one central loop and two outside loops. Triple smooth loops linked at the saturation points $b$ are built according to the formulas ($\alpha$=0…2π)

$$\overset{\equiv}{x}(\alpha) = (\operatorname{rect}\alpha + \operatorname{rect}(\alpha-\pi))x_1\left(3\alpha-\frac{\pi}{2}\right) + (\operatorname{rect}\alpha+\operatorname{rect}(\alpha-\pi)-1)\left[x_2\left(\pm 3\alpha-\frac{\pi}{2}\right) - (b_{1x}+b_{2x})\operatorname{sgn}(\pi-\alpha)\right],$$
$$\overset{\equiv}{y}(\alpha) = (\operatorname{rect}\alpha + \operatorname{rect}(\alpha-\pi))y_1\left(3\alpha-\frac{\pi}{2}\right) + (\operatorname{rect}\alpha+\operatorname{rect}(\alpha-\pi)-1)\left[y_2\left(\pm 3\alpha-\frac{\pi}{2}\right) - (b_{1y}+b_{2y})\operatorname{sgn}(\pi-\alpha)\right], \quad (73)$$

where $x_1(\alpha)$, $y_1(\alpha)$ are equations of the central loop; $x_2(\alpha)$, $y_2(\alpha)$ are equations of the outside loops; $b_{1x}$, $b_{1y}$ are coordinates of the saturation point of the central loop; $b_{2x}$, $b_{2y}$ are coordinates of the saturation points of the outside loops; rect$\alpha$=H($\alpha$)-H($\alpha$-π/3) is a π/3-wide rectangular pulse.

When assembling loops (73), the condition $\gamma_1=\gamma_2$ is usually met, where $\gamma_1$, $\gamma_2$ are slope angles of tangents to the unsplit ($a_1$=0) central loop and the unsplit ($a_2$=0) outside loops, respectively, at the saturation point $b_1$. The slope angles $\gamma_1$, $\gamma_2$ of the tangents are defined by the formulas[6]

$$\gamma_1 = \arctan\frac{b_{1y}}{n_1 b_{1x}},$$
$$\gamma_2 = \arctan\frac{b_{2y}}{n_2 b_{2x}}. \qquad (74)$$

In case the argument $3\alpha$ of functions $x_2(\alpha)$, $y_2(\alpha)$ in equations (73) is used with the plus sign, a non-self-crossing loop is obtained. In case of the minus sign – a self-crossing loop. Non-self-crossing and self-crossing loops have the same appearance. The central and the outside loops may differ in type. Fig. 31(a) shows a triple smooth non-self-crossing loop built by equations (73). The loop consists of the central Classical loop (1) ($a_1$=$a$) and a pair of outside Classical loops (1) ($b_{2x}$=($b_x$-$b_{1x}$)/2, $b_{2y}$=$n_2b_{2x}$tan$\gamma_1$, $\gamma_2$=$\gamma_1$). The formulas for triple loops linked in the $x_{max}$ points are presented in the supplementary material.

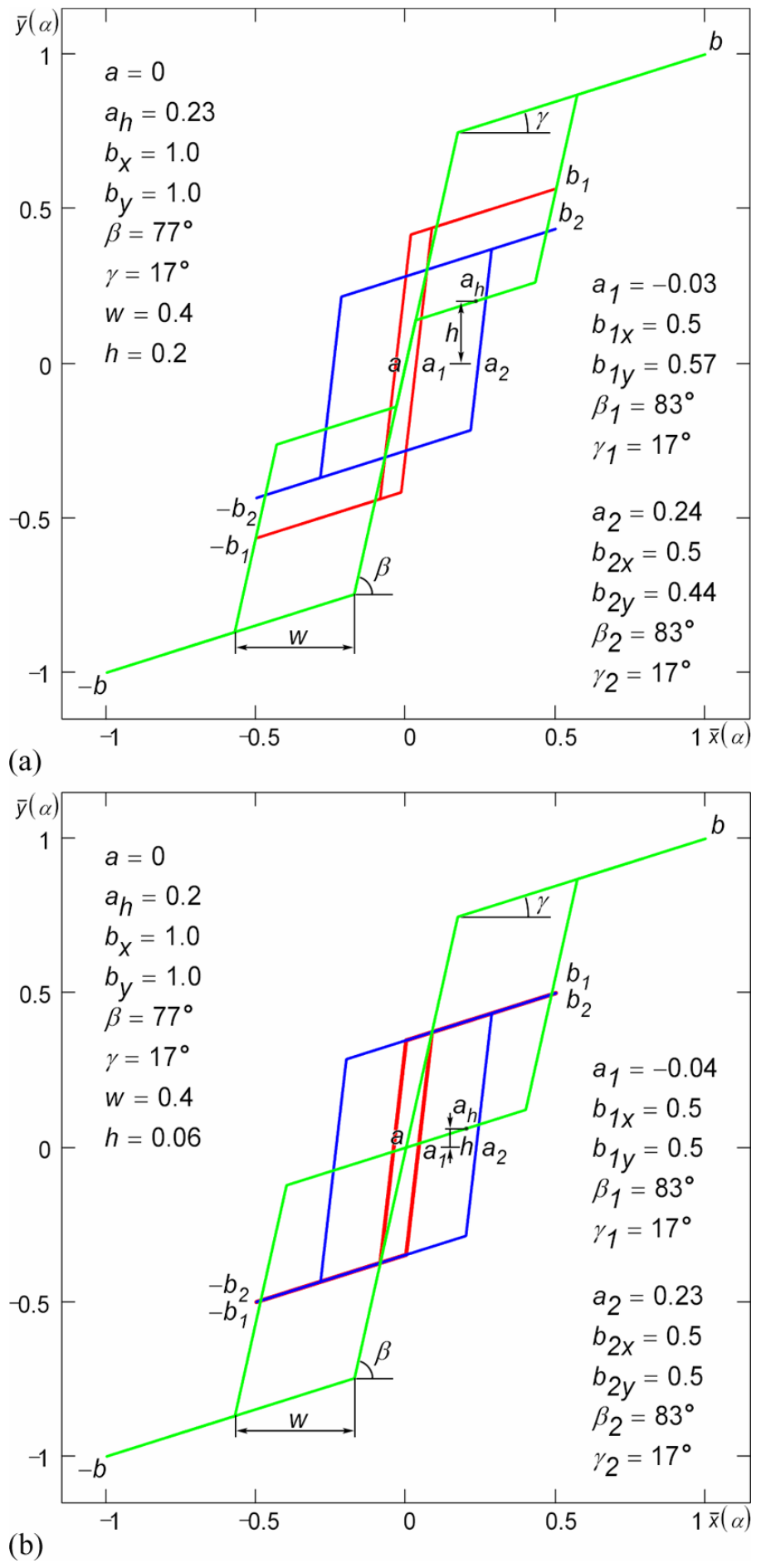

Fig. 30. Double piecewise-linear non-self-crossing hysteresis loops Play with Gain (a) bridge-connected, (b) no bridge-connected. The loops are obtained as a result of addition of two single loops Play with Gain.

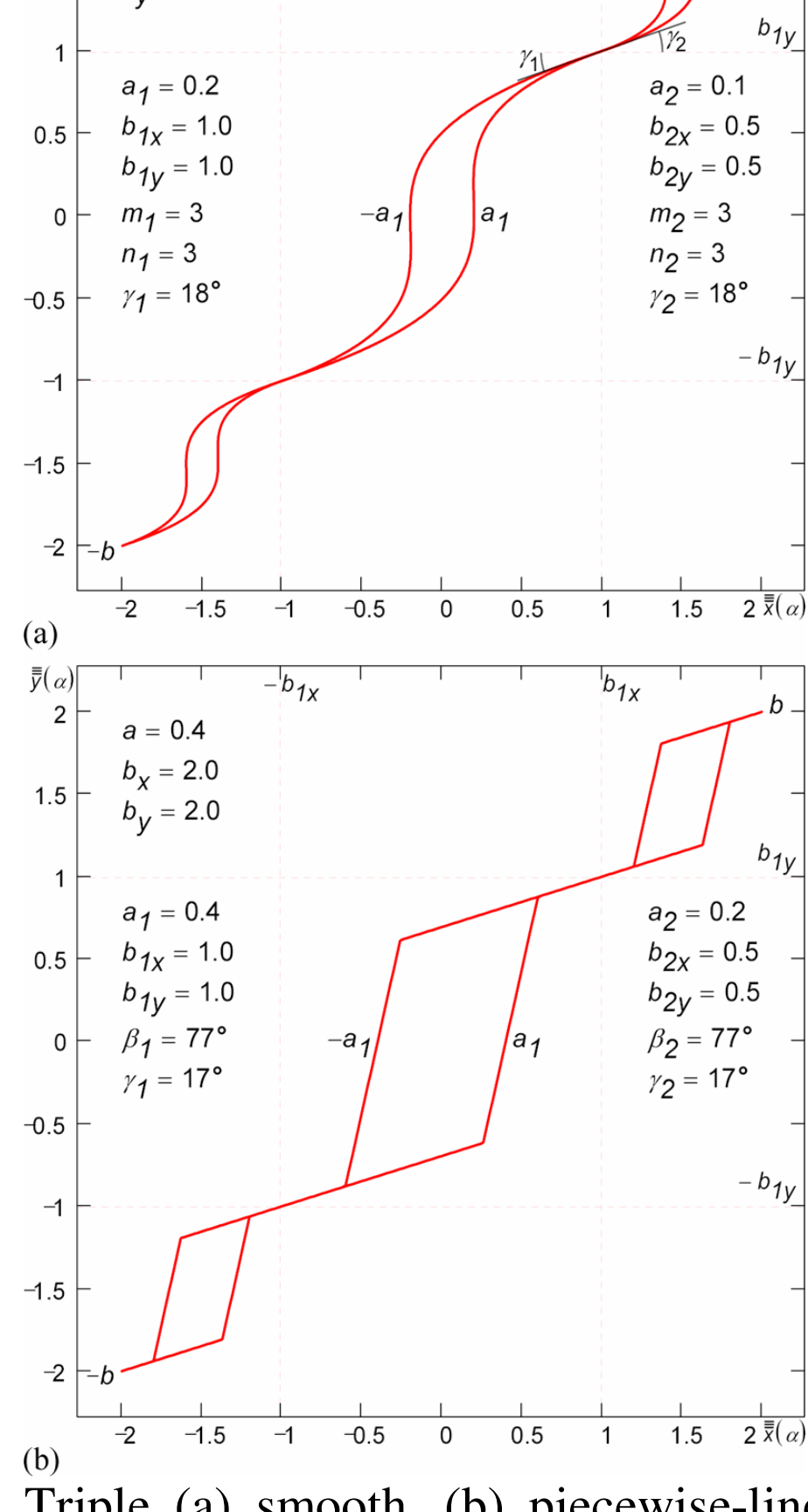

Fig. 31. Triple (a) smooth, (b) piecewise-linear hysteresis loop formed by linking three loops (a) Classical, (b) Play with Gain in the saturation points $b$. In the linking points, the loop can be made both non-self-crossing and self-crossing.

Triple non-self-crossing and self-crossing piecewise-linear loops[28] are built by formulas (73) by replacing π with $T$/2 ($\alpha$=0…$T$). Fig. 31(b) shows a triple piecewise-linear non-self-crossing loop built by equations (73). The loop consists of the central Play with Gain loop (47) ($a_1$=$a$) and a pair of outside Play with Gain loops (47) ($b_{2x}$=($b_x$-$b_{1x}$)/2, $b_{2y}$=($b_y$-$b_{1y}$)/2, $\gamma_2$=$\gamma_1$).

***2. Loops with arbitrarily long whiskers***

Triple hysteresis loops (73) are useful for producing single smooth[25, 29] or single hybrid[15] hysteresis loops having long whiskers. As a reminder, long whiskers in the model (1) can be obtained by increasing the power $m$.[6]

However, the loop curvature is also changing considerably at the same time. Fig. 32(a) shows a simulation of a single smooth loop Tilted Classical with long whiskers by using a triple non-self-crossing smooth loop (73). The triple loop consists of the central Tilted Classical loop (30) ($a_1$=$a$, $\theta_1$=$\theta$) and whiskers which are formed from a pair of the outside unsplit loops of Leaf type (1) ($a_2$=0, $b_{2x}$=($b_x$-$b_{1x}$)/2, $b_{2y}$=($b_y$-$b_{1y}$)/2, $m_2$=3, $n_2$=1) oriented at the angle $\gamma_2$=arctan($b_{2y}$/$b_{2x}$).

The required curvature $\kappa_1$ of the loop is determined by the formula

$$\kappa_1 = \arctan\frac{b_{1y}\left[\tan\theta_1(1-n_1)\left(b_y-b_{1y}\right)-b_x+b_{1x}\right]+n_1 b_{1x}\left(b_y-b_{1y}\right)}{b_{1x}\left(n_1-1\right)\left[\tan\theta_1\left(b_y-b_{1y}\right)-b_x+b_{1x}\right]}. \tag{75}$$

The formula (75) is derived from the condition $\gamma_1$=$\gamma_2$, where $\gamma_1$, $\gamma_2$ are slopes of the tangents at the unsplit ($a_1$=0) central loop and the unsplit ($a_2$=0) outside loops, respectively, at the saturation point $b_1$. From formula (75), for some arbitrary value $\kappa_1$ one can determine the corresponding value $\theta_1$ or for arbitrary values $b_x$, $\kappa_1$, $\theta_1$ – the corresponding value $b_y$.

Fig. 32(b) shows an example of building hybrid tilted Classical loop with whiskers[7, 15] out of three loops. The target loop consists of the central hybrid tilted Classical whiskerless loop (38) ($a_1$=$a$, $b_{1y}$=$b_y$-2$b_{2y}$, $\theta_1$=$\theta$, $\gamma_1$=$\gamma$, $\kappa_1$=$\kappa$) and whiskers which are formed from a pair of outside unsplit piecewise-linear loops Play with Gain without Whiskers (40) ($a_2$=0, $b_{2x}$=($b_x$-$b_{1x}$)/2, $b_{2y}$=$b_{2x}$tan$\gamma_2$, $m_2$=$n_2$=1, $\gamma_2$=$\gamma_1$).[6] Due to the relocation of the saturation point $b$ from position $\alpha$=$T$/4 to position $\alpha$=$T$/8 ($\Delta\alpha$=$T$/8) in (31), besides replacing $\pi$ with $T$/2 in formulas (73), value $T$/8 should be subtracted from the arguments of the functions $x_1$ and $y_1$. Piecewise-linear loops with whiskers of any length can be built directly by equations (47).

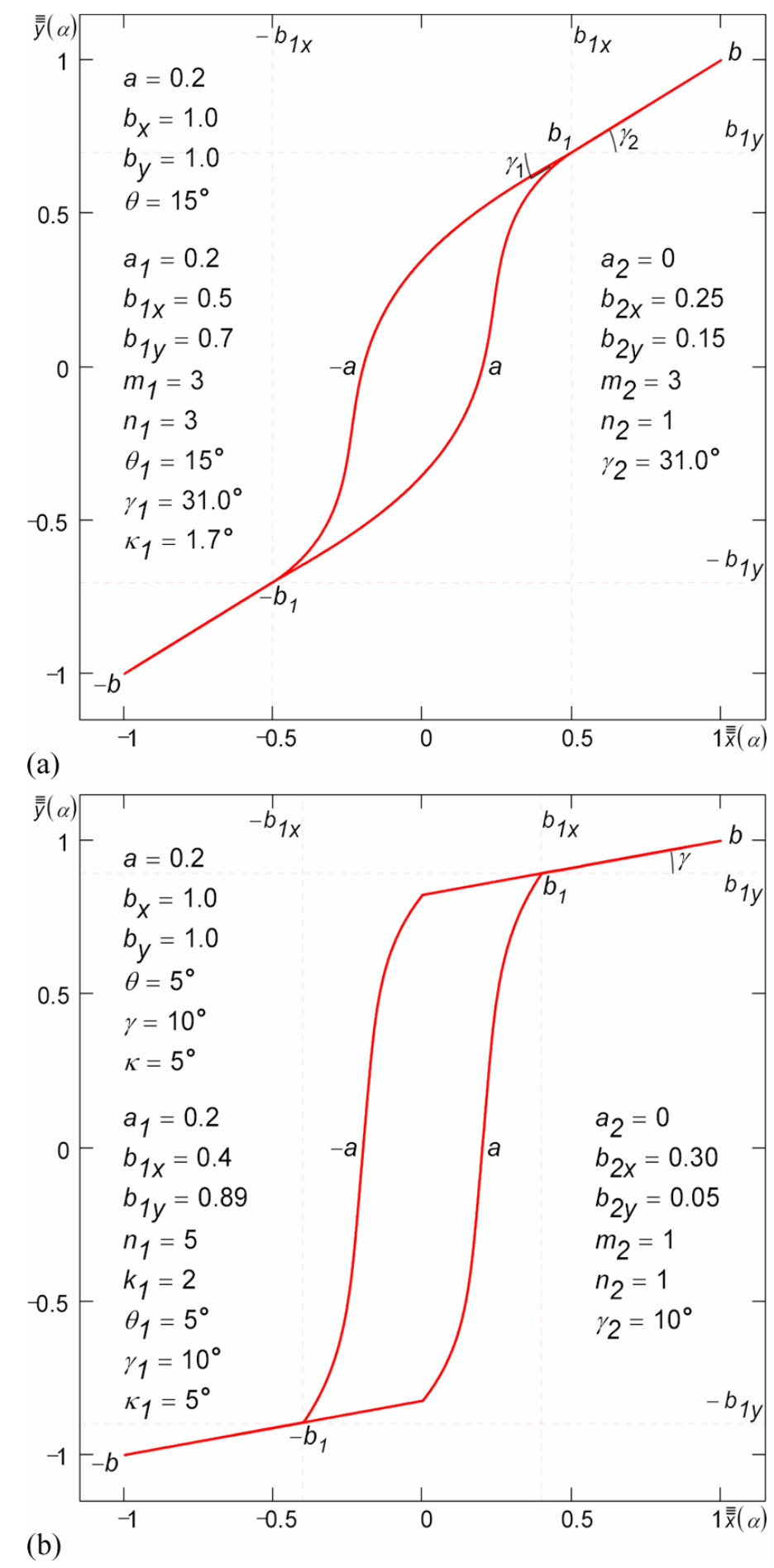

Fig. 32. Simulation of a single Classical (a) smooth, (b) hybrid loop with long whiskers by means of a triple non-self-crossing hysteresis loop. Whiskers are the outside pair of (a) smooth, (b) piecewise-linear unsplit loops of Leaf type.

### *3. Adding vertical splitting*

The simplest way to get triple self-crossing loops of the Classical type is based on representation (2), in which apart from the splitting $a_x$ along $x$, the additional splitting $a_y$ along $y$ is performed. The splitting curve equations in this case are as follows: $x_2(\alpha)$=$a_x\cos^m\alpha$, $y_2(\alpha)$=$a_y\cos^m\alpha$. Taking into account the additional splitting, the equations of the triple loop are written as

$$\begin{aligned} x(\alpha) &= a_x\cos^m\alpha + \hat{b}_x\sin^n\alpha, \\ y(\alpha) &= a_y\cos^m\alpha + \hat{b}_y\sin\alpha, \end{aligned} \tag{76}$$

where parameters $a_x$, $a_y$, $\hat{b}_x$, $\hat{b}_y$ are determined according to the formulas ($m$=1):[6]

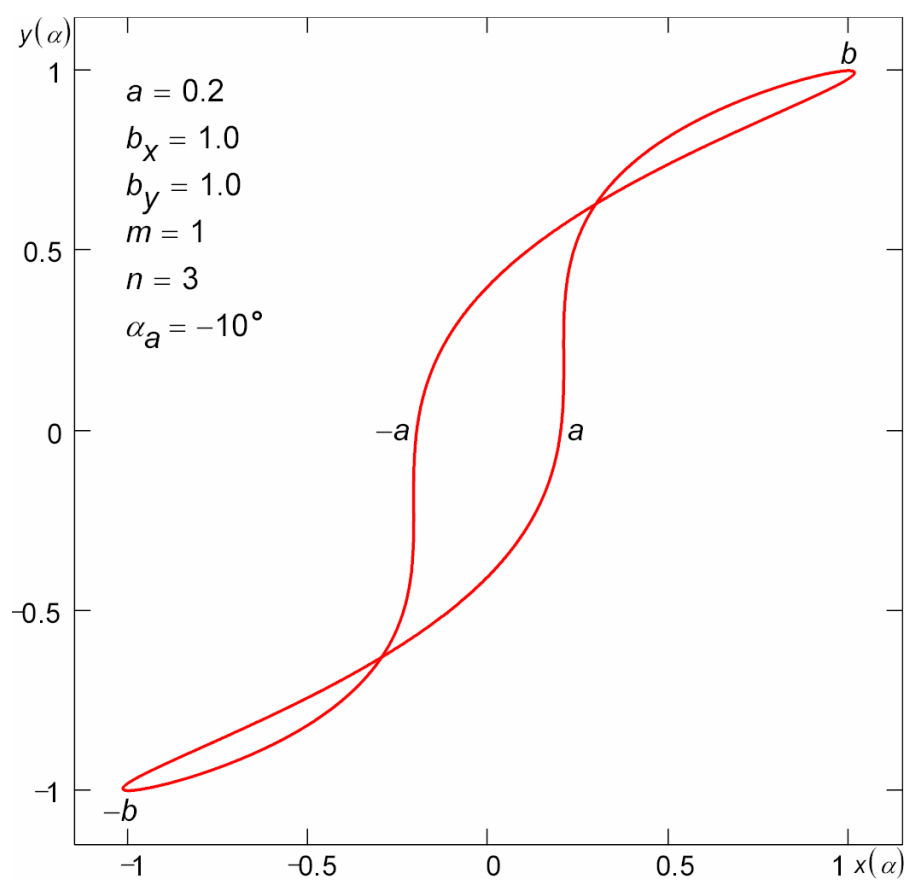

Fig. 33. Triple smooth self-crossing hysteresis loop of the Classical type formed as a result of an additional vertical splitting.

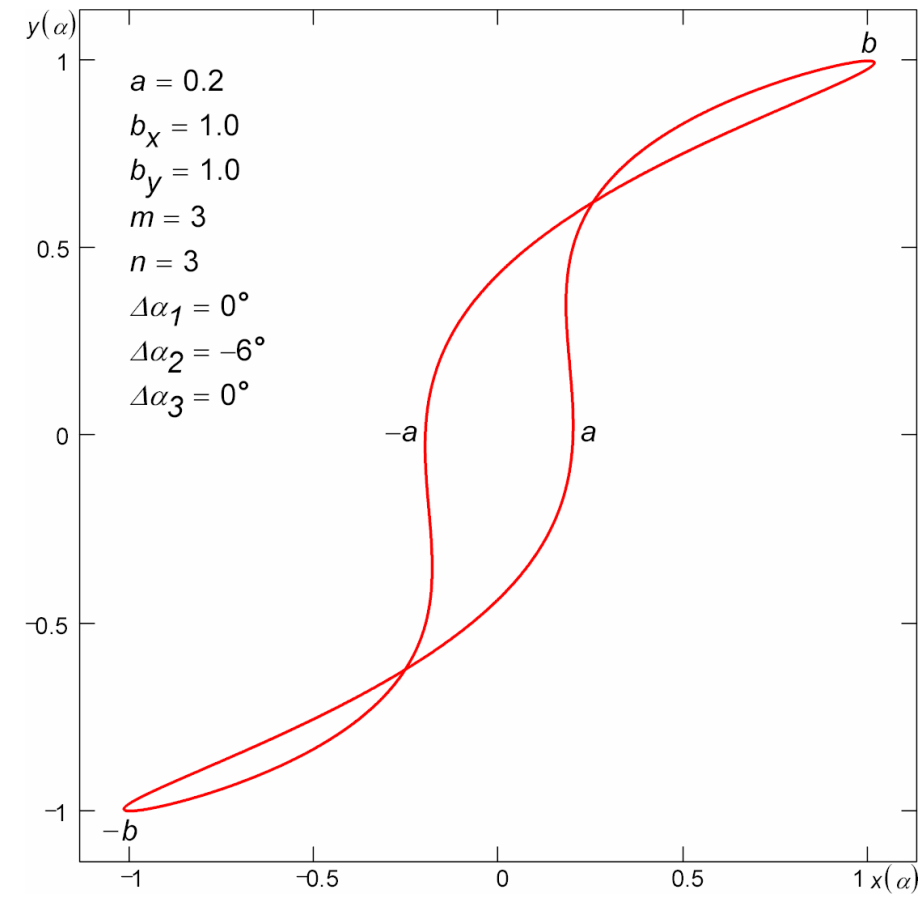

Fig. 34. Triple smooth self-crossing hysteresis loop of the Classical type formed as a result of "squeezing" by the phase shift $\Delta\alpha_2$.

$$
\begin{aligned}
a_x &= \frac{a\cos^n\alpha_a - b_x\sin^n\alpha_a}{\sin^{n+1}\alpha_a + \cos^{n+1}\alpha_a},\\
a_y &= -b_y\sin\alpha_a,\\
\hat{b}_x &= \frac{a\sin\alpha_a + b_x\cos\alpha_a}{\sin^{n+1}\alpha_a + \cos^{n+1}\alpha_a},\\
\hat{b}_y &= b_y\cos\alpha_a,
\end{aligned}
\tag{77}
$$

where $\alpha_a$ is the value of parameter $\alpha$ at the split point $a$. Negative $\alpha_a$ sets a value of foldover of the triple loop. The working expressions for the case $m$=3, 5, … are given in the supplementary material. As it follows from formulas (77), the triple loop degenerates into a regular Classical loop with $\alpha_a$=0. By setting $\alpha_a$>0, it is possible to change the shape and the curvature of a single loop.[6]

Fig. 33 shows a triple loop of the Classical type[25, 29] obtained according to the described method. It is worth noting that, assuming splitting $a$=0, the triple self-crossing loop degenerates into a double non-self-crossing loop of the Propeller type (see also Sec. II.C.3).[6]

***4. Squeeze causing a foldover***

According to the improved model (13), the triple self-crossing loops are formed by setting up a negative phase shift $\Delta\alpha_2$ (or positive $\Delta\alpha_3$) that "squeezes" the loop so tight that a foldover appears. Fig. 34 shows an example of a triple loop of the Classical type[25, 29] obtained by the described method. It is worth noting that, assuming splitting $a$=0, the triple self-crossing loop degenerates into a double non-self-crossing loop of the Propeller type (see Sec. II.C.4).

***5. Triple piecewise-linear loop as a special case of the piecewise-linear shifted loop***

A triple self-crossing loop (see Fig. 35) can be obtained from the shifted loop (see Sec. II.B.3). To do this, the splitting $a$ must be negative and the horizontal size $w$ of the shifted section of the loop must exceed the value $w$ defined by formula (61). Triple self-crossing loops Play, Relay with Gain, Relay as well as variations of these loops with attenuation and/or without whiskers are given in the supplementary material. Due to the way of formation of this triple loop, it follows that detection of such loop during physical measurements may serve as a sign of simultaneous action of two separate hysteresis processes superimposed on each other in the system under consideration.

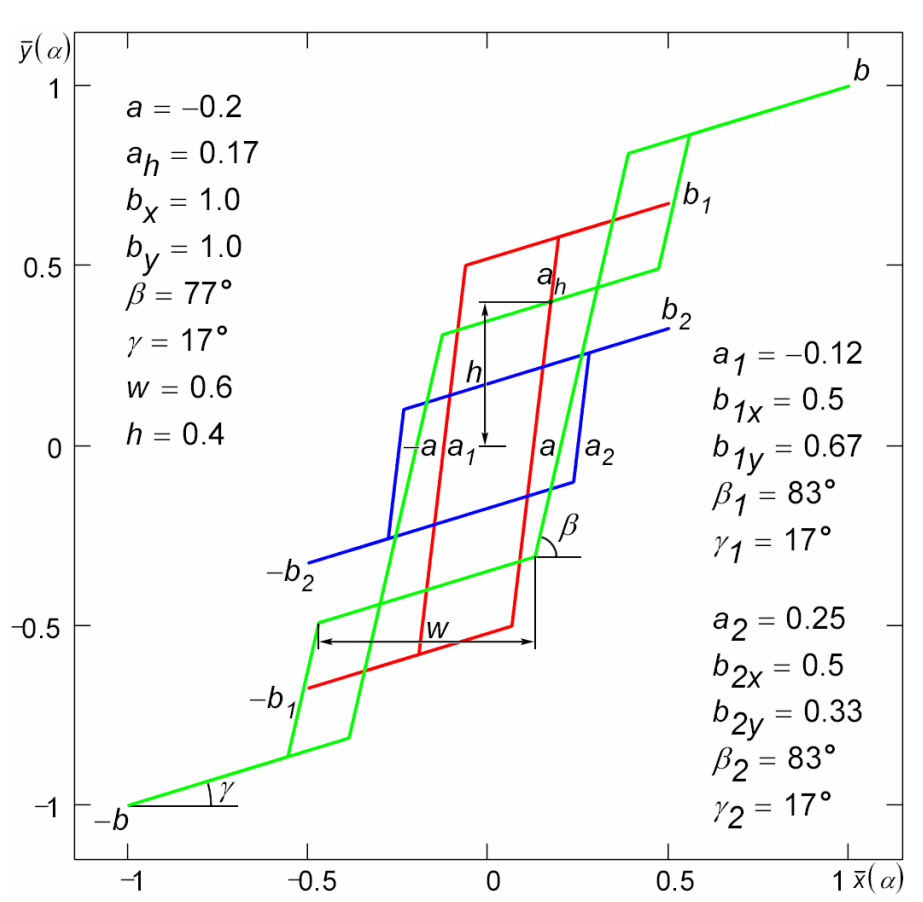

Fig. 35. Triple piecewise-linear self-crossing hysteresis loop of type Play with Gain obtained as a result of addition of two single loops Play with Gain.

## E. Area of a hysteresis loop

The area of a hysteresis loop characterizes energy losses in a piezoelectric/ferromagnetic material while applying an alternating electric/magnetic field. To find the hysteresis loop area, the following general formula is used

$$S=\oint x(\alpha)\frac{dy(\alpha)}{d\alpha}d\alpha=-\oint y(\alpha)\frac{dx(\alpha)}{d\alpha}d\alpha=\frac{1}{2}\oint\left[x(\alpha)\frac{dy(\alpha)}{d\alpha}-y(\alpha)\frac{dx(\alpha)}{d\alpha}\right]d\alpha. \qquad (78)$$

According to (78), as a hysteresis loop is being scaled as $\bar{x}(\alpha)=Ax(\alpha)$, $\bar{y}(\alpha)=By(\alpha)$, its area changes proportionally to the product of the scale factors $A$ and $B$, i. e., $\bar{S}=ABS$. A series of three hysteresis loops, where the areas are doubling with scaling, is shown in Fig. 36(a).

***1. Smooth loops***

By substituting the expressions $x(\alpha)$, $y(\alpha)$ of the improved model (13) along with their derivatives into (78) and integrating, one can obtain[6, 30]

$$\begin{aligned}S&=\left[\frac{\hat{a}}{2^{m-1}}C_m^{\frac{m-1}{2}}\cos(\Delta\alpha_1-\Delta\alpha_3)+\frac{\hat{b}_x}{2^{n-1}}C_n^{\frac{n-1}{2}}\sin(\Delta\alpha_2-\Delta\alpha_3)\right]\pi b_y\\&=\frac{Am_1(\cos\Delta\alpha_3-\tan\varphi_1\sin\Delta\alpha_3)}{\sqrt{\tan^2\varphi_1+1}}\pi b_y=(A_1\cos\Delta\alpha_3-B_1\sin\Delta\alpha_3)\pi b_y,\end{aligned} \qquad (79)$$

where $n$ is an odd number; $A_1$, $B_1$ are the Fourier coefficients for the 1st harmonic (see Eqs. (25)). It follows from formula (79) that the area of the loop (13) is defined by the amplitude $Am_1$ and the phase $\varphi_1$ of the 1st harmonic only; the other harmonics of the expansion (7) of the generating function $x(\alpha)$ do not affect the area $S$ in any way.

The loop tilt produced by a phase shift $\Delta\alpha_1$ ($m$>1, $\Delta\alpha_2$=$\Delta\alpha_3$=0) results in an increase in the loop area $S$ (see Fig. 3). The area $S$ of a hysteresis loop increases with an increase in the phase shift $\Delta\alpha_2$ (see Fig. 4). With an increase in the phase shift $\Delta\alpha_3$, the area $S$ decreases (see Fig. 5).

It follows from formula (79) that with $\Delta\alpha_2$=$\Delta\alpha_3$, loop area

$$S=\frac{C_m^{\frac{m-1}{2}}\pi a b_y}{2^{m-1}\cos^{m-1}(\Delta\alpha_1-\Delta\alpha_3)}. \qquad (80)$$

Thus, under the specified condition the loops will have the same area $S$ regardless of their type ($n$ is absent in the formula) and saturation $b_x$ provided that the parameters $a$, $b_y$ and $m$ are the same (see Fig. 36(b)).

It follows from formula (80) that in case of $\Delta\alpha_2$=$\Delta\alpha_3$=0, the area of a loop having a negative slope ($\Delta\alpha_1$>0) at the split point equals to the area of a loop having a positive slope ($\Delta\alpha_1$<0) at the split point, i. e., $S\big|_{+\Delta\alpha_1}=S\big|_{-\Delta\alpha_1}$. The inequality $S\big|_{|\Delta\alpha_1-\Delta\alpha_3|>0}>S\big|_{\Delta\alpha_1-\Delta\alpha_3=0}$ also follows from formula (80). Thus, for any loop type in case of $\Delta\alpha_2$=$\Delta\alpha_3$=0 ($m$>1), the energy losses in the tilted loops ($\Delta\alpha_1\neq0$, see Fig. 3) are greater than in the "upright" ones ($\Delta\alpha_1$=0).

It also follows from formula (80) that with $m$=1 ($\Delta\alpha_2$=$\Delta\alpha_3$), loop area

$$S = \pi a b_y. \tag{81}$$

Thus, under the specified conditions the loops will have the same area regardless of their type ($n$ is absent in the formula), saturation $b_x$ and phase shifts provided that the parameters $a$ and $b_y$ are the same. Fig. 36(c) shows an example of loops of the same area having different phase shifts $\Delta\alpha_1$. The formula of areas of the Leaf ($n$=1) loops with $m$=1 and $\Delta\alpha_1$=$\Delta\alpha_3$ looks the same as (81); the areas of these loops are also independent of the values of saturation $b_x$ and the values of phase shifts. Area (81) is numerically equal to the area of an ellipse having semi-major axis $a$ and semi-minor axis $b_y$.

To calculate the area of the loop based on the model (1), the zero phase shifts $\Delta\alpha_1$=$\Delta\alpha_2$=$\Delta\alpha_3$=0 should be substituted into formula (79). As a result, we can obtain

$$S = \frac{1}{2^{m-1}} C_m^{\frac{m-1}{2}} \pi a b_y = \frac{\pi A m_1 b_y}{\sqrt{\tan^2 \varphi_1 + 1}} = \pi A_1 b_y. \tag{82}$$

There is an inaccuracy in Ref. 1 in formula (27) and in the accompanying text to this formula. Formula (27) is valid for any $n$ (even or odd). Since formula (27) does not have a solution with $m$=1, formula (82) presented in this paper should be used instead of it. By the same reason instead of formula (28) from Ref. 1, the corresponding formula presented in the supplementary material should be used.

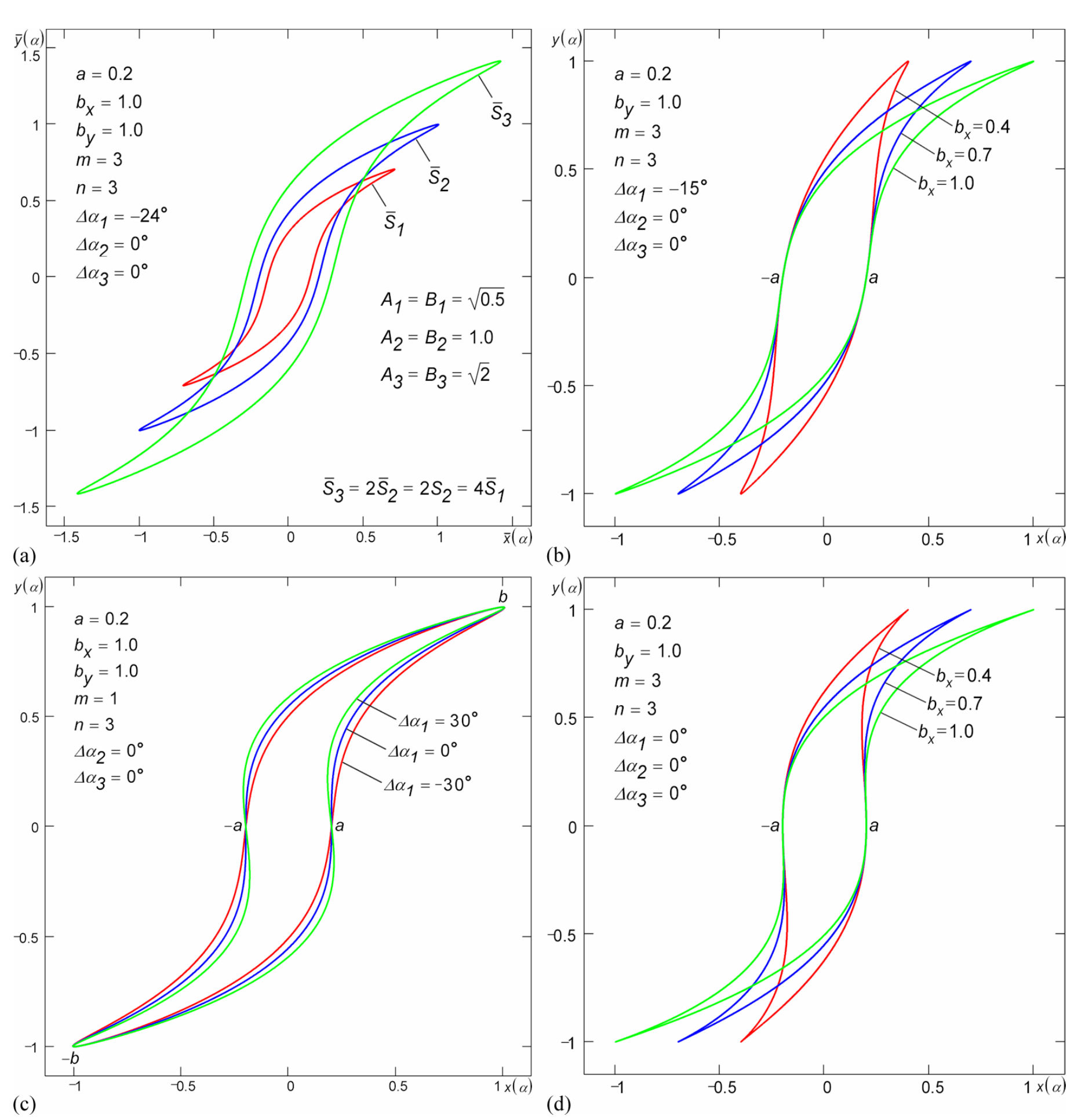

Fig. 36. The hysteresis loop area (a) is proportional to the product of the scale factors $A$ and $B$ along axes $x$ and $y$, respectively, (b) is invariant to the saturations $b_x$ (restrictive condition is: $\Delta\alpha_2$=$\Delta\alpha_3$), (c) is invariant to the phase shifts $\Delta\alpha_1$ (restrictive conditions are: $m$=1, $\Delta\alpha_2$=$\Delta\alpha_3$), (d) is invariant to the saturations $b_x$.

Since $n$ does not enter into formula (82), all three loop types Leaf, Crescent, and Classical of the model (1) have the same loop area $S$ provided that the other parameters $a$, $b_y$, and $m$ of these loops are the same (see Fig. 1). Since parameter $b_x$ does not enter into formula (82) as well, all the loops (1) having the same parameters $a$, $b_y$, and $m$ but different $b_x$ have the same area $S$ (see Fig. 36(d)).

The area of smooth loops (29) tilted at the split point by angle $\theta$ by skewing the coordinate system can also be calculated by formula (82).[6] Since the oblique angle $\theta$ is not included in formula (82), all loops tilted by skewing at any angle $\theta$ have the same area $S$ provided that the parameters $a$, $b_y$, and $m$ of these loops are the same (see Fig. 8).

Area of the loop Classical (30) tilted and curved by skewing is calculated by the general formula[6]

$$S=\frac{\pi a}{2^{m-1}}C_m^{\frac{m-1}{2}}\left\{b_x\tan\kappa\left[1-(m+1)\prod_{k=0}^{\frac{n-1}{2}}\frac{2k+1}{m+n-2k}\right]+b_y\right\}, \qquad (83)$$

where $n$ is an odd number. According to (83), the area does not depend on the oblique angle $\theta$, therefore, loops tilted by skewing by any angle $\theta$ have the same areas provided that all other loop parameters are the same (see Fig. 9). For $\kappa$=0 (see Fig. 8), formula (83) transforms into (82). In particular case, for example, when $m$=$n$=3 (see Fig. 9), the formula (83) is as follows

$$S=\frac{3}{8}\pi a\left(b_x\tan\kappa+2b_y\right). \qquad (84)$$

It is noteworthy that with $\Delta\alpha_1=\Delta\alpha_2=\Delta\alpha_3=\Delta\alpha\neq 0$, where $\Delta\alpha$ is an arbitrary real number, formula (79) also takes the form (82). In this case, loops built according to models (1) and (13) will have equal areas. This is because under the above conditions, these loops are of exactly the same shape (see Sec. II.A.2.d).

Area of the loop Classical (14) tilted by rotation is calculated by the general formula[6]

$$S=\frac{\pi a}{2^{m-1}}C_m^{\frac{m-1}{2}}\left\{\sin\theta\left(b_x\cos\theta-b_y\sin\theta\right)\left[1-(m+1)\prod_{k=0}^{\frac{n-1}{2}}\frac{2k+1}{m+n-2k}\right]+b_y\right\}, \qquad (85)$$

where $n$ is an odd number. According to (85), the classical loop areas are different for $\theta$ rotation angles having opposite signs. For $\theta$=0, formula (85) transforms into (82). In particular case, for example, when $m$=$n$=3, the formula (85) is as follows

$$S=\frac{3}{8}\pi a\left[\sin\theta\left(b_x\cos\theta-b_y\sin\theta\right)+2b_y\right]. \qquad (86)$$

***2. Piecewise-linear and hybrid loops***

The areas of the simplest piecewise-linear loops having a parallelogram shape (see Fig. 15, Fig. 19, Fig. 20, Fig. 21) are easily determined using formula of the area of a parallelogram ($m$=$n$=1)[6]

$$S=4a\left[\frac{(a-b_x)\tan\gamma+b_y}{\tan\beta-\tan\gamma}\tan\beta-\frac{a\sin\beta\sin\gamma}{\sin(\beta-\gamma)}\right]. \qquad (87)$$

With no gain ($\gamma$=0), formula (87) becomes much simpler

$$S=4ab_y. \qquad (88)$$

Since area (88) corresponds to the area of parallelogram with the side $2a$ and the height $2b_y$, the areas of piecewise-linear loops of Play and Non-ideal Relay types (see Figs. 20(e), (f), (k), (l)) can be determined by this formula. Since the hybrid loops (37) with zero phase shifts (see Fig. 11, Fig. 17) have translational symmetry, their area is equal to the area of parallelogram with side $2a$ and height $2b_y$ and thus, it can also be calculated by formula (88).

Area of the hybrid Classical loop with gain/attenuation ($\gamma\neq 0$) (38) (see Fig. 18) is calculated by the following general formula[6]

$$S=4a\left\{\left[(b_x-a)\tan\kappa-b_x\tan\gamma+b_y\right]\left[\frac{k(n-1)}{(k+1)(n+k)}\tan\theta\tan\gamma+1\right]+(b_x-a)\left(\frac{k}{n+k}\tan\gamma-\tan\kappa\right)\right\}, \qquad (89)$$

where $n$ is an odd number, $k$ is an even number. With no gain ($\gamma$=0), formula (89) takes the simple form (88). As follows from (88), the zero gain hybrid loop areas do not depend on $b_x$, $n$, $\theta$, or $\kappa$.

Table. Average relative approximation error $\langle\delta\rangle$ (%) of the original and the improved parametric models of hysteresis loop.

| Model | Loop Leaf, Fig. 37(a) | Loop Classical, Fig. 37(b) | Loop Classical, Fig. 37(c) |
|---|---|---|---|
| Original | 0.8 | 2.9 | 1.7 |
| Improved | 0.5 | 1.0 | 1.0 |

Since the $\alpha$ parameter values at all corner points of the piecewise-linear loops are known, their Cartesian coordinates are known also. Therefore, the area of the polygonal loops, such as Play-Relay-Play, Play-Play, Play-Relay (see. Figs. 12-14, Fig. 16, Fig. 23); Play-Play-Play, Play-Relay-Play (see Fig. 22); Shifted Play, Shifted Non-ideal Relay (see Fig. 24) and similar, can be calculated using the formula for the area of a polygon specified by coordinates of its vertices.

In a number of particular cases, the polygonal loops degenerate into such polygons, which can be represented as a set of trapeziums and/or parallelograms/triangles. For example, the loop Play-Relay-Play with Gain in Fig. 22(c) can be represented as two identical trapeziums and a rectangle; and the loop Shifted Play with Gain in Fig. 24(a) – as three parallelograms, two of which are identical. Thus, to calculate areas of these and similar loops, the corresponding formulas of the areas of trapeziums, parallelograms or triangles can be used.

## III. APPLICATION OF THE IMPROVED MODEL

The developed model enables building smooth, piecewise-linear, hybrid, minor, mirror-reflected, inverse, reverse, double and triple loops. Calculation of derivatives, searching for the harmonically linearized transfer function of a hysteresis element and the inverse function are performed in the improved model similarly to Ref. 1.

The usage of phase shifts $\Delta\alpha_1$, $\Delta\alpha_2$, $\Delta\alpha_3$ and a number of other transformations which are changing the loop tilt and curvature, let us reduce the hysteresis loop approximation error by several times. Error analysis conducted according to the method described in Ref. 1 has shown that the average relative approximation error

$$\langle\delta\rangle = \frac{100\%}{2nb_y}\sum_{i=1}^{n}\left|y_m(x_i) - y_e(x_i)\right| \qquad (90)$$

(where $y_m(x)$ are the model data, $y_e(x)$ are the experimental data, $n$ is the number of points on the ascending or descending section of an experimental loop) of the improved model does not exceed 1% (see Table).[6] For the comparison to be correct, the experimental hysteresis loops from Ref. 1 were used during the error determination. Fig. 37 shows approximating loops built according to the existing (1) and the improved (13) models overlaid with the experimental loops from Ref. 1. It is easy to see that the approximation accuracy of the improved model is notably higher than the approximation accuracy of the existing model. To minimize error, it is recommended to draw the approximating loop by using least squares.[31]

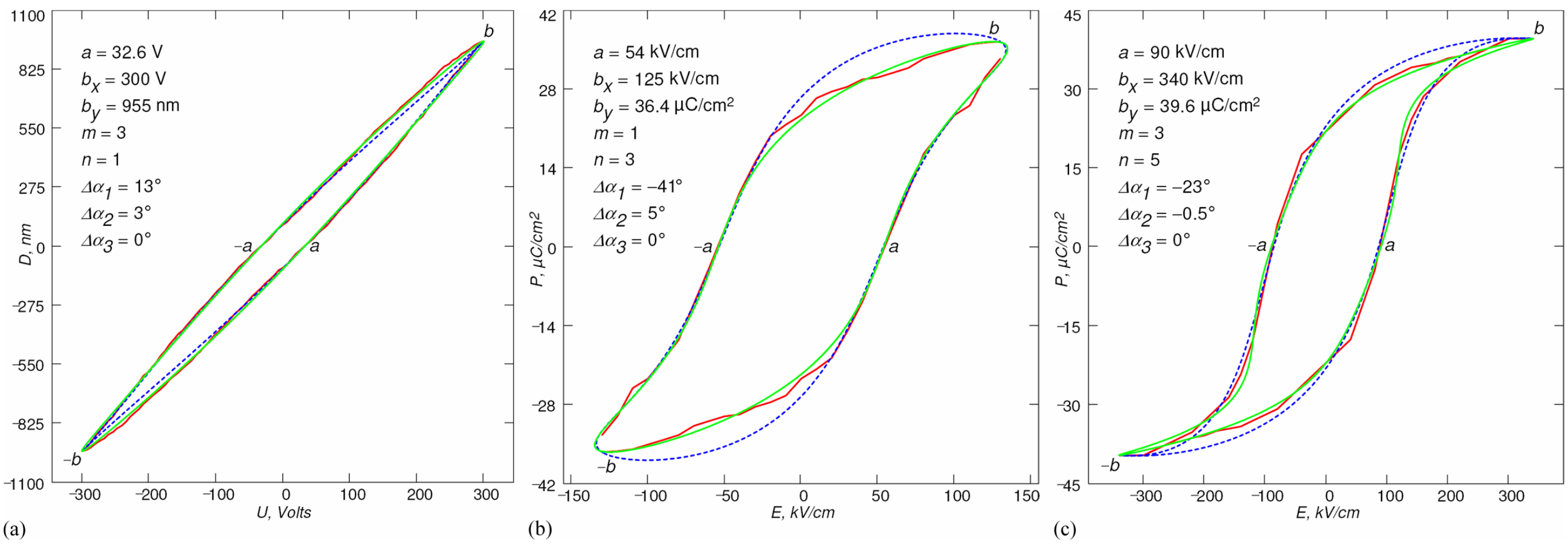

Fig. 37. Approximation of real smooth hysteresis loops of types (a) Leaf, (b), (c) Classical. The red uneven loop is experimental; the blue dotted loop is the existing model; the green smooth one is the suggested improved model. Approximation error $\langle\delta\rangle$ of the improved model does not exceed 1%.

Like the original model, the improved one can be applied for correction of the distortions caused by hysteresis of piezomanipulators (see Fig. 37(a)) of the scanning probe microscope (SPM). To do this, phase-shifting elements providing the phase shifts $\Delta\alpha_1$, $\Delta\alpha_2$ (or $\Delta\alpha_1$, $\Delta\alpha_3$) should be added to the existing hardware suggested in Ref. 1.

Instead of calculating the $m$th and $n$th powers of sine signal using the hardware multipliers, as it was suggested in Ref. 1, it is possible to simply sum the sine and cosine signals of multiple frequencies according to (3). For example, to compensate for the distortions caused by the hysteresis loop Leaf ($m$=3, $n$=1) shown in Fig. 37(a), two cosine signals of frequencies $\omega$ and $3\omega$, an element shifting the phase by a quarter of period, three amplifiers with gains $3/4a$, $a/4$, $b_x$ and two summing amplifiers are required.[6]

According to (7), in order to compensate distortions of the loop Leaf, two cosine signals of frequencies $\omega$ and $3\omega$, an element shifting the phase by value $\varphi_1=\arctan(4b_x/3a)$, two amplifiers with gains $Am_1=\sqrt{9a^2/16+b_x^2}$, $Am_3=a/4$ and a summing amplifier are required.[6]

## IV. SUMMARY

Improvements have been made to the existing model of a hysteresis loop built on parametric equations; several formulas derived earlier were refined. As a result of the improvements, the approximation accuracy has increased several times. It is shown that the generating function of the hysteresis loop can be represented either as a sum of an unsplit loop and a splitting curve, or in the form of a frequency spectrum, or in an exponential form.

A general formula has been derived for building piecewise-linear hysteresis loops of Play and Non-ideal Relay types as well as their numerous variations widely used in simplified models of hysteresis phenomena. The possibility of composing/decomposing various hysteresis loops is demonstrated. For a number of loop types, the results of the composition/decomposition suggest the simultaneous existence in the system under consideration not one but two different hysteresis processes superimposed on each other. Several new formulas describing various types of single, double, and triple loops have been obtained.

In addition to the presented parametric equations of the parallelogram and the rectangle, the parametric equations of a rhombus, square, regular hexagon and regular octagon are found.[6] A more general expression has been suggested to determine the hysteresis loop area, which allows for evaluation of energy losses in piezoelectric/ferromagnetic materials. In the course of the study, several identities related to the binomial coefficients were found.[6]

The hysteresis loop model developed is especially suitable for solving the tasks of simulation of cyclically operating instruments that include hysteresis elements.[32] In addition, the model permits to build the output signals $x(t)$ and $y(t)$ with simple hardware components thus allowing easy hardware implementation of both direct and inverse hysteresis loops.

## SUPPLEMENTARY MATERIAL

Supplementary material includes zip-archive with Mathcad 2001i worksheets, where all aspects of the original and the improved parametric models of hysteresis loop are considered in detail. Those who do not have Mathcad software may take advantage of the enclosed readable Mathcad worksheets as a PDF-document. Due to restriction on the length of the article, it presents only the most common hysteresis loops. If the required loop is absent in the article, it makes sense to search in the supplementary material.

## DATA AVAILABILITY

The data that supports the findings of this study are available within the article and its supplementary material.

## ACKNOWLEDGMENTS

I thank Oleg E. Lyapin, David W. Waddell, Oleg V. Obyedkov and Larisa B. Sharova for critical reading of the manuscript and checking of the supplementary materials; Assoc. Prof. Eugene A. Fetisov for support and stimulation.